\documentclass[journal,onecolumn]{IEEEtran} 

\usepackage{cite}

\ifCLASSINFOpdf
 \usepackage[pdftex]{graphicx}
  \graphicspath{ {../Figures/}, {Figures/} }
\else
  \usepackage[dvips]{graphicx}
  \graphicspath{{../Figures/}}
\fi

\usepackage{amsmath}
\usepackage{algorithmic}

\usepackage{array}

\ifCLASSOPTIONcompsoc
  \usepackage[caption=false,font=normalsize,labelfont=sf,textfont=sf]{subfig}
\else
  \usepackage[caption=false,font=footnotesize]{subfig}
\fi

\usepackage{fixltx2e}

\usepackage{url}

\usepackage{mathtools,amssymb,amsfonts,dsfont,relsize,exscale,commath} 
\usepackage{graphicx,epstopdf}
\usepackage[margin=2cm]{geometry}
\usepackage{adjustbox,multirow,cleveref}
\usepackage[dvipsnames]{xcolor}
\usepackage{bm} 
\usepackage{tikz}

\DeclareMathOperator{\sign}{sgn}

\DeclarePairedDelimiter\floor{\lfloor}{\rfloor}

\newcommand{\defeq}{\vcentcolon=}

\renewcommand{\arraystretch}{1.3}

\def\levy{{L\'{e}vy}}

\def\w{{w}} 

\def\c{{c}} 

\def\barg{\bar{g}(1)}

\newtheorem{theorem}{Theorem}
\newtheorem{remark}{Remark}
\newtheorem{lemma}{Lemma}

\usepackage{color}

\newcommand{\RL}{{\mathbb R}}

\def\ba{\begin{align}}
\def\ea{\end{align}}
\def\ban{\begin{align*}}
\def\ean{\end{align*}}

\def\be{\begin{eqnarray}}
\def\ee{\end{eqnarray}}
\def\ben{\begin{eqnarray*}}
\def\een{\end{eqnarray*}}

\def\bqq{\begin{equation}}
\def\eqq{\end{equation}}
\def\bqqn{\begin{equation*}}
\def\eqqn{\end{equation*}}

\def\elabel#1{\label{e:#1}}

\def\sq{$\Box$}

\def\qed{\ifmmode\sq\else{\unskip\nobreak\hfil
\penalty50\hskip1em\null\nobreak\hfil\sq
\parfillskip=0pt\finalhyphendemerits=0\endgraf}\fi\par\medbreak}

\def\sign{{\rm sign}}

\newsavebox{\junk}
\savebox{\junk}[1.6mm]{\hbox{$|\!|\!|$}}

\def\limsup{\mathop{\rm lim\ sup}}

\def\til={{\widetilde =}}

 \def\eq#1/{(\ref{#1})}

\newtheorem{corollary}[theorem]{Corollary}
\newtheorem{proposition}[theorem]{Proposition}

\def\eq#1/{(\ref{e:#1})}

\newcommand{\beqn}[1]{\notes{#1}%
\begin{eqnarray} \elabel{#1}}

\newcommand{\eeqn}{\end{eqnarray} } 

\newcommand{\beq}[1]{\notes{#1}%
\begin{equation}\elabel{#1}}

\newcommand{\eeq}{\end{equation}} 

\def\bdes{\begin{description}}
\def\edes{\end{description}}

\newcommand{\oo}{\overline}

\def\barg{{\oo {g}}}

\def\notes#1{}

\definecolor{mag}{rgb}{0.7,0,0.3}
\definecolor{dgreen}{rgb}{0.1,0.5,0.1}
\definecolor{dred}{rgb}{.8,0,0}
\definecolor{gray}{rgb}{.8,.8,.8}
\definecolor{brown}{rgb}{0.6451,0.3706,0.1745}

\begin{document}

\title{Nonasymptotic Gaussian Approximation\\
for Inference with Stable Noise}

\author{Marina~Riabiz,~\IEEEmembership{Student Member,~IEEE,}
        Tohid~Ardeshiri,~\IEEEmembership{Member,~IEEE,}     
        Ioannis~Kontoyiannis,~\IEEEmembership{Fellow,~IEEE,}    
        and~Simon~Godsill,~\IEEEmembership{Member,~IEEE}


\thanks{M. Riabiz, T. Ardeshiri, S. Godsill  and I. Kontoyiannis  are with the Department of Engineering,
	University of Cambridge, Cambridge, CB2 1PZ, UK. 
	Email addresses: mr622@cam.ac.uk, ta417@cam.ac.uk, 
	ik355@cam.ac.uk, sjg30@cam.ac.uk.
	Preliminary versions of some of these results
	were presented in the conference 
	papers~\cite{Riabiz2017,Riabiz2017b,Riabiz2018}. This work was in part supported through the  EPSRC-DTG2013/EU scheme.}%
}

%
%
%

\maketitle

\begin{abstract}
The results of a series of theoretical studies are reported, examining the convergence rate for different approximate representations of $\alpha$-stable distributions. Although they play a key role in modelling random processes with jumps and discontinuities, the use of $\alpha$-stable distributions in inference often leads to analytically intractable problems. The LePage series, which is a probabilistic representation employed in this work, is used to transform an intractable, infinite-dimensional inference problem into a conditionally Gaussian parametric problem. A major component of our approach is the approximation of the tail of this series by a Gaussian random variable. Standard statistical techniques, such as Expectation-Maximization, Markov chain Monte Carlo, and Particle Filtering, can then be applied. In addition to the asymptotic normality of the tail of this series, we establish explicit, nonasymptotic bounds on the approximation error. Their proofs follow classical Fourier-analytic arguments, using Ess\'{e}en's smoothing lemma. Specifically, we consider the distance between the distributions of: $(i)$~the tail of the series and an appropriate Gaussian; $(ii)$~the full series and the truncated series; and $(iii)$~the full series and the truncated series with an added Gaussian term. In all three cases, sharp bounds are established, and the theoretical results are compared with the actual distances (computed numerically) in specific examples of symmetric $\alpha$-stable distributions. This analysis facilitates the selection of appropriate truncations in practice and offers theoretical guarantees for the accuracy of resulting estimates. One of the main conclusions obtained is that, for the purposes of inference, the use of a truncated series together with an approximately Gaussian error term has superior statistical properties and is likely a preferable choice in practice.
\end{abstract}


%
\IEEEpeerreviewmaketitle

\section{Introduction}

\IEEEPARstart{S}{tatistical} modelling and inference 
for time series and random processes are 
of central importance in many areas 
of science and engineering. In applications, the
time- or space-evolution of quantities of interest
is often described through regression models
that include random `noise' components. These components
may represent the inherent randomness in 
the underlying system, 
or the noise introduced by the observation process,
or both.


Consider, for example, a simple discrete-time linear regression model
for a time series $\mathbf{x} := [x_1, \ldots, x_N]'$, 
expressed as, 
\begin{IEEEeqnarray}{rCl} \label{eq:AR}
	{\bf x}={\bf G}{\bm\lambda}+{\bf u},
\end{IEEEeqnarray}
where the $P$-dimensional parameter vector
$\bm\lambda := [\lambda_1, \ldots, \lambda_P]'$
and the $N\times P$ matrix of known regressors
$\bf G$ describe the
deterministic part of the system,
and the random process $\{u_n\}$ describes
the random noise component,
$\mathbf{u} := [u_1, \ldots, u_N]'$.
This encompasses many models of current interest,
including
Fourier, wavelet and other expansions used in 
compressive sensing, communication systems, genomics,
and signal processing.  

Another common class of motivating examples
is that of  state-space models where the state evolves 
over time with random disturbances, as,
\begin{IEEEeqnarray}{rCl}
x_n=Ax_{n-1}+u_n, \qquad
1\leq n\leq N,
\label{state_space}
\end{IEEEeqnarray}
\noindent 
where $A$ is the autoregressive parameter, 
and observations may in addition contain noise: 
\begin{IEEEeqnarray}{rCl} \label{eq:AR:noisy}
{y_n}={b x_n}+{v_n}, \qquad 1\leq n\leq N.
\end{IEEEeqnarray}
Here $b$ is an observation parameter and $\{v_n\}$ is
the observation noise process.

Depending on the application at hand, there are many possible
inference objectives;
for example,
state inference or prediction for $x_n$ in the state-space model, 
parameter estimation for $\bm \lambda$, $A$ and $b$, 
and model choice to determine the structure and dimensionality of the model. 
A common modelling choice is to assume that
the processes $\{u_n\}$
and $\{v_n\}$ are Gaussian:  
Since the driving noise process $\{u_n\}$ can often be thought 
of as the sum of many small independent contributions, 
the Gaussian assumption is a natural consequence 
of the central limit theorem (CLT). 
Similarly, the measurement noise process $\{v_n\}$ typically 
is the result of the sum of small independent perturbations, 
which again justifies the Gaussian assumption via the CLT.
In these cases, standard methods are available for likelihood-based 
or Bayesian inference, using closed-form results combined with, 
for example, variational Bayes 
or Monte Carlo sampling.

However, many real-world cases exhibit extreme values 
much more frequently than the Gaussian model of \eqref{eq:AR} or \eqref{state_space} would allow. 
Examples of such abrupt changes include variations presented by stock prices 
or insurance gains/losses in financial applications, 
as studied extensively since the seminal works 
\cite{Mandelbrot1963} and \cite{Fama1965a}; {we refer to \cite{Rachev2000} for a more recent review. }
Further applications can be found in various fields of engineering, such as communications ({see \cite{Azzaoui2010} for statistical modelling of channels,  \cite{Fahs2012, FreitasEganClavierEtAl2017} for capacity bounds, }
\cite{LiebeherrBurchardCiucu2012} for delay bounds in networks 
with $\alpha$-stable noise, and \cite{ShevlyakovKim2006, WarrenThomas1991} 
for signal detection), 
signal processing \cite{Nikias1995}, image analysis \cite{Achim2001, Achim2006}
and audio processing \cite{Lombardi2006}. 
Sudden changes are studied also in the climatological sciences 
\cite{Katz1992, Katz2002}, and in the medical sciences; see,
e.g.,
\cite{ChenWangMcKeown2010} on brain connectivity representations. Moreover, in the field of sparse modelling and Compressive Sensing, a noise distribution is required that leads to sparse solutions (in transformed domains), a case much better dealt with using heavy-tailed models than the Gaussian case; see, e.g., \cite{UnserTaftiAminiEtAl2014, UnserTaftiSun2014, Unser2014, AminiUnser2014, CarrilloRamirezArceEtAl2016} for a detailed review of modelling with {sparse} signals, and a connection between sparsity and heavy-tailed distributions, 
\cite{Lopes2016, ZhouYu2017} for the estimation of the degree of sparsity, and
\cite{Tzagkarakis2009, AchimBuxtonTzagkarakisEtAl2010} for compressed-sensing Bayesian methods based on heavy-tailed assumptions.


In many of these situations, the random phenomena considered can be still thought 
of as emerging from the combination of many independent perturbations. 
According to the generalized CLT 
\cite[p.~162]{Gnedenko1968}\cite[p.~576]{Feller1966}, 
whenever the sum of independent identically distributed (i.i.d.) random 
variables (RVs) converges in distribution, it converges to 
a member of the class of $\alpha$-stable distributions;
this class is central to this paper, and it is introduced in detail 
in Section~\ref{susection:intro:stable}. The Gaussian 
is a special member of this class, the only one with finite variance. 
Hence, using non-Gaussian $\alpha$-stable distributions for the 
vector $\bf u$ in \eqref{eq:AR} offers a way of modelling time series with 
large (extreme) values. 

The main motivation for this work, as well as the main driving force
for the large attention that the $\alpha$-stable laws have received
in applications (see the extensive bibliography listed in \cite{nolan_web}),
both stem from the 
key role of the $\alpha$-stable distribution in the generalized CLT,
and from the modelling flexibility offered by the class of $\alpha$-stable
laws.

\subsection{$\alpha$-stable distributions}\label{susection:intro:stable}

We adopt the standard notation of \cite{Samoradnitsky1994}. 
We write
	$X\sim\mathcal{S}_\alpha(\sigma,\beta,\mu)$
to denote that the RV $X$ has an $\alpha$-stable distribution 
with parameters $\sigma,\beta$ and $\mu$,
where $\alpha \in (0,2)$, is the tail parameter. 
Indeed, as consequence of the generalized CLT 
\cite[Theorem~XVII.5.1]{Feller1966}, 
when $\mu=0$, the probability density function (PDF) $p(x)$
of $X$ has tails that decay like $|x|^{1+\alpha}$,
{
\begin{IEEEeqnarray*}{rCl}
	\lim_{|x| \rightarrow \infty}\frac{p(x)}{|x|^{-1-\alpha}}= C(\alpha, \sigma, \beta),
\end{IEEEeqnarray*}
}
for some finite constant $C(\alpha, \sigma, \beta)$.
This asymptotic behaviour of the PDF corresponds to the presence of extreme values in the distribution, with more extreme  values (and hence heavier tails) appearing  more frequently for smaller values of~$\alpha$. 
The parameter $\beta\in[-1,1]$
is a measure of skewness: $\beta=0$ corresponds to symmetric stable laws, 
while $\beta = \pm 1$ 
corresponds to the fully left or right skewed cases.
Finally, 
$\mu\in(-\infty,\infty)$ and $\sigma >0$ are the location and scale parameters,
respectively. 


The characteristic function~(CF) 
$\phi_X(s):=\mathbb{E}\left[ \exp{\left( isX\right) }\right]$,
for $s\in\mathbb{R}$,           
of an $\alpha$-stable RV 
$X\sim\mathcal{S}_\alpha(\sigma,\beta,\mu)$ 
can be expressed 
\cite{Gnedenko1968}
as,
\begin{IEEEeqnarray}{ll} 
	\log (\phi_X(s)) 
	=&  
	\begin{cases} 
		-\sigma^{\alpha}{|s|}^{\alpha}\left\lbrace 1-i\beta\sign(s)\tan\frac{\pi\alpha}{2}\right\rbrace + i\mu s, 
&		\alpha \neq 1, 
		\\ 
		-\sigma|s|\left\lbrace 1+i\beta\sign(s)\frac{2}{\pi}\log|s|\right\rbrace + i\mu s,   
		&
\alpha = 1.
	\end{cases} 
	\label{eq:cf}
	\IEEEeqnarraynumspace
\end{IEEEeqnarray}
Throughout, $\log$ denotes the natural logarithm.
Notice that this CF has a pole for $\alpha =1$.
For the sake of simplicity, throughout the paper we assume
that $\alpha\neq 1$.
Also, it is easy to see from~\eqref{eq:cf} that the
class of $\alpha$-stable distributions includes the 
Gaussian ($\alpha=2$), Cauchy
($\alpha =1,\beta =0$),
and \levy~($\alpha =1/2,\beta =1$) families.
Unlike the CF, the PDF of $\alpha$-stable distributions cannot be
expressed in closed form, except in the three special cases mentioned.
This presents significant complications in the development
of effective methodological tools for inference, when models
involve $\alpha$-stable distributions. Nevertheless, as we
describe next, 
a wide variety of relevant statistical tools have been 
proposed in the literature and have been applied in practice.

\subsection{Motivation: Inference with stable distributions}
\label{subsect: param_inference}

The simplest and most common class of inferential procedures 
is probably that of parameter estimation. For systems 
governed by $\alpha$-stable noise this may for example involve 
estimating the parameters $\bm \lambda$, $b$ and $A$ in 
\eqref{eq:AR}-\eqref{eq:AR:noisy} above, as described
in more detail in Section \ref{sec:AR}. Other inference examples 
can be found in the references to specific application 
domains provided above. 

Numerous techniques have been developed 
for estimating the parameter vector $\boldsymbol{\theta} = (\alpha,\, 
\sigma, \, \beta,\, \mu)$ from batch data. 
Common frequentist approaches include
those based on the quantiles of the 
distribution \cite{McCulloch1986}, its logarithmic moments 
\cite{Kuruoglu2001}, the empirical CF \cite{Koutrouvelis1980}, 
approximate maximum likelihood estimators \cite{Nolan2001}, 
or block-maxima scaling \cite{StoevMichailidisTaqqu2011}.
However, since the $\alpha$-stable PDF is not available in closed form,
all the above approaches are {\em approximate.} This issue
similarly affects corresponding Bayesian methods aiming 
at computing the posterior distribution of the parameters,
for which the likelihood function needs to be evaluated; see, e.g., 
\cite{Lombardi2007}.

On the other hand, $\alpha$-stable distributions admit representations
involving {\em latent variables}, enabling (asymptotically) exact Bayesian 
inference.
Schemes such as those based on marginal 
representations of the stable likelihood or on the conditional and 
pseudo-marginal samplers \cite{Buckle1995, Qiou1998, Riabiz2015}, 
belong to this class. Also, the product property 
\cite{Feller1966} and the scale mixture of normals representation 
of symmetric stable distributions can be used, as in 
\cite{Godsill1999, Tsionas1999, GodsillKuruoglu1999, Godsill2000}.
The central object of interest in the present work is yet
another latent variable model, the so-called
{\em Poisson series representation} (PSR) of stable laws,
previously employed in 
\cite{Lemke2014, Lemke2014a, Lemke2015}.

Although the PSR is an exact representation, 
it is an infinite series which itself needs to 
be approximated. For effective inference, 
it is then necessary to quantify the error incurred by such 
an approximation. The study of this approximation error is the 
main aim of this paper. Therefore, the present work 
provides a more firm theoretical
foundation for the Bayesian parameter samplers 
mentioned above, and more generally for Bayesian 
inference in models involving $\alpha$-stable noise.

Studying the error incurred by an approximation to
the PSR is also relevant to state inference in continuous-time  
stochastic dynamics through Bayesian methods such as sequential Monte-Carlo (SMC) 
\cite{Doucet2000, Cappe2007, DoucetJohansen2011}. 
When part of the state is (conditionally) linear and Gaussian, combining  Kalman-filter steps \cite{Kalman1960} with SMC filters results in more efficient samplers in terms of Monte Carlo variance \cite{Doucet2000, Schon2005}. 
As discussed in Section~\ref{section:CLT:st:integral}, the PSR  
extends to $\alpha$-stable \levy~processes and state space models driven by these processes, and it enables efficient SMC inference methods as implemented in 
\cite{Lemke2011}, \cite{Lemke2014a}, \cite{Lemke2015a}, \cite{Riabiz2017a}. 
Note, however, that
the results established here only pertain 
to $\alpha$-stable RVs and \levy~processes; 
continuous-time models driven by $\alpha$-stable \levy~processes
will be examined in future work;
also see Section~\ref{section:CLT:st:integral} for some additional remarks.

\subsection{Main contributions and paper organisation}

The central object of interest in this work is the 
Poisson series representation (PSR) of an $\alpha$-stable RV,
mentioned above.
As described in Section~\ref{sect:PSR}, the PSR is 
an infinite sum of RVs involving the arrival times of 
a Poisson process. Since it is impossible to compute
the entire infinite series in practice,
only approximate versions of the PSR can be employed
for simulation and inference purposes.
The starting point of our approach is the 
truncation of the PSR, followed by the 
approximation of the tail of the series (to which 
we refer as the {\em residual} series) 
by an appropriately chosen Gaussian RV.
We recently noted that such an approximation is asymptotically 
exact, as the truncation point becomes larger \cite{Riabiz2017}.
This CLT is the first main contribution of this work,
given in Section~\ref{sec:CLT_rv_scalar}:
Theorem~\ref{th:CLT}	
provides a precise version of the CLT together with a complete proof,
under conditions weaker than those stated in \cite{Riabiz2017}. 

We then investigate the {\em nonasymptotic} accuracy of the above 
approximation, as a function of the truncation point.
The main tool in this investigation is Ess\'{e}en's smoothing lemma;
this is a classical Fourier-inversion inequality, used to 
translate information on the distance between two characteristic 
functions (CFs) to information about the distance between 
the two corresponding probability distributions;
see, e.g., \cite{RachevKlebanovStoyanovEtAl2013}, 
and the discussion in Section~\ref{section:smoothing}. 
The derivation of explicit expressions for the
CFs of several quantities of interest is our second main
contribution,
given in Section~\ref{section:cf:equations}. 

Let $c$ denote the truncation parameter for the PSR.
In Section~\ref{section:residual_convergence} we establish
nonasymptotic, strong, explicit upper bounds 
on the distance between the distribution of the
PSR residual and an appropriate Gaussian distribution.
In the symmetric ($\beta=0$) case, these results 
in particular imply that the convergence of the CLT
in Section~\ref{sec:CLT_rv_scalar} takes place
at a rate $O(1/c)$, and that it is faster when $\alpha$ approaches 2;
this is consistent with the numerical findings reported
in \cite{Riabiz2017b}. We also establish a different bound
that decays like $1/\sqrt{c}$ asymptotically, but which is
tighter than our previous bound for 
relatively small values of $c$ and $\alpha$.
This is the third main contribution of this paper. 

The Gaussian approximation of the PSR residual suggests an 
elegant approximate representation of a stable RV.
However, bounds on the distance between the residual and 
a Gaussian do {\em not} immediately translate to 
corresponding bounds between this approximate representation
and the corresponding stable law.
Obtaining such bounds is the fourth and perhaps 
most important contribution of this work. 
When using this approximate representation in the 
context of inference procedures with $\alpha$-stable models, 
having explicit bounds facilitates the selection of 
an appropriate value for the truncation parameter, 
in a way that also provides estimation error guarantees 
for the results. For the case of symmetric ($\beta=0$)
stable laws, our bounds are stated in
Section~\ref{section:residual_contribution}.
These results, as well as the numerical
study performed in \cite{Riabiz2017b},
indicate that the approximated PSR is closer to 
the stable distribution for smaller values of $\alpha$, 
and that the rate of convergence now depends on $\alpha$.

There is extensive earlier work on the analysis of the convergence 
rate of the truncated PSR to the corresponding stable law --
as opposed to the convergence of the distribution of the residual studied here; 
see, e.g., \cite{JanickiWeron1994,janicki1992computer,JanickiKokoszka1992,%
LedouxPaulauskas1996, Bentkus1996, Bentkus2001}. 
In Section~\ref{section:residual_contribution} we review the most
relevant of these results, and we derive bounds indicating
that representation proposed in this work
should typically yield a better approximation 
to the stable distribution than 
simply truncating the PSR.

Finally, in Section~\ref{section:linear_models} we illustrate
the utility of our main results with an example of statistical inference.
We recall an MCMC-based inference scheme for the parameters 
of the discrete-time linear models \eqref{eq:AR} and \eqref{state_space},
and we discuss the use of our Gaussian approximation
bounds in this setting. We also briefly 
describe potential extensions of our results
to continuous-time systems and to 
multivariate stable distributions.
In each of these cases, analogues of the PSR
representation have already been established,
and having a CLT for the residual and nonasymptotic 
bounds on the induced approximation, like the ones 
established in this work, would potentially be of 
significant interest in applications.  

The proofs of most of the main results in the paper, together
with the more technical lemmas, are given in the appendices.

\subsection{Notation}

Capital letters, e.g., $X,Y$, are used for RVs, and `hats' denote
approximate versions, e.g., $\hat{X}$ denotes a RV with a distribution 
which is close to that of $X$.
The following notation
is used for some common distributions: $\mathcal{N}(\mu_W, \sigma^2_W)$ 
is the normal distribution with mean $\mu_W$ and variance $\sigma^2_W$; 
$\mathcal{U}(a,b)$ is the uniform distribution on the interval $(a,b)$; 
and ${\rm Poisson}(t)$ is the Poisson distribution with mean $t$. 

Throughout, the sequence $\{\Gamma_j\}$ will denote the successive
arrival times of a unit rate Poisson process. If a RV $X$ is defined
as a series of random terms involving the sequence $\{\Gamma_j\}$,
then $X_{(c,d)}$ will denote the sum of those terms 
corresponding to indices $j$ such that $\Gamma_j \in (c,d)$,
for $0 \leq c < d\leq \infty$. The number of terms in $X_{(c,d)}$ 
is denoted by $N_{(c,d)}$, with the convention that $X_{(c,d)} =0$ 
if $N_{(c,d)}=0$. A subscript notation is used for the moments of such 
RVs, e.g., $\mu_W$ is the mean of $W$, and $m_{(c,d)}$ is the mean 
of $X_{(c,d)}$.

The Kolmogorov distance between two RVs $S$ and $T$ with distribution
functions $F$ and $G$, respectively, is denoted by
$\Delta(S,T):=\sup_x|F(x)-G(x)|$. Upper bounds on $\Delta(S,T)$
derived from the smoothing lemma with a finite smoothing parameter
will be denoted as $I(S,T)$, and when the smoothing parameter
goes to infinity the corresponding bounds will be denoted by
$\bar{I}(S,T)$. Numerically computed values for these bounds
will be denoted by $Q(S, T)$ and $\bar{Q}(S,T)$, respectively. 

Some of the CFs considered in the paper are complex valued and, 
for a fixed argument, can be expressed in polar form 
as $z = r e^{i \theta}$, with
$r>0$ and $\theta \in \mathbb{R}$. To avoid
any ambiguity, we adopt the convention that in all
such expressions $\theta$ is assumed to 
lie in the interval $(-\pi, \pi]$. Then it is possible to 
uniquely invert the exponential of $z$, 
obtaining the \textit{principal value} complex logarithm 
of a CF, $\log(z) = \log(r) + i \theta.$ 
Although not always necessary, for the sake of clarity we 
will always work with principal value complex logarithms. 



Two complex-analytic functions that appear repeatedly in our
analysis are the lower incomplete gamma function, 
\begin{IEEEeqnarray}{rCl} 
	\gamma(s,x)& := & \int_0^x t^{s-1}e^{-t}\dif t, 
	\qquad -s \notin \mathbb{N},x>0,
	\IEEEeqnarraynumspace
	\label{eq:lower_inc_gamma}
\end{IEEEeqnarray}
and the upper incomplete gamma function, 
\begin{IEEEeqnarray}{rCl} 
	\Gamma(s,x)& := & \int_x^\infty t^{s-1}e^{-t}\dif t, 
	\qquad -s \notin \mathbb{N},x>0.
	\IEEEeqnarraynumspace
	\label{eq:upper_inc_gamma}
\end{IEEEeqnarray}
Then, for any $x>0$,
\begin{IEEEeqnarray}{rCl}
	\Gamma(s)& = & \gamma(s,x) + \Gamma(s,x),
	\qquad -s \notin \mathbb{N},		
	\IEEEeqnarraynumspace
	\label{eq:gamma}
\end{IEEEeqnarray}
is the regular (complete) gamma function.   

We use the symbol `$\sim$' to denote the fact that
a RV $X\sim \mathcal{D}$ has distribution $\mathcal{D}$,
but also to denote the following asymptotic relationship:
for two (real) functions  $f_1$ and $f_2$,
we say $f_1(x)\sim f_2(x)$ when 
	$\lim_{x\rightarrow \infty} f_1(x)/f_2(x) =1.$
Finally, $f_1(x)= O( f_2(x))$ as usual means that,
$\limsup_{x\rightarrow \infty}|f_1(x)/f_2(x) | < \infty.$



\newpage

\section{The Poisson Series Representation (PSR)}
\label{sect:PSR}

Let $X\sim\mathcal{S}_\alpha(\sigma,\beta, 0)$ be an
$\alpha$-stable RV for some $\alpha \in (0,2)$, 
$\alpha \neq 1$.
The PSR, originally introduced by 
\levy~and formalised by LePage et al.~\cite{LePage1981,LePage1981a,LePage1989}
states that $X$ admits the representation,
\begin{IEEEeqnarray}{c} 
	X
	\overset{\mathcal{D}}{=} 
	\sum_{j=1}^{\infty} \Gamma_j^{-1/\alpha} W_j
	-\mathbb{E}[W_1]b_j^{(\alpha)}, 
	\label{eq:PSR_RV}
\end{IEEEeqnarray}
where $\overset{\mathcal{D}}{=}$ denotes equality in distribution,
and:
\begin{itemize}
	\item $\{\Gamma_j\}_{j=1}^\infty$ are the arrival times of a unit rate Poisson 
	process, so that the differences $\Gamma_j-\Gamma_{j-1}$,
	$j=1,2,\ldots$, are i.i.d.\ exponential RVs with mean~1;
	\item $\left\lbrace W_j \right\rbrace_{j=1}^{\infty}$ are i.i.d.\ RVs 
	independent of 
	$\{\Gamma_j\}_{j=1}^\infty$, with,
	\begin{IEEEeqnarray}{c} 
		\mathbb{E}[|W_1|^\alpha]< \infty; 
		\label{eq:condition_W}
	\end{IEEEeqnarray}
	\item  $\{b_j^{(\alpha)}\}_{j=1}^\infty$ are constants
	that are non-zero only if $\alpha \in (1,2)$,
	given by:
	\begin{IEEEeqnarray*}{c} 
		b_j^{(\alpha)} =
		\dfrac{\alpha}{\alpha-1 }\left( j^{\frac{\alpha-1}{\alpha}} 
		-(j-1)^{\frac{\alpha-1}{\alpha}} \right).
	\end{IEEEeqnarray*}
\end{itemize}
The exact representation in~\eqref{eq:PSR_RV} 
can be found in \cite[Theorem~1.4.5]{Samoradnitsky1994}.
Observe that this is only 
valid for a \textit{strictly stable} RV~$X$, i.e.,
when the location parameter $\mu =0$. But a stable RV 
$Y\sim\mathcal{S}_\alpha(\sigma,\beta,\mu)$ with $\mu\neq 0$
can simply be obtained from
$X\sim\mathcal{S}_\alpha(\sigma,\beta, 0)$ 
as $Y\overset{\mathcal{D}}{=}X+\mu$. 
We also observe that the constants $b_j^{(\alpha)}$ have
a telescoping nature, so that, for $N\geq 1$,
	\begin{IEEEeqnarray}{c} 
		\sum_{j=1}^{N}b_j^{(\alpha)} =
		\dfrac{\alpha}{\alpha-1 } N^{\frac{\alpha-1}{\alpha}}\mathds{1}{\left(\alpha \in (1,2)\right)},
		\IEEEeqnarraynumspace
		\label{eq:sum_bi}
	\end{IEEEeqnarray}
where $\mathds{1}{\left(\cdot\right)}$ denotes the indicator 
function, equal to 1 if $(\cdot)$ is satisfied, and 0 otherwise. 

Although the distribution of the RVs $\{W_j\}$ above has
not been explicitly described, the $\alpha$-th absolute
moment of $W_1$ can be expressed in terms
of $\alpha,\beta$ and $\sigma$ as follows:
\begin{IEEEeqnarray}{rCl}
	\sigma^\alpha & = &\frac{\mathbb{E}\left[ {\abs{W_1}}^\alpha\right] }{C_\alpha}, 
	\IEEEyesnumber\IEEEyessubnumber
	\label{eq:mapping_sigma}
	\\
	\beta & = &\frac{\mathbb{E}\left[ {\abs{W_1}}^\alpha \sign{W_1}\right]}{\mathbb{E}\left[ {\abs{W_1}}^\alpha\right] },
	\IEEEyessubnumber
	\label{eq:mapping_beta}
\end{IEEEeqnarray}
\noindent
where,
\begin{IEEEeqnarray}{c} 
	C_\alpha = \left( \int_0^\infty x^{-\alpha} \sin{x} \dif x\right)^{-1} = \frac{1-\alpha}{\Gamma(2-\alpha)\cos(\pi \alpha /2)}.
	\label{eq:C_alpha}
	\IEEEeqnarraynumspace
\end{IEEEeqnarray}
See Appendix~\ref{app:sigma_beta} for some more
details about \eqref{eq:mapping_sigma} and \eqref{eq:mapping_beta}, 
particularly when $W_1 \sim \mathcal{N}(\mu_W, \sigma_W^2)$,
a case of special interest since taking
$W_1 \sim \mathcal{N}(0, \sigma_W^2)$
corresponds to the symmetric stable distribution,
namely, the case $\beta =0$. Also note that,
in view of the relations~(\ref{eq:mapping_sigma}),~(\ref{eq:mapping_beta})
and~(\ref{eq:C_alpha}), 
when $W_1 \sim \mathcal{N}(\mu_W, \sigma_W^2)$,
there is a 1-1 relationship between the parameters
$(\alpha,\mu_W,\sigma_W)$ and $(\alpha,\beta,\sigma)$;
Table~\ref{tab:param_fig_PSR} contains a few numerical 
examples to which we will return later.

\begin{table}[ht!]
	\caption{Numerical examples of the correspondence 
between the parameters $(\mu_W,\sigma_W)$ of $W_1$,
and the parameters $(\sigma,\beta)$ of the corresponding
$\alpha$-stable distribution.}
	\label{tab:param_fig_PSR}
	\centering
	\begin{tabular}{|c ||c c||c c|}
		\hline
		& $\mu_W$ & $\sigma_W$ & $\sigma$  &  $\beta$\\
		\hline
		\hline
		\multirow{3}{*}{$\alpha$ = 0.8} & 0 & 1 &  1.16 & 0 \\
		\cline{2-5}
		& 1 & 1 &  1.71 & 0.84 \\
		\cline{2-5}
		& 1 & 0 & 1.42 & 1 \\
		\hline
		\hline
		\multirow{3}{*}{$\alpha$ = 1.2} & 0 & 1 & 1.37  &  $0$\\
		\cline{2-5}
		& 1 & 1 & 1.99  &  $0.99$\\
		\cline{2-5}
		& 1 & 0 & $1.80$  &  1\\
		\hline
	\end{tabular}
\end{table}

\begin{figure*}[ht!] 
	\centerline{
		\includegraphics[width=0.25\linewidth]{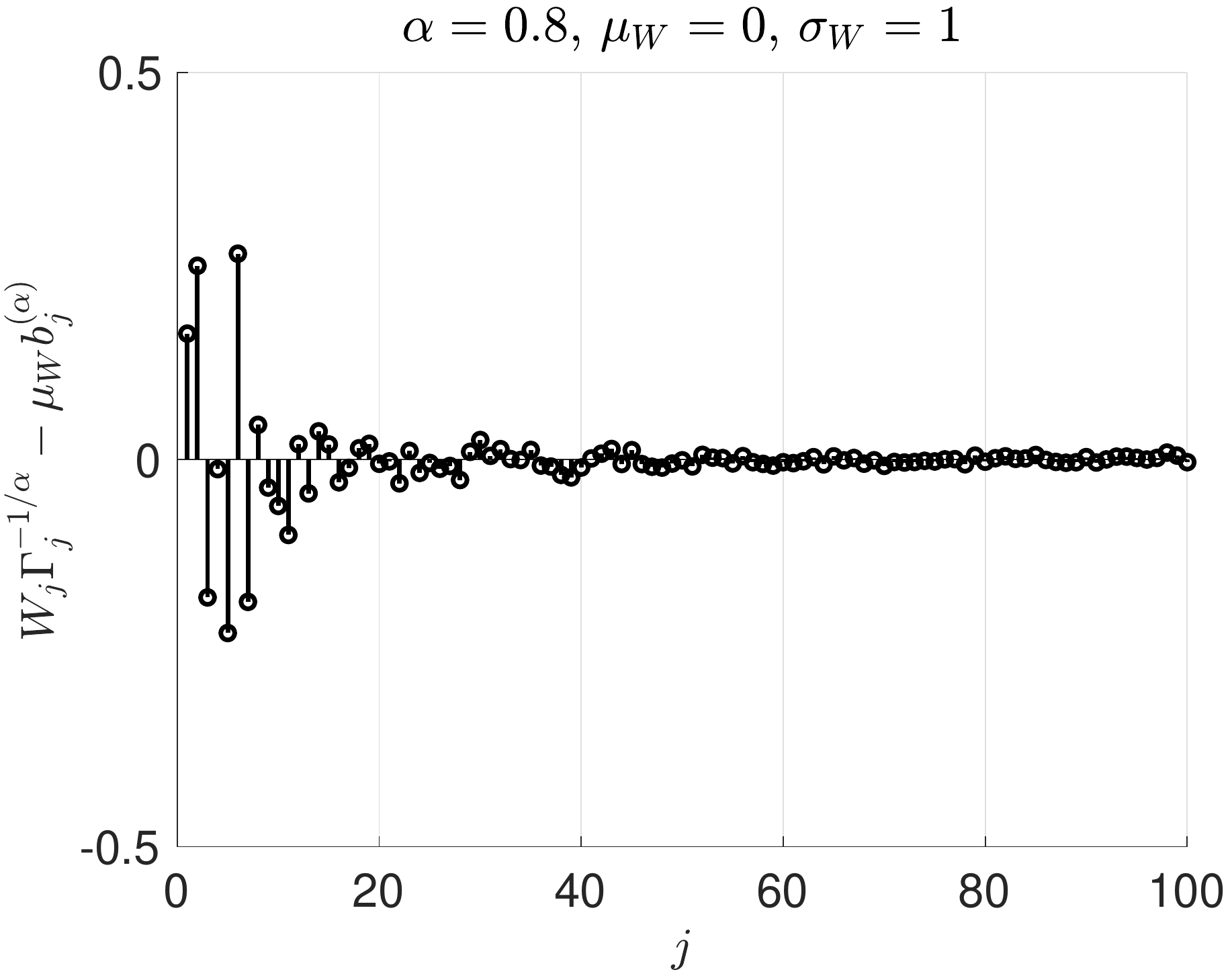}		
		\includegraphics[width=0.25\linewidth]{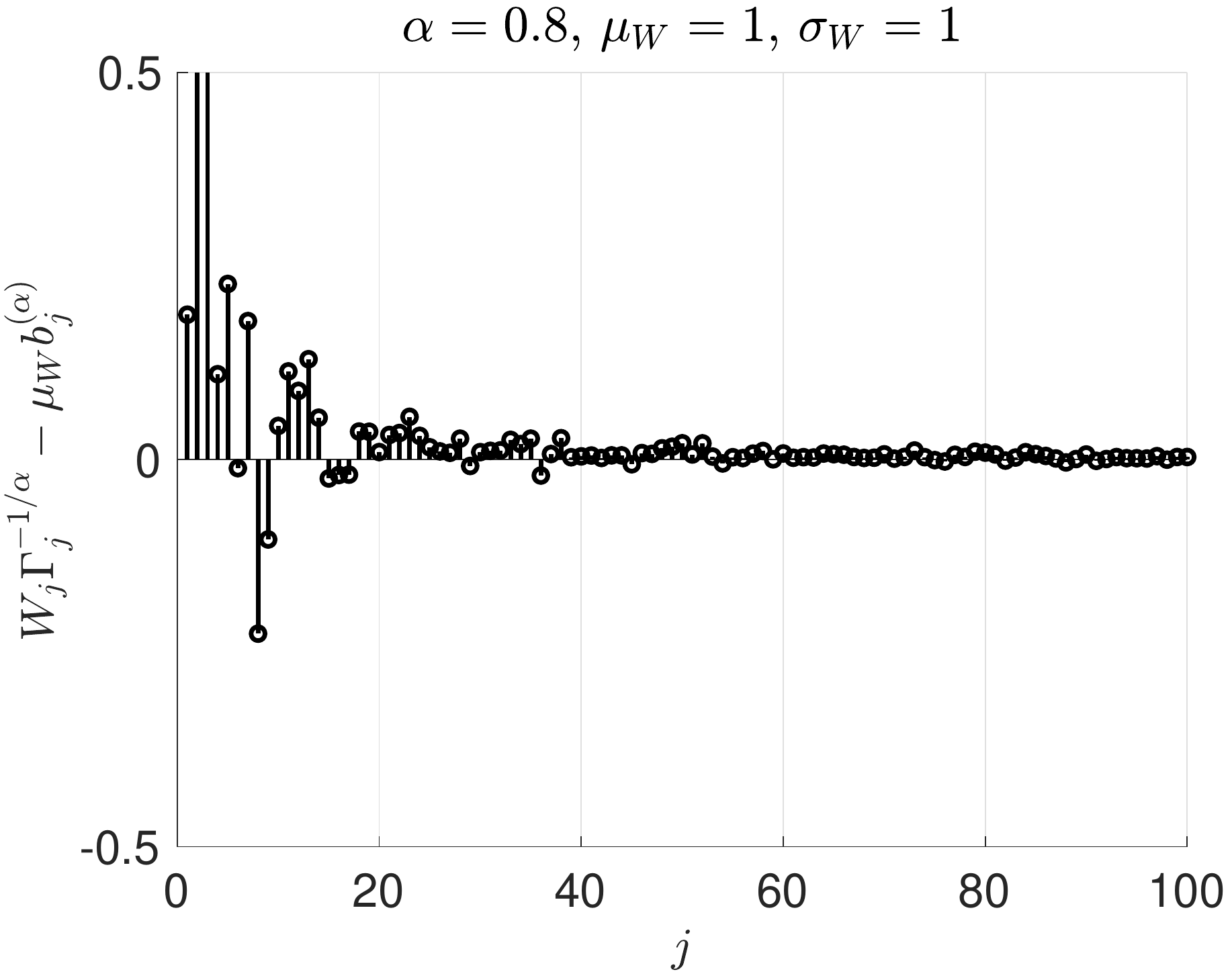}	
		\includegraphics[width=0.25\linewidth]{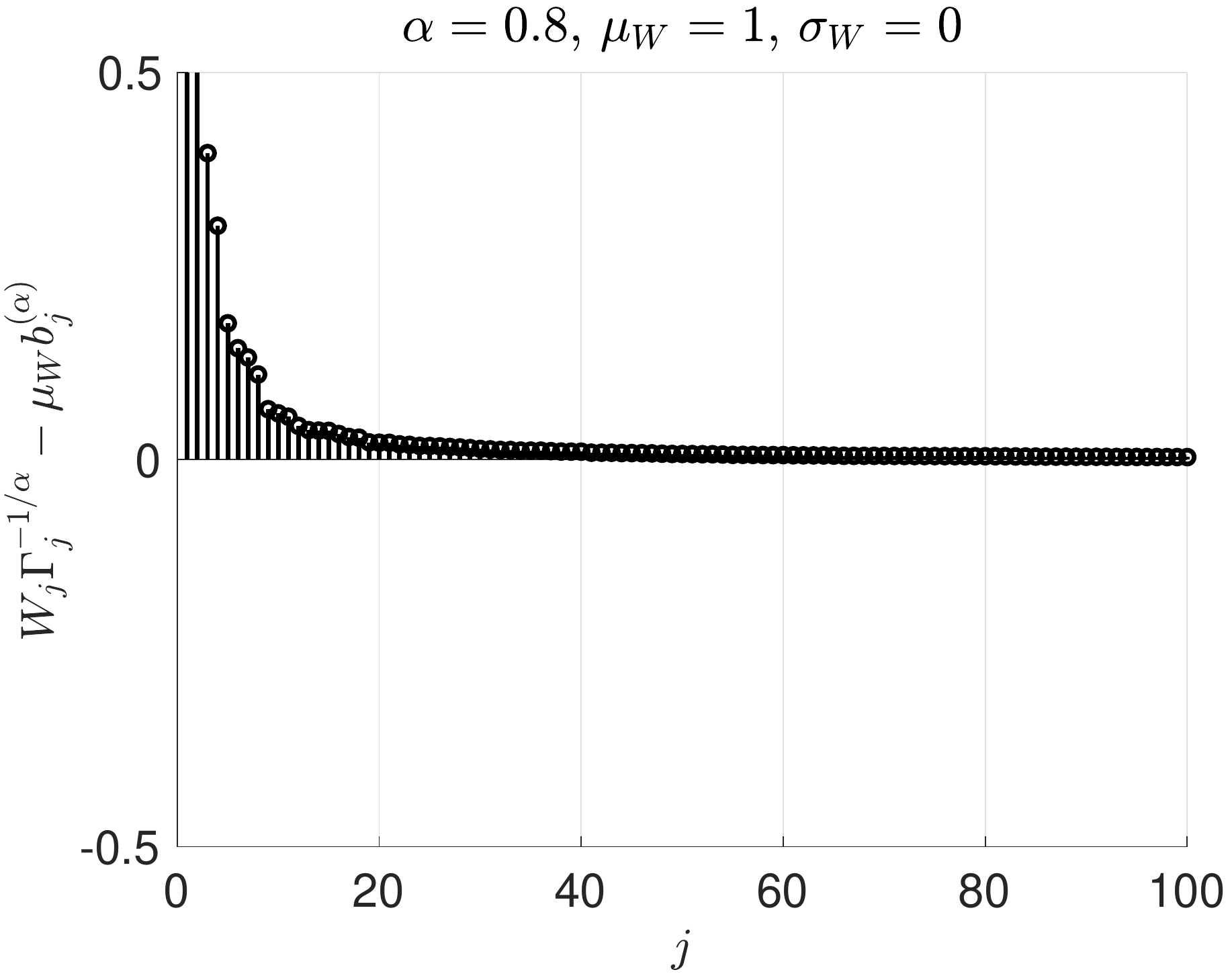}
	}
	\centerline{
		\includegraphics[width=0.25\linewidth]{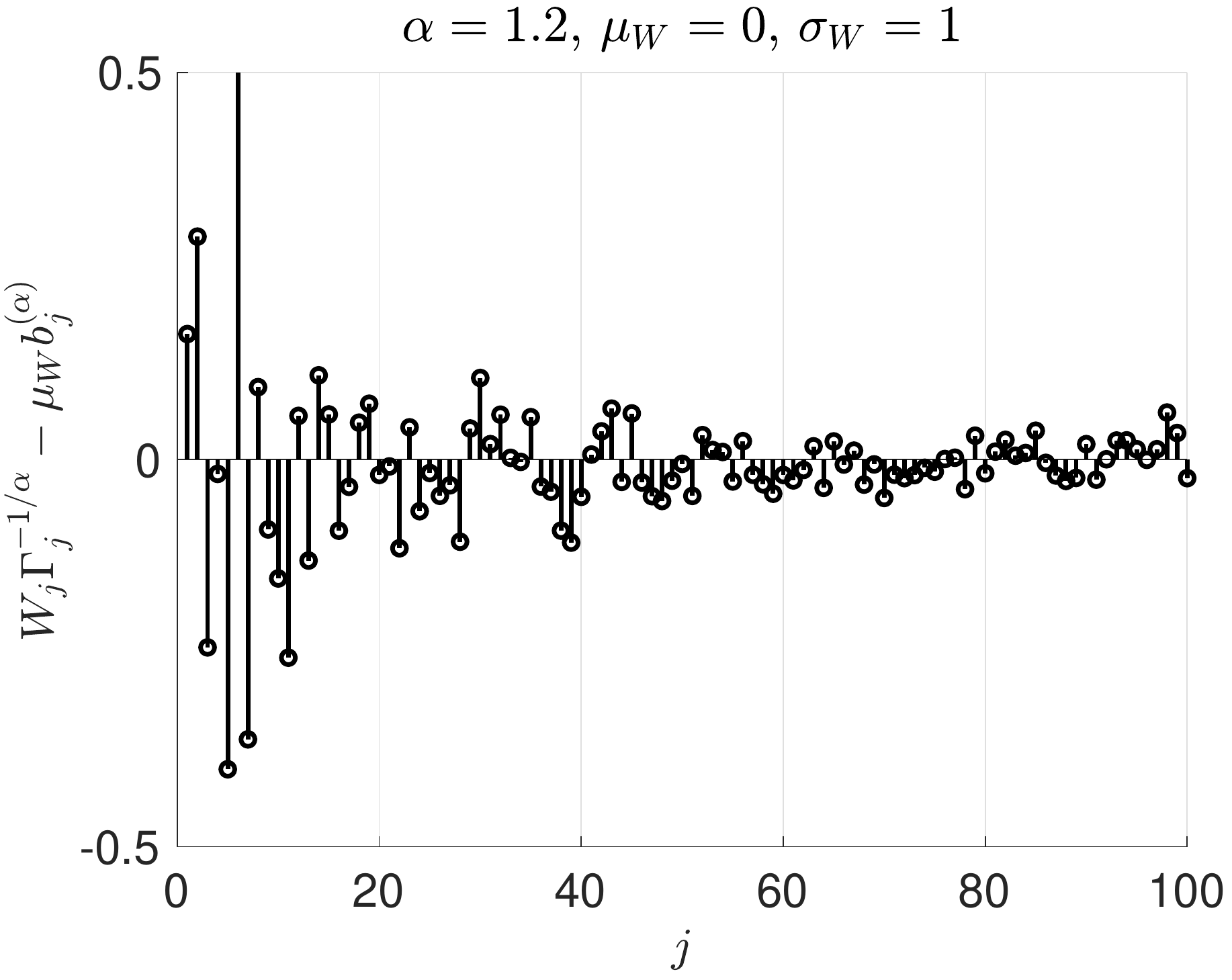}		
		\includegraphics[width=0.25\linewidth]{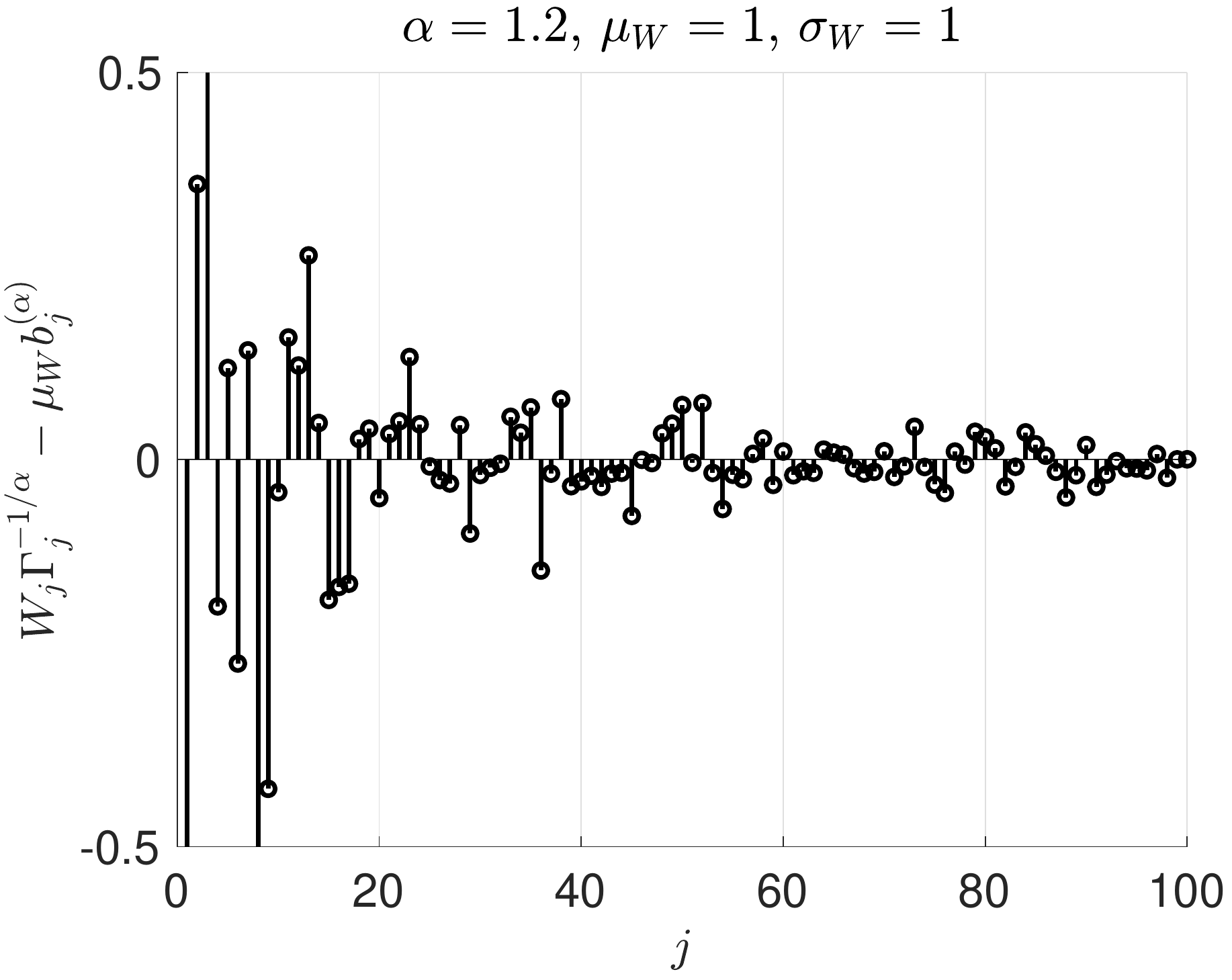}	
		\includegraphics[width=0.25\linewidth]{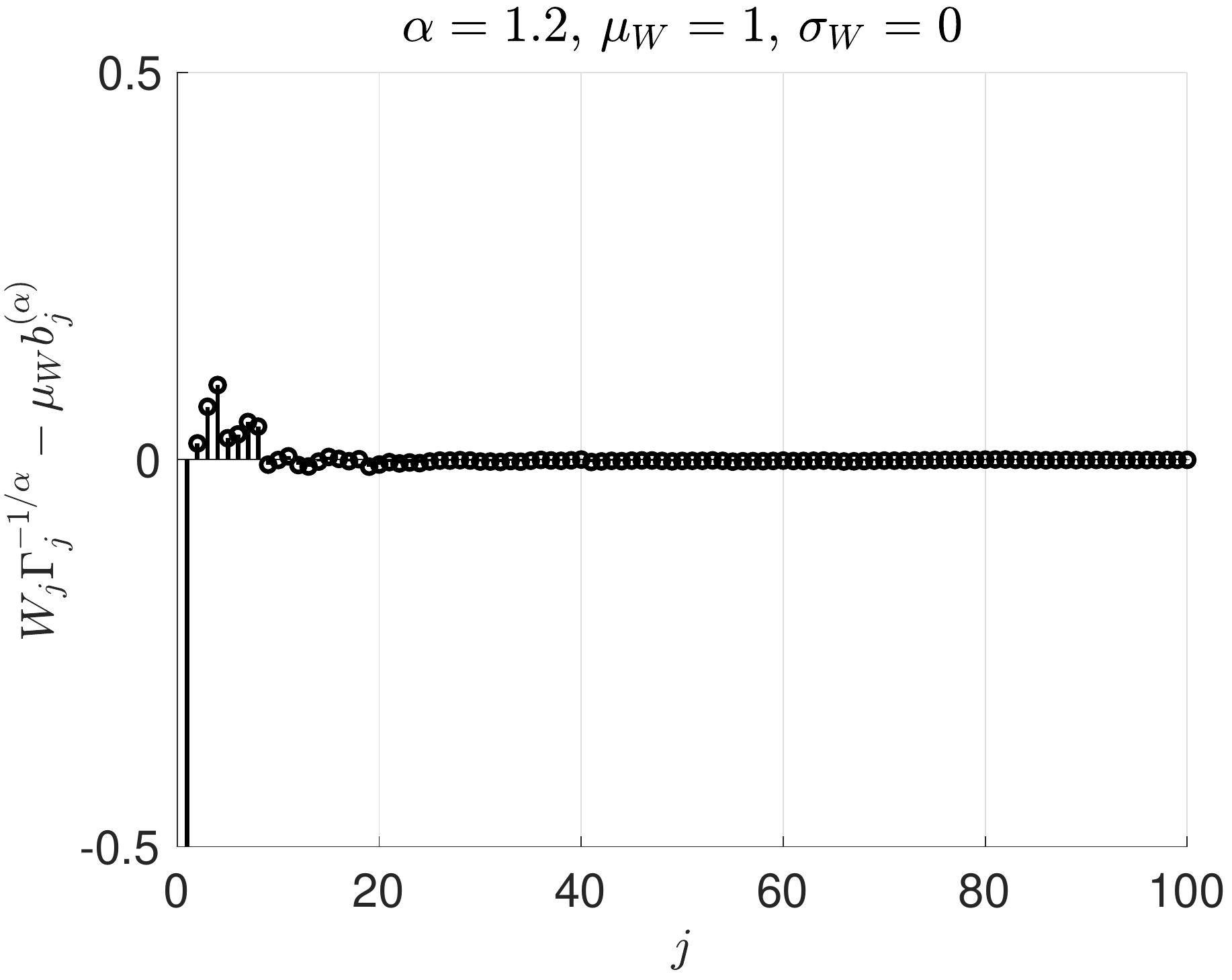}
	}
	\caption{First 100 terms of a PSR realization with 
	$W_1 \sim \mathcal{N}(\mu_W, \sigma^2_W)$,
	for different values of the parameters $\alpha,\mu_W$ and $\sigma^2_W$.
	Top row: $\alpha = 0.8$ (heavier tails); 
	bottom row: $\alpha = 1.2$ (lighter tails).
	}
	\label{fig:jumps_w}
\end{figure*}

Figure~\ref{fig:jumps_w} shows a realization of the first 
100 terms {$W_j \Gamma_j^{-1/\alpha} - \mu_Wb_j^{(\alpha)}$} of the PSR
with $W_1 \sim \mathcal{N}(\mu_W, \sigma_W^2)$,
for different values of the parameters $\alpha,\mu_W$
and $\sigma_W$. The corresponding values of $\sigma$
and $\beta$ are shown in Table~\ref{tab:param_fig_PSR}.
Since the sequence
$\{\Gamma_j\}_{j=1}^\infty$ is increasing with probability one,
the terms $\{\Gamma_j^{-1/\alpha}\}_{j=1}^{\infty}$ are
decreasing, and the summands in the PSR are therefore
stochastically decaying (in absolute value). 
This is indeed confirmed by the sample paths shown
in Figure~\ref{fig:jumps_w}.

In most of the paper we will
focus on the case $W_1 \sim \mathcal{N}(\mu_W, \sigma^2_W)$, 
since this allows a conditionally Gaussian representation of 
the $\alpha$-stable distribution which is useful for the inference 
tasks mentioned earlier.
In fact, 
from the PSR \eqref{eq:PSR_RV}  it follows that, 
if $W_1 \sim \mathcal{N}(\mu_W,\sigma_W^2)$, 
we can  write an auxiliary variable model for $X$ as,
\begin{IEEEeqnarray}{rCl} 
	X | \{\Gamma_j\}_{j=1}^\infty  & \sim & 	
	\mathcal{N}
	\left(
	\mu_W  m,\sigma^2_W {S}^2	\right), \label{eq:inf_cond_gauss_X}\\
	m & := & \sum_{j=1}^\infty\Gamma_j^{-1/\alpha}-b_j^{(\alpha)}, \nonumber\\
	S^2 & := & \sum_{j=1}^\infty\Gamma_j^{-2/\alpha}. \label{eq:inf_variance_X}
\end{IEEEeqnarray}
In this model $m$ and $S^2$  are treated as auxiliary RVs, and $X$ has 
a conditionally Gaussian structure. This means that standard auxiliary 
variable methods for conditionally Gaussian models may be readily 
applied, for example blocked and collapsed Gibbs samplers \cite{Dyk2008}, 
Monte Carlo EM \cite{Tanner1997} and Rao-Blackwellised particle 
filters \cite{Doucet2000, Schon2005}. We note that the general framework here is a  {continuous} {\em scale\/} and {\em mean\/} mixture of normals, since the density {$p$}
of $X$, can be expressed in obvious notation as:
\begin{IEEEeqnarray}{c} 
	p=\int_{S\in\mathbb{R}^+}\int_{m\in\mathbb{R}}\mathcal{N}
	\left(
	\mu_W m,\sigma^2_W S^2\right)p(m,S)\dif m \dif S.
	\label{eq:gauss_mixture}
\end{IEEEeqnarray}

\newpage

\section{Truncation of the PSR}
\label{sec:PSR:truncation}

While the exact representations of the stable law 
in~(\ref{eq:PSR_RV}) and~\eqref{eq:gauss_mixture} are very appealing, 
they are computationally intractable because of the infinite 
summations involved. Since
the summands of the PSR~(\ref{eq:PSR_RV}) are stochastically decaying,
a first intuitive approach would be to consider a large but finite number 
of summands so that the distribution of the truncated PSR,
\begin{IEEEeqnarray}{c} 
	X_{(0,c)}:=\sum_{j:\Gamma_j \in [0,c]} W_j\Gamma_j^{-1/\alpha},
	\label{eq:X_0c}
\end{IEEEeqnarray}
is `close enough' to that of $X$, e.g., in terms of their
Kolmogorov distance. A number of works 
in the literature {have} been devoted to the study 
of this approximation, 
and  {in 
Section~\ref{section:residual_contribution} we will recall the main results that are used as a comparison with ours. Here we mention that the truncated PSR has been used in engineering applications
by \cite{Azzaoui2010} for the task of generating  
stable variables, and the choice of the truncation parameter is also based on analysis in the frequency domain.}

However, more accurate results can be obtained by taking into
account the residual part of the series, at least approximately.
Making this idea precise is the main aim of this paper:
For a given truncation parameter $c>0$, let
$X_{(0,c)}$ be defined as in (\ref{eq:X_0c}) above,
so that,
\begin{IEEEeqnarray}{c} 
	X		
	\,{\buildrel \over =}\,  X_{(0,c)} 
	+ R_{(c,\infty)}, 
	\label{exact_PSR_split}
\end{IEEEeqnarray}	
\noindent
where the PSR residual $R_{(c,\infty)}$ is,
\begin{IEEEeqnarray}{c} 
	R_{(c,\infty)} := \lim_{d\rightarrow\infty}R_{(c,d)},
	\label{eq:true_residual}
\end{IEEEeqnarray}
with,
\begin{IEEEeqnarray}{c} 
	R_{(c,d)}:=
	\sum_{j:\Gamma_j \in(c,d)} W_j\Gamma_j^{-1/\alpha}
	-\mathbb{E}[W_1]\sum_{j=1}^{\floor{d}}b_j^{(\alpha)},
	\IEEEeqnarraynumspace
	\label{eq:reminder_PSR_finite}
\end{IEEEeqnarray}
and $\floor{\cdot}$ denoting the lower integer part.  
Our first main result will establish that the approximation
of $X$ by `$X_{(0,c)}$ plus an appropriately chosen Gaussian'
is asymptotically (as $c\to\infty$) exact, as long as
{$\mathbb{E}[W_1^2]<\infty$}.

To simplify the notation, from now on we will assume that $d$ is integer 
so that $\floor{d} =d$. Notice that the limit in \eqref{eq:true_residual} 
denotes convergence in distribution, and is guaranteed to exist by 
the fact that the full series is known to 
converge to an $\alpha$-stable RV 
\cite[Theorem~1.4.5]{Samoradnitsky1994}.



Next, we compute the first and second order moments
of $R_{(c,d)}$  and of $R_{(c,\infty)}$. 
Since,
conditioned on the Poisson number of events,
the ordered arrival times $\{\Gamma_j\}$ of a unit rate
Poisson process may 
equivalently be written as an unordered set of i.i.d.\ 
{\em uniformly\/} distributed random variables 
$\{U_j\}$ \cite{Kingman1992}, a generative model for $R_{(c,d)}$ 
is as follows,
\begin{IEEEeqnarray}{rCl} 
	R_{(c,d)} &=& \sum_{j=1}^{N_{(c,d)}} Y_j - B, \label{eq:Rcd_new}   \\
	N_{(c,d)} &\sim & {\rm Poisson}(d-c), \nonumber\\
	Y_j &:=&W_j U_j^{-1/\alpha}, \label{eq:Y_j}\\
	W_j &\sim & \ F, \nonumber\\
	U_j &\sim & \mathcal{U}(c,d), \label{eq:U_j_c_d}\\
	B&:=& \mathbb{E}[W_1] \dfrac{\alpha}{\alpha-1 } d^{\frac{\alpha-1}{\alpha}}\mathds{1}{\left(\alpha \in (1,2)\right)},
	\label{eq:B}
\end{IEEEeqnarray}
where ${\rm Poisson}(d-c)$ is the Poisson distribution with mean $d-c$, 
$F$ is the distribution of $W_1$ which is assumed to satisfy
\eqref{eq:condition_W}, and the
expression for $B$ is obtained from \eqref{eq:sum_bi}.
In other words, we can think of $R_{(c,d)}$ as a compound Poisson process
containing two sources of randomness: The 
random Poisson number of arrivals,~$N_{(c,d)}$, and the values of
the RVs $Y_j$ being summed. Based on this observation,
the following lemma is proved in Appendix~\ref{app:finite_res_momemts} 

\medskip

\begin{lemma}	
\label{lem:finite_res_momemts}
Suppose $\mathbb{E}[W_1^2] <\infty$, let $Y_1$ be the RV
defined in \eqref{eq:Y_j} with CF $\phi_{Y_1}(s)$,
and let $B$ be the constant in \eqref{eq:B}. 
Then $\phi_{R_{(c,d)}}(s)$, the CF of $R_{(c,d)}$, is:
	\begin{IEEEeqnarray}{c} 
		\phi_{R_{(c,d)}}(s) =	
		\exp\left((d-c)\big(\phi_{Y_1}(s)-1\big)
		-
		is B\right). 
		\IEEEeqnarraynumspace
		\label{eq:cf_R_cd}
	\end{IEEEeqnarray}
Moreover, the mean and variance of the residual $R_{(c,d)}$ are:
	\begin{IEEEeqnarray}{rCl} 
		{m_{(c,d)}} & = &
		\mathbb{E}[W_1]\frac{\alpha}{\alpha-1}\left(d^{\frac{\alpha-1}{\alpha}}-c^{\frac{\alpha-1}{\alpha}}\right) - B, 
		\IEEEeqnarraynumspace
		\label{eq:m_cd}
		\\
		S^2_{(c,d)} & = &
		\mathbb{E}[W_1^2]
		\frac{\alpha}{\alpha-2}\left(d^{\frac{\alpha-2}{\alpha}}-c^{\frac{\alpha-2}{\alpha}}\right).  
		\IEEEeqnarraynumspace
		\label{eq:S_cd}		
	\end{IEEEeqnarray}
\end{lemma}	

\medskip

Note that, we can take the limit as $d\to\infty$ in the last two expressions 
in the lemma to obtain,
\begin{IEEEeqnarray}{rCl} 
	m_{(c,\infty)} & := & \lim_{d\rightarrow \infty}m_{(c,d)}
	=
	\mathbb{E}[W_1]\frac{\alpha}{1-\alpha}  c^{\frac{{\alpha-1}}{\alpha}},
	\IEEEeqnarraynumspace
	\label{eq:m_c_infty}
	\\
	S^2_{(c,\infty)}& := &  \lim_{d\rightarrow \infty} S^2_{(c,d)}
	= \mathbb{E}[W_1^2]\frac{\alpha}{2-\alpha} c^{\frac{{\alpha-2}}{\alpha}}.
	\label{eq:S_c_infty}
	\IEEEeqnarraynumspace
\end{IEEEeqnarray}
\noindent 
In the following section we will take a similar limit for the
CF of $R_{(c,\infty)}$,
\begin{IEEEeqnarray}{rCl} 
\phi_{R_{(c,\infty)}}(s) = \lim_{d\rightarrow \infty} \phi_{R_{(c,d)}}(s),
\label{eq:lim_CF_res}
\end{IEEEeqnarray}
where the existence of the limit is guaranteed by the existence of the PSR. 


\newpage

\section{Asymptotic Normality of the PSR Residual}
\label{sec:CLT_rv_scalar}

Although it is easy to see that $R_{(c,\infty)}$ is not Gaussian,
the following CLT-like result states its asymptotic normality as $c \rightarrow \infty$. A first version of this result for the special case $W_j\equiv1$ was 
presented in \cite{Lemke2014}, and the general case with random $W_j$ was 
stated in~\cite{Riabiz2017}. 
Here 
we provide a precise statement of the claim together with a complete proof,
under milder moment conditions than those in \cite{Riabiz2017}. 

\medskip

\begin{theorem}
	\label{th:CLT}	
	Let $R_{(c,\infty)}$, $m_{(c,\infty)}$ and $S^2_{(c,\infty)}$ 
be defined as in \eqref{eq:true_residual},
\eqref{eq:m_c_infty} and~\eqref{eq:S_c_infty}, respectively.
	If {$\mathbb{E}[W_1^2]< \infty$},
	then,
	\begin{IEEEeqnarray}{c} 
		Z_{(c,\infty)} \defeq
		\frac{R_{(c,\infty)} - m_{(c,\infty)}}{S_{(c,\infty)}}
		\underset{c\rightarrow \infty}{\overset{\mathcal{D}}{\longrightarrow}} Z,
		\IEEEeqnarraynumspace
		\label{eq:CLT}
	\end{IEEEeqnarray}
	where $Z\sim \mathcal{N}({0},1)$ and $		\underset{c\rightarrow \infty}{\overset{\mathcal{D}}{\longrightarrow}}$ denotes convergence in distribution, as $c\rightarrow\infty$.	
\end{theorem}

\medskip

\begin{IEEEproof}
The proof is based on the \levy~continuity theorem
\cite{Feller1966}: We will show that, for any fixed 
$s \in\mathbb{R}$, the CF of $Z_{(c, \infty)}$, $\phi_{Z_{(c,\infty)}}(s)$, 
converges to the CF of $Z$, $\phi_Z(s)=\exp(-s^2/2)$,
as~$c\rightarrow \infty$. 
First we express the CF
	of $Z_{(c,\infty)}$
	in terms of the CF 
	of $Y_1$, defined in \eqref{eq:Y_j}.
	Using~\eqref{eq:CLT}, by a change of variables we have,
	\begin{IEEEeqnarray}{rCl} 
		\phi_{Z_{(c,\infty)}}(s) 
		& = &
		\exp\left(-i\frac{ {m_{(c,\infty)}}}{{S_{(c,\infty)}}}s\right) \phi_{R_{(c,\infty)}}\left(\frac{{s}}{{S_{(c,\infty)}}}\right),
		\label{eq:cf_res_standardized}
	\end{IEEEeqnarray}
	and taking the limit as in \eqref{eq:lim_CF_res} and 
	using~\eqref{eq:cf_R_cd},  
	\begin{IEEEeqnarray*}{rCl} 
		\phi_{R_{(c,\infty)}}\left({{s}}\right)
		& = &
		\lim_{d \rightarrow \infty}\phi_{R_{(c,d)}}\left({{s}}\right) \nonumber
		\\
		& =& \lim_{d \rightarrow \infty} \exp\left(  (d-c)\left(\phi_{Y_1}\left({s}\right)-1\right) - {i B s}
		\right), 
		\IEEEeqnarraynumspace
	\end{IEEEeqnarray*}
	\noindent
	or, equivalently,
	\begin{IEEEeqnarray}{rCl} 
		\log(\phi_{R_{(c,\infty)}}\left({{s}}\right))
		& =& \lim_{d \rightarrow \infty} \left(  (d-c)\left(\phi_{Y_1}\left({s}\right)-1\right) - {i B s}
		\right). 
		\label{eq:log_cf_Z}
		\IEEEeqnarraynumspace
	\end{IEEEeqnarray}
	{
	By Lemma 3.3.19 in \cite{Durrett2019},
		\begin{IEEEeqnarray}{rCl}
				\left| e^{ix} - \sum_{k=0}^n \frac{(ix)^m}{m!} \right| \leq \min\left\lbrace \frac{|x|^{n+1}}{(n+1)!}, \, \frac{2|x|^n}{n!}\right\rbrace. 
			\label{eq:durrett} 
			\end{IEEEeqnarray}
	Therefore, 
	we have the following bound 
	on the difference between  $\phi_{Y_1}(s)$ and its 
	second-order Taylor expansion at zero
based on \eqref{eq:durrett} and Jensen's inequality (applied to the absolute value function),
	\begin{IEEEeqnarray}{rCl} 
		\left| \phi_{Y_1}\left({s}\right) 
		- \sum_{k=0}^{2}    \frac{i^k \mathbb{E}[Y_1^k]}{k!}{s}^k \right| 
		& \leq &  
		\mathbb{E} \left[ \min \left\lbrace  \frac{|s|^3 |Y_1|^3}{3!} ,\, { s^2 Y_1^2} \right\rbrace\right] \nonumber
		\\
		& = &  
		\mathbb{E}_{W_1} \mathbb{E}_{U_1\sim \mathcal{U} (c,d)} \left[ \min \left\lbrace  \frac{|s|^3|W_1|^3}{3!}{U_1^{-3/\alpha}} ,\, {s^2 W_1^2}{U_1^{-2/\alpha}} \right\rbrace\right] \nonumber
		\\
		& \leq & 
		\mathbb{E}_{W_1}  \left[ \min \left\lbrace  \frac{|s|^3|W_1|^3}{3!}\mathbb{E}_{U_1\sim \mathcal{U} (c,d)}\left[ {U_1^{-3/\alpha}}\right]  , \, {s^2 W_1^2}\mathbb{E}_{U_1\sim \mathcal{U} (c,d)}\left[ {U_1^{-2/\alpha}}\right]  \right\rbrace\right], 		
		\label{eq:taylor_third} 
	\end{IEEEeqnarray}
where in the equality we have made explicit the distributions with respect to which we are taking the expected value, based on the definition of $Y_1$ \eqref{eq:Y_j}, and the second inequality trivially follows from $\min\{a, b\} \leq a$ and $\min\{a, b\} \leq b$ , for any $a, \, b \in \mathbb{R}$. 
}
	\noindent
In order to further bound the above right-hand side, we recall that from
\eqref{eq:non_asympt_mean} and 
	\eqref{eq:non_asympt_var} in the proof of
	 Lemma~\ref{lem:finite_res_momemts},
	\begin{IEEEeqnarray}{C} 
		(d-c) {\mathbb{E}[Y_1^k]}= 
		\begin{cases}
			{m_{(c,d)}+B}, & k = 1,
			\\
			S^2_{(c,d)}, & k = 2,
		\end{cases}
		\IEEEeqnarraynumspace
		\label{eq:moments_Y}
	\end{IEEEeqnarray}
{
	while, from \eqref{eq:moments_uniform},
	\begin{IEEEeqnarray*}{C} 
		(d-c)\mathbb{E}{[}U_1^{-k/\alpha}{]}= \frac{\alpha}{\alpha-{k}}\left(d^{\frac{\alpha-k}{\alpha}}-c^{\frac{\alpha-k}{\alpha}}\right).
	\end{IEEEeqnarray*}}
	\noindent	
Multiplying~\eqref{eq:taylor_third} by $(d-c)$ and substituting these yields,
	 \begin{IEEEeqnarray*}{l}			
		\left| (d-c) \left(\phi_{Y_1}\left({s}\right)-1\right) - i(m_{(c,d)} + B)s + \frac{s^2S^2_{(c,d)}}{2}\right| 
		\\
		\hspace{1 cm}
	 \leq
	 {
		\mathbb{E}  \left[ \min \left\lbrace  \frac{|s|^3|W_1|^3}{3!}\frac{\alpha}{\alpha-3}\left(d^{(\alpha-3)/\alpha}-c^{(\alpha-3)/\alpha}\right), \, { s^2 W_1^2}\frac{\alpha}{\alpha-2}\left(d^{(\alpha-2)/\alpha}-c^{(\alpha-2)/\alpha}\right) \right\rbrace\right], 		
	}
	\end{IEEEeqnarray*}
and taking the limit as $d \rightarrow \infty$
and using \eqref{eq:log_cf_Z}, gives,
	\begin{IEEEeqnarray}{rCl}			
		\left|\log(\phi_{R_{(c,\infty)}}\left({{s}}\right)) - im_{(c,\infty)} s + \frac{s^2S^2_{(c,\infty)}}{2}\right| 
		& \leq &
		 {
			\mathbb{E}  \left[ \min \left\lbrace  \frac{|s|^3|W_1|^3}{3!}\frac{\alpha}{3-\alpha}c^{(\alpha-3)/\alpha}, \, { s^2 W_1^2}\frac{\alpha}{2- \alpha}c^{(\alpha-2)/\alpha} \right\rbrace\right].
		}
		\label{eq:inequality_phi_Y}	
	\end{IEEEeqnarray}
	\noindent
Finally, replacing $s$ by 
$s/S_{(c, \infty)}$ in \eqref{eq:inequality_phi_Y},
and using \eqref{eq:cf_res_standardized} and \eqref{eq:S_c_infty}, 
we obtain,
	\begin{IEEEeqnarray*}{rCl} 
		\left|\log(\phi_{Z_{(c,\infty)}}\left({{s}}\right))  + \frac{s^2}{2}\right| 
		& \leq &
				 {
				 	s^2
			\mathbb{E}  \left[ \min \left\lbrace  \frac{|s|}{3!}			
			\frac{|W_1|^3}{\mathbb{E}[W_1^2]^{3/2}}
				\frac{\frac{\alpha}{3-\alpha}}
			{\left(\frac{\alpha}{2-\alpha}\right)^{3/2}} c^{-1/2}
			,\,   \frac{W_1^2}{\mathbb{E}[W_1^2]} \right\rbrace\right]
		}
	\end{IEEEeqnarray*}
{	
We can apply the dominated convergence theorem to the argument of the expectation,  
given  that this is bounded by  $\frac{W_1^2}{\mathbb{E}[W_1^2]}$ (which  is integrable  by assumption) 
and it vanishes as $c\to\infty$. Thus the limit of the expectation exists and it is also zero.  Hence
\begin{eqnarray}	
\log(\phi_{Z_{(c,\infty)}}\left({{s}}\right)) \to  
\log\phi_Z(s):=- \frac{s^2}{2},
\label{eq:Z:lim_dc}
\end{eqnarray}
	so that
	$\phi_{Z_{(c,\infty)}}\left(s\right) \to \phi_Z(s),$
as required.
}
\end{IEEEproof}

\medskip


\subsection{Gaussian approximation of the residual}
Theorem \ref{th:CLT} offers an asymptotic justification
for the Gaussian approximation of the PSR residual,
\begin{IEEEeqnarray}{c} 
	\hat{R}_{(c,\infty)}\sim \mathcal{N}
	\big(
	{m_{(c,\infty)}}, {S^2_{(c,\infty)}}
	\big),
	\IEEEeqnarraynumspace
	\label{eq:tot_gauss}
\end{IEEEeqnarray}
discussed earlier in the 
context of practical inference procedures.
Notice that this  approximation matches the first two moments of 
	$\hat{R}_{(c,\infty)}$
to those of 
the exact residual $R_{(c,\infty)}$ for any value of the truncation parameter $c$, and that $R_{(c,\infty)}$ converges in distribution to its 
Gaussian approximation as $c\rightarrow \infty$.
Then, by analogy with \eqref{exact_PSR_split}, we can introduce the 
following RV,
\begin{IEEEeqnarray}{rCl} 
	\hat{X}		
	&{\buildrel \over :=} &  X_{(0,c)} 
	+ \hat{R}_{(c,\infty)}, 
	\IEEEeqnarraynumspace
	\label{proxy_PSR_split}
\end{IEEEeqnarray}
that
converges in distribution to $X\sim \mathcal{S}_\alpha(\sigma, \beta,0)$, as~$c\rightarrow\infty$. 

Note that Theorem~\ref{th:CLT} does not assume that the $W_j$ are Gaussian. 
However, when they are, we have the following overall approximate 
conditionally Gaussian structure for the model, which in part
justifies our focus on the case $W_j \sim\mathcal{N}(\mu_W, \sigma^2_W)$ 
in the rest of the paper. 

\subsection{Approximate conditionally Gaussian representation}
\label{subsect:approx_cond_gauss}

Suppose $W_j \sim \mathcal{N}(\mu_W, \sigma^2_W)$ and that only 
the finite collection of values $\{\Gamma_j\leq c\}$ is known, 
a much more realistic assumption than knowing all the
values in the infinite sequence $\{\Gamma_j\}$. 
Then the $\alpha$-stable distributed 
RV $X\sim \mathcal{S}_\alpha(\sigma, \beta, 0)$ has the approximate 
conditionally Gaussian representation,
\begin{IEEEeqnarray}{rCl} 
	X |\{\Gamma_j\leq c\}  &\overset{\text{approx}}{\sim}
	& \mathcal{N}
	\left( 
	m_{(0,c)} + m_{(c,\infty)}, \,\,   S^2_{(0,c)}  +
	S^2_{(c,\infty)}
	\right), 
	\label{eq:proxy_gauss_model}
	\\
	m_{(0,c)}&:= &	\mu_W \sum_{j:\Gamma_j \in[0,c]}\Gamma_j^{-1/\alpha}, \nonumber
	\\
	S^2_{(0,c)}&:= & \sigma^2_W \sum_{j:\Gamma_j \in[0,c]}\Gamma_j^{-2/\alpha}, \nonumber
\end{IEEEeqnarray}
where $m_{(c,\infty)}$ and $S^2_{(c, \infty)}$ are given in \eqref{eq:m_c_infty} and \eqref{eq:S_c_infty}.
{Thus, by analogy with \eqref{eq:gauss_mixture}, the density $p$ of $X$ can be approximately represented~as  
	\begin{IEEEeqnarray}{c} 
		p{\approx}\int_{S_{(0,c)} \in \mathbb{R}^+}\int_{m_{(0,c)} \in\mathbb{R}}\mathcal{N}
		\left( 
		m_{(0,c)} + m_{(c,\infty)}, \,\,   S^2_{(0,c)}  +
		S^2_{(c,\infty)}
		\right){p}_{({m}_{(0,c)}, {S}_{(0,c)})}({m}_{(0,c)}, {S}_{(0,c)})\dif {m}_{(0,c)} \dif {S}_{(0,c)},
		\label{eq:gauss_MS_mixture_proxy}
	\end{IEEEeqnarray}
implying that $X$ can also be approximated by a continuous mean and scale mixture of normals. 
}

\begin{figure*}[ht!] 
	\centerline{
		\includegraphics[width=0.3\linewidth]{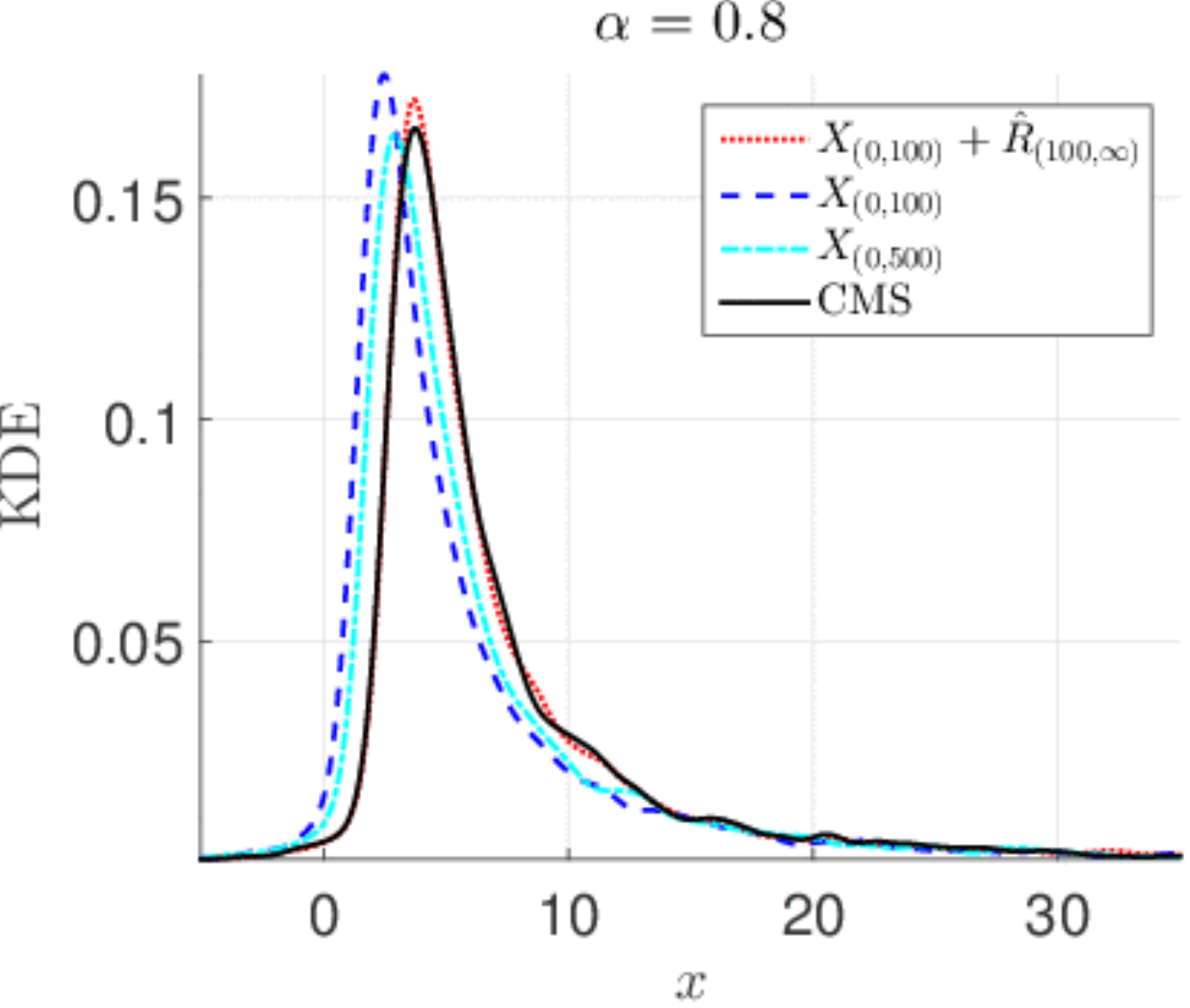}		
		\includegraphics[width=0.3\linewidth]{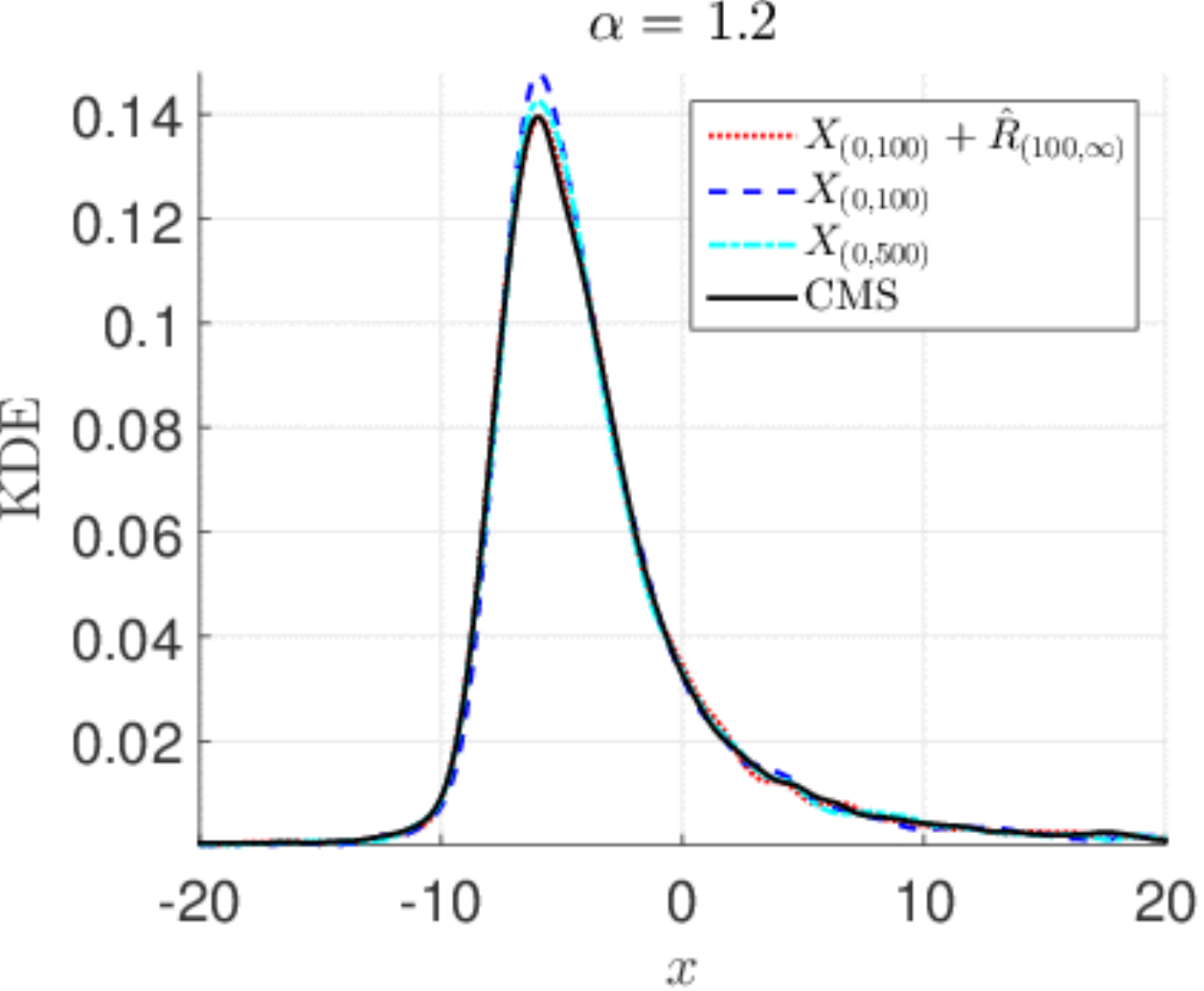}	
		\includegraphics[width=0.3\linewidth]{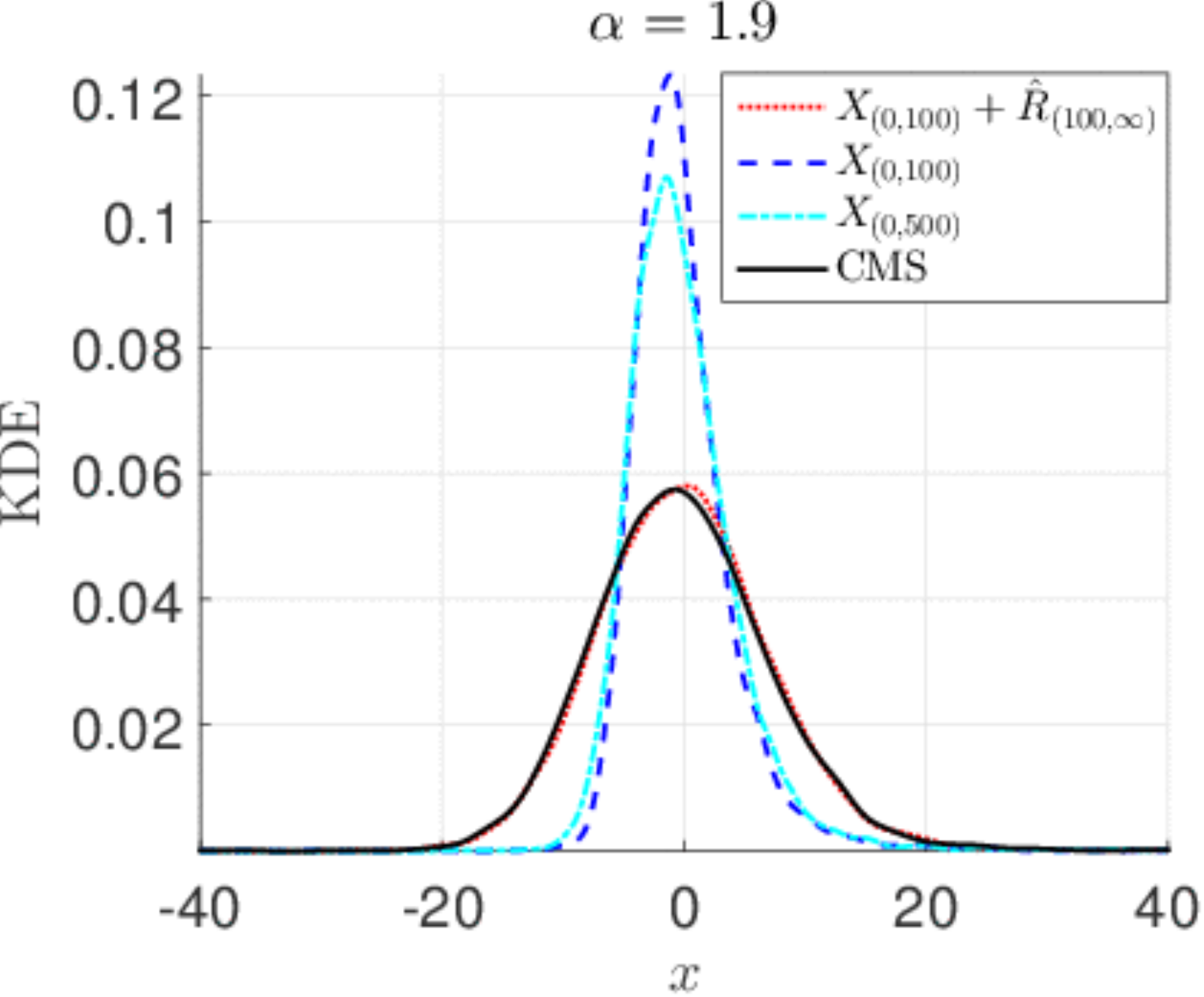}
	}
	\caption{Kernel density estimates based on $10^4$ samples
generated via three different methods: 
$X_{(0,c)}$ is the truncated PSR \eqref{eq:X_0c}, 
$X_{(0,c)} + \hat{R}_{(c, \infty)}$ is our proposed 
approximation for the PSR \eqref{proxy_PSR_split};
and CMS is the Chambers-Mallows-Stuck method. We compare results
for three values of the tail parameter $\alpha =0.8$ (left), 
$\alpha = 1.2$ (centre), $\alpha = 1.9$ (right), and for
two values of the truncation parameter, $c =100$ and $c = 500$.   
	}
	\label{fig:truncations_c}
\end{figure*}

In Figure \ref{fig:truncations_c} we compare 
kernel density estimates for the density of
$X\sim\mathcal{S}_\alpha(\sigma, \beta, 0)$ 
obtained by three different sampling methods:
$(i)$~$X_{(0,c)}$ is the obvious approximation
of $X$ by the truncated PSR \eqref{eq:X_0c};
$(ii)$~$X_{(0,c)} + \hat{R}_{(c, \infty)}$ is our proposed 
approximation for $X$ \eqref{proxy_PSR_split}; and
$(iii)$ `CMS' is the benchmark Chambers-Mallows-Stuck method for 
generating exact samples of stable RVs \cite{Chambers1976, Weron1996}. 
The results shown, indicate that
adding the Gaussian approximation of the residual to $X_{(0,c)}$ 
produces an approximate distribution that is much closer to the true stable law 
than that obtained by simple truncation of the PSR. 

Therefore, in view of Theorem~\ref{th:CLT}, existing
inference methods for the exact PSR can be used
for the approximation \eqref{proxy_PSR_split} and,
indeed, the
inference schemes in \cite{Lemke2012,  Lemke2014a, Lemke2015, Lemke2015a, 
Riabiz2017} mentioned in Section~\ref{subsect: param_inference} are based on 
\eqref{eq:proxy_gauss_model}.
Moreover, the
quality of this approximation is 
controlled directly by the truncation parameter $c$,
therefore,
it is important to have some
quantitative measure of the accuracy of the resulting approximation,
and also of the nature of its dependence on 
$c$ and on the parameters 
$\{\alpha,\, \beta, \, \sigma\}$. These issues
are addressed is the following sections.

\subsection{Choice of the truncation parameter $c$}

In order to quantify the approximation error in the
representation $X_{(0,c)}+\hat{R}_{(c,\infty)}$, 
and also in order to be able to choose appropriate
values for the truncation parameter $c$, the following
considerations should be kept in mind:
\begin{itemize}
	\item 
The distribution of $\hat{R}_{(c,\infty)}$ is closer
to that of $R_{(c,\infty)}$ when $c$ is large, according to 
Theorem \ref{th:CLT}. 
	\item 
	On the other hand,
	the computational complexity of the approximate conditionally Gaussian 
model  \eqref{eq:proxy_gauss_model} increases with $c$; in fact, 
the expected cardinality of the set of latent RVs $\{\Gamma_j <c \}$ needed to compute $m_{(0,c)}$ and $S^2_{(0,c)}$ is $O(c)$. 
	\item 
	Even when the distribution 
	of $R_{(c,\infty)}$ is far from Gaussian 
	(in particular, when $c$ is small),
	its contribution to the PSR might be relatively small 
	when compared to that of $X_{(0,c)}$. 	
\end{itemize} 

\noindent 
We will consider a choice of $c$ to be `good' if it makes the distribution 
of $\hat{X}$  close to that of $X$.
Quantifying the distance between $X$and $\hat{X}$ involves 
computing the distance between $R_{(c,\infty)}$ and $\hat{R}_{(c,\infty)}$,
so we proceed by first estimating how far the distribution
of the PSR residual is from the corresponding
Gaussian, for finite $c$.   

In view of the above discussion, our aim in the rest of the paper is to 
provide accurate bounds that can guide us in choosing
appropriate values of $c$, given the distribution parameters. 
The main tools that we employ in the derivation of such bounds
are based on classical Fourier-analytic techniques, 
summarized in the following section.

\newpage

\section{The Smoothing Lemma}
\label{section:smoothing}

As before, let $c \geq 0$ be the value of the
truncation parameter.
Suppose $S$ and $T$ are RVs with CFs $\phi_{S}(s)$ and $\phi_T(s)$, $s \in \mathbb{R}$, respectively, let
$F_{S}(x)$ and $F_T(x)$, $x\in\mathbb{R}$, be the 
corresponding CDFs, and assume that 
$\mathbb{E}[S]=\mathbb{E}[T]=0$. 
Furthermore, assume that $F_T(x)$ has derivative $p_T(x)$ such that $|p_T(x)|\leq m < \infty, \, \forall x \in \mathbb{R}$. We write,
\begin{IEEEeqnarray*}{rCl} 
	\Delta(S,T):=	\sup_{x \in \mathbb{R}}\left| F_{S}(x) - F_T(x)\right|,
	\label{eq:delta_cdf}
\end{IEEEeqnarray*}
\noindent
for the Kolmogorov distance between the distributions of $S$ and $T$.
Ess\'{e}en's smoothing lemma \cite[Lemma~XVI.4.2]{Feller1966}
states that, for any $\Theta>0$,
\begin{IEEEeqnarray}{rCl} 
	\Delta(S,T)
	&\leq &   
	\frac{1}{\pi}\int_{-\Theta}^\Theta\frac{ \left|\phi_{S}(s) -  \phi_T(s)\right|}{|s|}  \dif s + \frac{24 m}{\pi\Theta} :=I(S,T),
	\IEEEeqnarraynumspace
	\label{eq:int_berry_finite}
\end{IEEEeqnarray}
and letting $\Theta\to\infty$, we also have,
\begin{IEEEeqnarray}{rCl} 
	\Delta(S,T)
	&\leq &   
	\frac{1}{\pi}\int_{-\infty}^\infty \frac{\left|\phi_{S}(s) -  \phi_T(s)\right|}{|s|}  \dif s := \bar{I}(S,T), 	
	\IEEEeqnarraynumspace
	\label{eq:int_berry_esseen}
\end{IEEEeqnarray}

\noindent
where \eqref{eq:int_berry_esseen} is meaningful only if the improper 
integral converges.
Note that the integrals above have a removable singularity at $s = 0$, 
due to the zero-mean assumption on $S$ and $T$.

We will use the smoothing lemma to obtain bounds on the following:
$(i)$~The distance between the distribution of the 
standardized PSR residual and the standard normal,
\begin{IEEEeqnarray}{rCl} \label{eq:kolmog_residual}
	\Delta(Z_{(c, \infty)},Z)
	&:= &   
	\sup_{x \in \mathbb{R}}\left| F_{Z_{(c, \infty)}}(x) - F_Z(x)\right|, 
\end{IEEEeqnarray}
where $F_{Z_{(c,\infty)}}(x)$ and $F_Z(x)$ denote the CDF of ${Z_{(c,\infty)}}$ 
and
the standard normal CDF, respectively. 
$(ii)$~The distance between 
the distribution of the approximately stable RV
$\hat{X}$ and the exact stable law,
\begin{IEEEeqnarray}{rCl} \label{eq:kolmog_stable}
	\Delta(X, \hat{X})
	&:= &   
	\sup_{x \in \mathbb{R}}\left| F_{X}(x) - F_{\hat{X}}(x)\right|, 
\end{IEEEeqnarray}
where $F_X(x)$ and $F_{\hat{X}}(x)$ are the CDFs of $X$  
and $\hat{X}$, respectively. And~$(iii)$ The distance between
the truncated PSR and the full PSR,
\begin{IEEEeqnarray}{rCl} \label{eq:kolmog_truncated}
	\Delta(X, X_{(0,c)})
	&:= &   
	\sup_{x \in \mathbb{R}}\left| F_{X}(x) - F_{X_{(0,c)}}(x)\right|, 
\end{IEEEeqnarray}
where $F_{X_{(0,c)}}(x)$ is the CDF of $X_{(0,c)}$~\eqref{eq:X_0c}. 

In order to apply the smoothing lemma to bound \eqref {eq:kolmog_residual}, 
\eqref{eq:kolmog_stable} and \eqref{eq:kolmog_truncated}, 
we need explicit expressions for the various CFs of interest. 
These are derived in the following section.

\newpage

\section{Characteristic Functions}
\label{section:cf:equations}

Recall that 
$X=X_{(0,c)}+R_{(c,\infty)}\sim\mathcal{S}_\alpha(\sigma,\beta,0)$ 
has CF 
$\phi_X(s)$ given in \eqref{eq:cf}, and that we approximate 
it by 
$\hat{X}=X_{(0,c)}+\hat{R}_{(c,\infty)}$ 
as in~\eqref{proxy_PSR_split}.  
Since $R_{(c,\infty)}$ and $\hat{R}_{(c,\infty)}$ 
are independent of $X_{(0,c)}$, we have,
\begin{IEEEeqnarray}{rCl} \phi_X(s)&=&\phi_{X_{(0,c)}}(s)\phi_{R_{(c,\infty)}}(s), \label{eq:true_stable_CF}
	\\
	\phi_{\hat{X}}(s)
	&=&\phi_{X_{(0,c)}}(s)\phi_{\hat{R}_{(c,\infty)}}(s). 
	\label{eq:proxy_stable_cf}
	\IEEEeqnarraynumspace
\end{IEEEeqnarray}
Also, the CFs for the true and approximated residuals can be expressed
in terms of the CFs of their normalised counterparts 
by a simple change of variable,
\begin{IEEEeqnarray}{rCl}
	\phi_{R_{(c,\infty)}}(s) &= \phi_{Z_{(c,\infty)}}(S_{(c,\infty)}s) \exp\left(i s m_{(c,\infty)}\right),
	\label{eq:res:unnormalized:exact}	
	\\
	\phi_{\hat{R}_{(c,\infty)}}(s) &= \phi_Z(S_{(c,\infty)}s) \exp\left(i s m_{(c,\infty)}\right), 
	\label{eq:res:unnormalized:proxy}
	\IEEEeqnarraynumspace
\end{IEEEeqnarray} 
where $m_{(c,\infty)}$ and $S_{(c,\infty)}$ are given in~\eqref{eq:m_c_infty} 
and~\eqref{eq:S_c_infty}, and $\phi_Z(\cdot)$ is the CF of the standard normal 
distribution \eqref{eq:Z:lim_dc}. 

Therefore, in order to use the smoothing lemma for 
\eqref{eq:kolmog_residual} and \eqref{eq:kolmog_stable},
we need to obtain explicit 
expressions for $\phi_{Z_{(c, \infty)}}(s)$, $\phi_{R_{(c, \infty)}}(s)$,
and $\phi_{X_{(0,c)}}(s)$. These are derived in the following
two subsections, in the 
case $W_1 \sim \mathcal{N}(\mu_W, \sigma^2_W)$.
The proofs are given in the Appendix.
For easy reference, the results
are summarized in Table~\ref{tab:CFS} at the end of this section.

\subsection{CF expressions when $W_1 \sim \mathcal{N}(\mu_W, \sigma^2_W)$} 
\label{section_CF_gaussian_case}

When $W_1$ is Gaussian, the following lemma shows
that it is possible to write the CF of the 
PSR residual in terms of an infinite series.

\medskip

\begin{lemma}\label{lemma:cf_series_R_cd}
	Let  $Z_{(c, \infty)}$ be defined as in  \eqref{eq:CLT}, and let 
$W_1 \sim \mathcal{N}(\mu_W, \sigma^2_W)$. Then, for $s\in\RL$ and $c>0$,
	\begin{IEEEeqnarray}{rCl} 
		\phi_{Z_{(c,\infty)}}(s) = 
		\exp\left(
		- \frac{s^2}{2} +\sum_{k=3}^\infty \bar{z}_k s^k\right),
		\label{eq:cf_Z_cinfty_series}
		\IEEEeqnarraynumspace
	\end{IEEEeqnarray}
	where,
	\begin{IEEEeqnarray}{
			rCl} 
		\bar{z}_k &:=& 
		\frac{i^k}{k!}\frac{\mathbb{E}[W_1^k]\frac{\alpha}{k - \alpha}}
		{\left(\mathbb{E}[W_1^2]\frac{\alpha}{2- \alpha}\right)^{k/2}} c^{1 - k /2}, \quad k \geq 3.
		\label{eq:z_bar_k}
		\IEEEeqnarraynumspace
	\end{IEEEeqnarray}
\end{lemma}

\medskip

An examination of the proof in Appendix \ref{app:series_cf} shows that, in the process of establishing
the lemma, we also obtained expressions for the CF of the unnormalized 
residuals $R_{(c, d)}$ and  $R_{(c, \infty)}$; these are shown in 
Table~\ref{tab:CFS}. We also note that these results hold not only
in the case
$W_1 \sim \mathcal{N}(\mu_W, \sigma^2_W)$, but also more generally 
for any distribution on $W_1$ that satisfies condition 
\eqref{eq:suff_cond_analiticity}.

Alternatively, performing direct computations when $W_1$ is normally 
distributed, we obtain the following integral expressions for the CF 
of the residual. 

\medskip

\begin{lemma}\label{lemma:cf_integral_R_cd}
	Let  $R_{(c, \infty)}$ be defined as in  \eqref{eq:true_residual}, and let $W_1 \sim \mathcal{N}(\mu_W, \sigma^2_W)$. Then, for $s\in\RL$,
	\begin{IEEEeqnarray}{rCl} 
		\log
		(\phi_{R_{(c,\infty)}} (s))  
		& = &
		{\alpha}\int_{0}^{c^{-1/\alpha}}\left( e^{(is t\mu_W-\sigma_W^2s^2t^2/2)} - 1 -is\mu_Wt\right) t^{-\alpha-1}\dif t	
		- is \mu_W \frac{\alpha}{\alpha -1}c^{\frac{\alpha-1}{\alpha}}. 
		\label{eq:cf_res_infty_integral}
	\end{IEEEeqnarray}
\end{lemma}

\medskip

Lemma~\ref{lemma:cf_integral_R_cd} is established in 
Appendix~\ref{app:integral_cf}, where we also derive 
analogous integral expressions for for the 
CFs of $X_{(0,c)}$ and of $X$.


\subsection{CF expressions when $W_1 \sim \mathcal{N}(0, \sigma^2_W)$} 
\label{section_PSR_cf_symm}

Next we obtain a more explicit expression for the CF
of the normalized residual $Z_{(c,\infty)}$ in the case
when the mean $\mu_W=0$. This expression was first derived,
by summing the series \eqref{eq:cf_Z_cinfty_series},
in~\cite{Riabiz2017}. A different proof, based on
the integral representation in Lemma~\ref{lemma:cf_integral_R_cd},
is given in Appendix~\ref{appendix:proofs_cf_symm}. 

\medskip

\begin{lemma}\label{lemma:cf_res_symm}
	Suppose $W_1\sim\mathcal{N}(0, \sigma^2_W)$,
	and denote,
	\begin{IEEEeqnarray}{c}
		a:= \frac{\alpha}{2}, \;\; 		\eta:=\frac{1-a}{a}, \;\; 	\w:=\frac{\eta s^2}{2c}, \;\; u := wS^2_{(c, \infty)},
		\IEEEeqnarraynumspace
		\label{eq:change_var}
	\end{IEEEeqnarray}
	for $\alpha\in(0,2), \alpha \neq 1$,
	where $S^2_{(c,\infty)}$ is given 
	as in (\ref{eq:S_c_infty}) with $E[W_1^2]=\sigma_W^2$.
	Then,
	\begin{IEEEeqnarray}{rCl} 
		\phi_{Z_{(c,\infty)}}(s) &=& 	\psi_{Z_{(c,\infty)}}(\w) 
		\nonumber \\
		&=&\exp 
		\left(c(1-e^{-\w}-\w^a\gamma\left(1-a,\w\right))
		\right),	
		\IEEEeqnarraynumspace
		\label{eq:asympt_cf_res_symm}
	\end{IEEEeqnarray}
	where $\gamma(s, x)$ is the lower incomplete gamma function \eqref{eq:lower_inc_gamma}.
	Moreover, 
	\begin{IEEEeqnarray}{rCl} 
		\phi_{X_{(0,c)}}(s) &=&	\omega_{{X_{(0,c)}}}(u) \nonumber \\
		& = & \exp(-c(1-e^{-u}+u^a\Gamma(1-a,u))),		
		\IEEEeqnarraynumspace
		\label{eq:cf_X0c_symm}
	\end{IEEEeqnarray}
	where $\Gamma(s,x) $ is the upper
	incomplete gamma function 	\eqref{eq:upper_inc_gamma}.
\end{lemma}

\medskip

Observe that, using the change of variables \eqref{eq:change_var}, the CF of the standard normal RV $Z$ \eqref{eq:Z:lim_dc} can be written,
\begin{IEEEeqnarray}{c} 
	\phi_Z(s) = \psi_Z(w) = \exp(-cw/\eta).
	\label{eq:cf_std_gauss_w}
\end{IEEEeqnarray}	 
\noindent 
Hence, 
using \eqref{eq:res:unnormalized:proxy}, when $\mu_W=0$ we have
\begin{IEEEeqnarray}{rCl} 
	\phi_{\hat{R}_{(c, \infty)}}(s) = \omega_{\hat{R}_{(c, \infty)}}(u) = \exp(-cu/\eta). 
	\label{eq:residual_proxy_u} 
\end{IEEEeqnarray}	
Then, as a consequence of Lemma~\ref{lemma:cf_res_symm} and 
equation \eqref{eq:proxy_stable_cf},
it follows that, when $\mu_W =0$, 
the CF 
$\phi_{\hat{X}}(s) = \omega_{\hat{X}}(u)$
of $\hat{X}$, the approximated stable distribution,
satisfies,
\begin{IEEEeqnarray*}{c} 
	\log \omega_{\hat{X}}(u) 
	=-c(1 - e^{-u} + u^a \Gamma(1-a,u) + u/\eta).
	\IEEEeqnarraynumspace
	\label{eq:proxy_stable_symmetric}
\end{IEEEeqnarray*}

From now on and for the rest of the paper we restrict attention
to the case $\mu_W =0$ of the symmetric stable law, for which we 
have the above 
closed-form expressions for $\phi_{Z_{(c, \infty)}}(s)$.

{\small{
		\begin{table*}[ht!]
			\caption{Summary of the logarithms of the CF expressions derived and used. Recall that $\alpha \in(0,2), \, \alpha \neq 1$, and $c >0$. 
			}
			\centering
		\bgroup
		\def\arraystretch{1.7}
		
			\begin{tabular}{|c || c ||c|| c || c| }
				\hline
				Distribution  $W_j$ & Skewness & RV & $\log(\phi(s))$ or $\log(\psi(w))$ or $\log(\omega(u))$, with $w$ and $u$ as in \eqref{eq:change_var}& Equations \\
				
				\hline
				\hline
				
				\multirow{6}{*}{Satisfying 
\eqref{eq:condition_W}, {\eqref{eq:suff_cond_analiticity}}}  & \multirow{6}{*}{$\beta \in [-1,1]$} & $X$ & 	$	-\sigma^{\alpha}{|s|}^{\alpha}\left\lbrace 1-i\beta\sign(s)\tan\frac{\pi\alpha}{2}\right\rbrace  $ & \eqref{eq:cf}, \eqref{eq:mapping_sigma}-\eqref{eq:mapping_beta}\\
				\cline{3-5}

				\cline{3-5}
				& & $R_{(c,d)}$ & $ ism_{(c, d)}
				- s^2S^2_{(c, d)}/2 +
				\sum_{k=3}^\infty r_k s^k$ &  \eqref{eq:S_cd},  \eqref{eq:cf_rcd_series},   \eqref{eq:r_k}\\
				
				\cline{3-5}
				& & $R_{(c,\infty)}$ & $ ism_{(c, \infty)}
				- s^2S^2_{(c, \infty)}/{2} +\sum_{k=3}^\infty \bar{r}_k s^k$& 
				\eqref{eq:cf_R_cd_series}, \eqref{eq:bar_r1}, \eqref{eq:S_c_infty}\\
				
				\cline{3-5}
				& & $Z_{(c,\infty)}$ & $- {s^2}/{2} + \sum_{k=3}^\infty \overline{z}_k s^k$ &  \eqref{eq:cf_Z_cinfty_series}   \eqref{eq:z_bar_k}\\

				\cline{3-5}
				& & ${Z}$ & $- {s^2}/{2} \quad $ or $\quad -cw/\eta$& \eqref{eq:Z:lim_dc}, \eqref{eq:cf_std_gauss_w} \\

				\cline{3-5}
				& & $\hat{R}_{(c,\infty)}$ & $-s^2 S^2_{(c,\infty)}/2 + ism_{(c, \infty)}$& \eqref{eq:m_c_infty}, \eqref{eq:S_c_infty}, \eqref{eq:res:unnormalized:proxy}\\

				\hline
				\hline
				\multirow{5}{*}{$ \mathcal{N}(\mu_W, \sigma^2_W)$}
				&\multirow{5}{*}{$\beta \in [-1,1]$} &  \multirow{2}{*}{$R_{(c,d)}$ } & 
				$c-d
				+    
				{\alpha}\int_{d^{-1/\alpha}}^{c^{-1/\alpha}}e^{(is t\mu_W-\sigma_W^2s^2t^2/2)}t^{-\alpha-1}\dif t - isB$
				& \eqref{eq:cf_Rcd_integral_1}, \eqref{eq:B}
				\\
				\cline{4-5}
				& & & ${\alpha}\int_{d^{-1/\alpha}}^{c^{-1/\alpha}}\left( e^{(is t\mu_W-\sigma_W^2s^2t^2/2)} - 1 -is\mu_Wt\right) t^{-\alpha-1}\dif t	
				- is \mu_W \frac{\alpha}{\alpha -1}c^{\frac{\alpha-1}{\alpha}}$
				&  \eqref{eq:cf_Rcd_integral_2}
				\\
				\cline{3-5}
				& & $R_{(c,\infty)}$ & 
				${\alpha}\int_{0}^{c^{-1/\alpha}}\left( e^{(is t\mu_W-\sigma_W^2s^2t^2/2)} - 1 -is\mu_Wt\right) t^{-\alpha-1}\dif t	 - is \mu_W \frac{\alpha}{\alpha -1}c^{\frac{\alpha-1}{\alpha}}$
				& \eqref{eq:cf_res_infty_integral}
				\\
				\cline{3-5}
				& & $X_{(0,c)}$   & 
				$-c
				+    
				\int_{c^{-1/\alpha}}^{\infty}e^{(is t\mu_W-\sigma_W^2s^2t^2/2)}t^{-\alpha-1}\dif t $
				& \eqref{eq:phi_X_0c_gauss}
				\\
				\cline{3-5}

				& & $X$ 
				& ${\alpha}\int_{0}^{\infty}\left( e^{(is t\mu_W-\sigma_W^2s^2t^2/2)} - 1 \right) t^{-\alpha-1}\dif t	 $
				& \eqref{eq:phi_X_0inf_gauss}
				\\
				\hline
				\hline
				\multirow{5}{*}{$ \mathcal{N}(0, \sigma^2_W)$}
				&\multirow{5}{*}{$\beta =0$}

				& $R_{(c,\infty)}$  & $c\big(1-\exp(-u)-u^a\gamma\left(1-a,u\right)\big)$
				& 
				\eqref{eq:phi_R_c_inf_symm}

				\\
				\cline{3-5}
				& & $Z_{(c,\infty)}$  & $c\big(1-\exp(-\w)-\w^a\gamma\left(1-a,\w\right)\big)$
				& \eqref{eq:asympt_cf_res_symm}
				
				\\
				
				\cline{3-5}

				& &  $X_{(0,c)}$  & 
				$-c(1-\exp(-u)+u^a\Gamma(1-a,u))$
				& \eqref{eq:cf_X0c_symm}

				\\

				\cline{3-5}
				& & $\hat{R}_{(c, \infty)}$
				& 
				$-cu/\eta$
				& 
				\eqref{eq:residual_proxy_u}

				\\
				\cline{3-5}
				& & $\hat{X}$
				& 
				$-c(1 - \exp{(-u)} + u^a \Gamma(1-a,u) + u/\eta)$
				& 
				\eqref{eq:proxy_stable_symmetric}
				\\	
				\hline
			\end{tabular}
		\egroup		
			\label{tab:CFS}
		\end{table*}
}}

\newpage

\section{Nonasymptotic Gaussian Bounds for the PSR Residual}
\label{section:residual_convergence}

\subsection{Nonasymptotic bound of order $O(1/c)$}
In this section we apply the smoothing lemma of Section~\ref{section:smoothing}
to derive 
explicit bounds on the distance
$\Delta(Z_{(c,\infty)}, Z)$, defined in~\eqref{eq:kolmog_residual}.
When $\mu_W =0$, the closed-form expression 
in~\eqref{eq:asympt_cf_res_symm} for $\phi_{Z_{(c,\infty)}}(s)$
can be used to further bound above the term
$\bar{I}(Z_{(c,\infty)},Z)$ in~(\ref{eq:int_berry_esseen}).
The resulting bounds, 
first presented in \cite{Riabiz2018}, 
are stated in the following theorem and proved
in Appendix \ref{appendix:proofs_th_lin}.

\medskip

\begin{theorem}
	\label{thm:linear}
	Let $W_j \sim \mathcal{N}(0, \sigma^2_W)$ and let $\Delta(Z_{(c,\infty)}, Z)$ be the Kolmogorov distance between $Z_{(c,\infty)}$ and $Z$, as in \eqref{eq:kolmog_residual}. Let $a = a(\alpha)$ and $\eta = \eta(\alpha)$ as in~\eqref{eq:change_var}, and define,
	\begin{IEEEeqnarray*}{rCl} 
		g(w) & := &  1-e^{-\w}-\w^a\gamma\left(1-a,\w\right), \qquad w\geq0.
	\end{IEEEeqnarray*} 
	Let 
	$\gamma(s,x)$ and $\Gamma(s,x)$ 
	be the lower and upper incomplete gamma functions,
	\eqref{eq:lower_inc_gamma} and \eqref{eq:upper_inc_gamma}
	respectively, and write,
	\begin{IEEEeqnarray}{rCl} 
		\bar{\gamma}(a) & := & \gamma(1-a,1),
		\label{eq:bar_gamma_a}	
		\\ 
		\barg &:=& g(1),
		\label{eq:bar_g}
		\\
		K(a) &:=& \frac{1}{\pi}\left[\frac{a}{2(2-a)} + \frac{1}{\eta^2} \right]. 
		\label{eq:K_c_a} 	
	\end{IEEEeqnarray} 	
	Then, for any $\c>1$, $\Delta(Z_{(c,\infty)}, Z)$  
	is bounded above by:

	\begin{IEEEeqnarray*}{rCl} 
		B_1(c, \alpha)
		:= &   
		\frac{cK(a)}{(c-1)} &
		\left[
		\frac{1}{(\c-1) \barg^2} +
		\left(\frac{1}{\barg}-\frac{1}{(\c-1) \barg^2}\right)
		\exp\left((\c-1)\barg\right)+
		\right. \nonumber
		\\
		& & 
		\left. +
		\frac{
			(\c-1)\exp\{(\c-1)
			(1-e^{-1})
			\}
		}
		{a[(\c-1)\bar{\gamma}(a)]^{2/a}}
		\Gamma\left({2}/{a}, (\c-1)\bar{\gamma}(a)\right) 
		\right].
		\IEEEeqnarraynumspace
	\end{IEEEeqnarray*}

	\noindent
\end{theorem}

{
\begin{remark}
	Observe that the upper bound on the Kolmogorov distance of the PSR residual from its Gaussian approximation
	coincides with the upper bound on the the distance of the PSR standardized residual from the standard Gaussian,
		\begin{IEEEeqnarray*}{rCl} 
			\bar{I}(R_{(c, \infty)}, \hat{R}_{(c, \infty)}) = \bar{I}(Z_{(c,\infty)}, Z).
		\end{IEEEeqnarray*} 
	In fact,
the Kolmogorov distance itself is invariant under monotone transformations, thus in particular for translation and scaling.
\end{remark}
}

\begin{remark}
	Even though the PSR split was introduced in \eqref{exact_PSR_split} 
for any $c \geq 0$, Theorem~\ref{thm:linear} and the results below
	hold only for $c>1$.
This is not a significant limitation because, in practice, 
we are indeed interested  
	in scenarios where there are at least a few terms in $X_{(0,c)}$.
\end{remark}

The following corollary is a simple consequence of Theorem~\ref{thm:linear};
its proof is given in Appendix~\ref{appendix:proof:corollary:lin} 

\medskip

\begin{corollary}
	\label{cor_alternative}
	Under the assumptions of Theorem \ref{thm:linear}, 
	as $c\to\infty$,
	$$B_1(c,\alpha)\sim 
	\frac{K(a)}{\bar{g}^2}\Big(\frac{1}{c-1}\Big),$$
so that,
	$$\Delta(Z_{(c,\infty)}, Z)
	=O\Big(\frac{1}{c}\Big).$$
\end{corollary}

\medskip

\begin{figure*}[ht!] 
	\centering
	\includegraphics[width=0.31\linewidth]{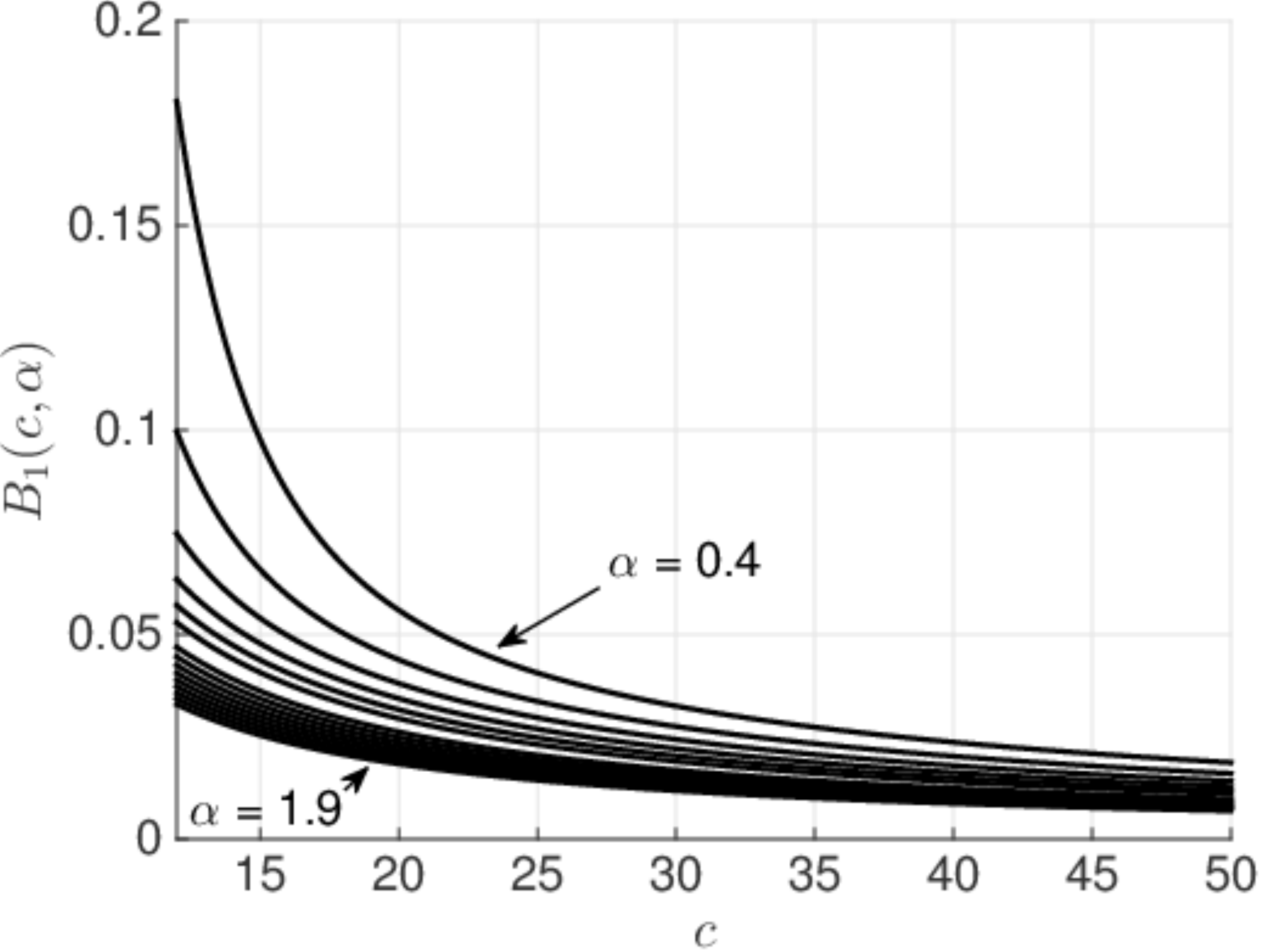} 
	\includegraphics[width=0.31\linewidth]{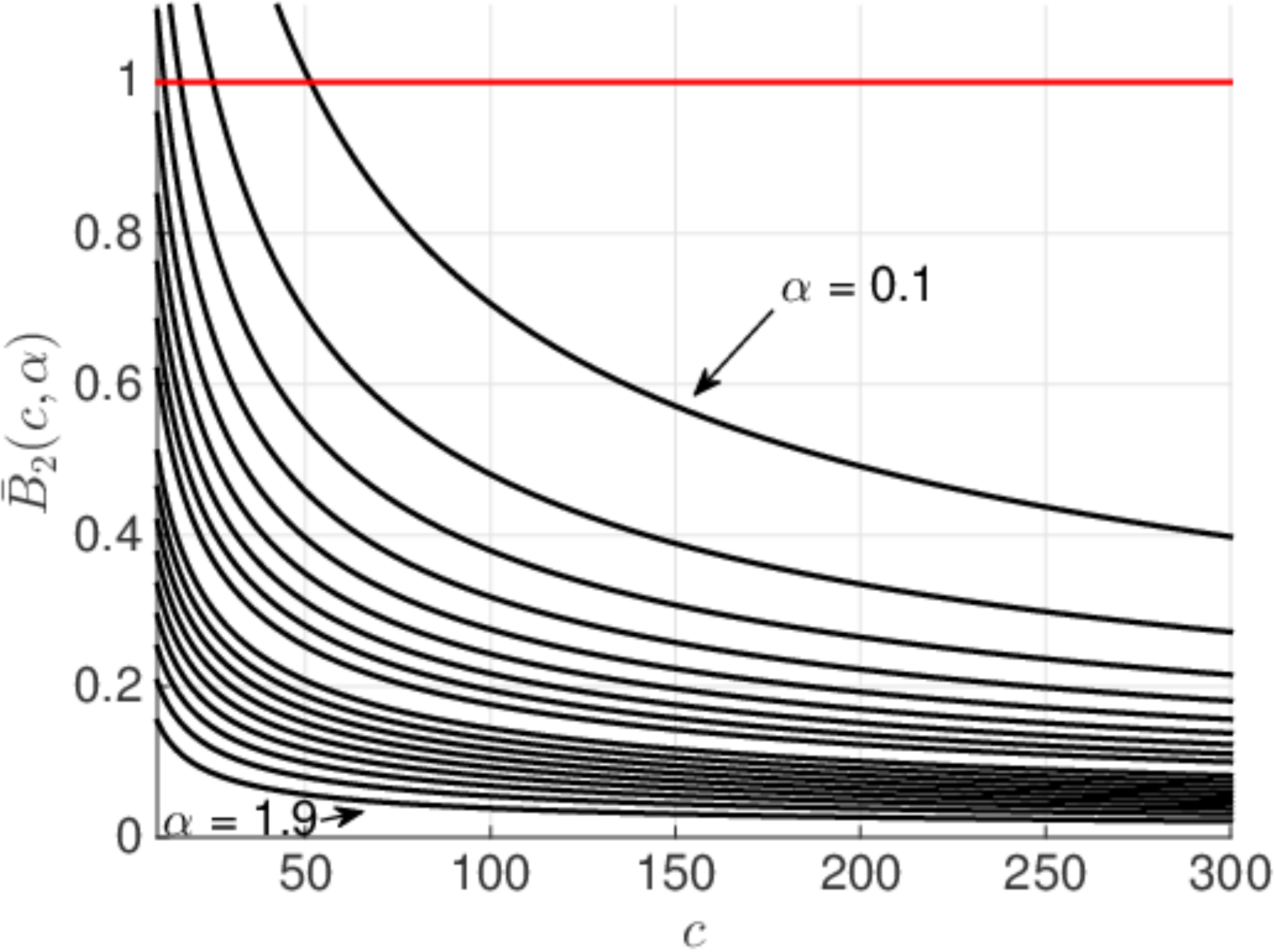} 
	\includegraphics[width=0.31\linewidth]{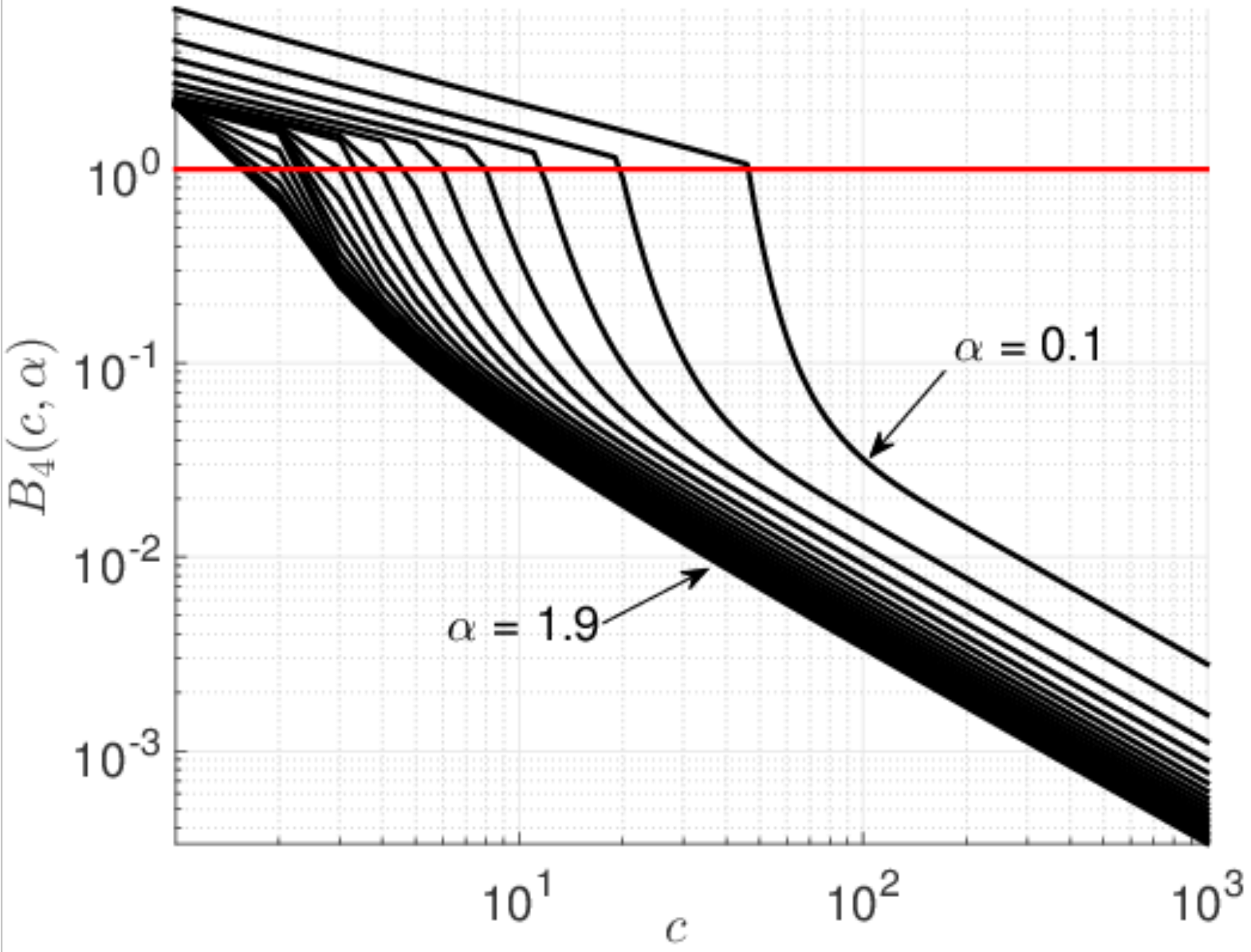} 
	\caption{
		{\small{
				Bounds on $\Delta(Z_{(c, \infty)}, Z)$. Left: each curve represents the values of the bound $B_1(c,\alpha)$
				for  $\alpha=0.4,\,0.9,\ldots,\,1.9$, $\alpha \neq 1$,
				plotted against $12\leq \c\leq 50$.
				Centre: each curve represents the values of the bound
				$\bar{B}_2(c,\alpha)$ 
				for $\alpha=0.1,\,0.2,\ldots,\,1.9$, $\alpha \neq 1$,
				plotted against $10\leq \c\leq 300$.
				Right: each curve represents $B_4(c, \alpha)$, for $\alpha=0.1,\, 0.2,\, \ldots,1.9$, $\alpha \neq 1$,
				plotted against $10\leq \c\leq 1000$.
[The red horizontal line at 1 shows the 
maximum possible value of the Kolmogorov distance.]
	}} }
	\label{f:1}
\end{figure*}

For values of $\alpha$ greater than $0.4$,
$B_1(c, \alpha)$ gives very good bounds, as shown on the left-hand side 
of Figure~\ref{f:1}. But for $\alpha$ below $0.4$, 
the results 
deteriorate significantly;
for example, for $\alpha=0.2$, $B_1(c, \alpha)$ is below~1 
(the maximum  possible Kolmogorov distance)
only for $\c>115$.

\subsection{Nonasymptotic bound of order $O(1/\sqrt{c})$}

The following result, obtained by bounding 
$I(Z_{(c, \infty)}, Z)$ in \eqref{eq:int_berry_finite}, gives an $O(1/\sqrt{\c})$
bound which is, of course, asymptotically inferior
to that in Theorem~\ref{thm:linear}, but which gives
sharper results for small~$\c$ and $\alpha<0.4$.

\medskip

\begin{theorem}
	\label{thm:root}
	Under the assumptions and in the notation 
	of
	Theorem \ref{thm:linear}, 
	for any $\delta\in(0,2)$,
	\begin{IEEEeqnarray*}{c} 
	\Delta(Z_{(c,\infty)}, Z)\leq
		B_2(c, \alpha, \delta):=\frac{9.6\sqrt{\eta}}{\pi\sqrt{2(2-\delta)\c}} + B_3(c, \alpha, \delta),
	\end{IEEEeqnarray*}
	where $B_3(c, \alpha, \delta)$ is the following $O(1/\c)$ term:
	\begin{IEEEeqnarray*}{rCl} 	
		B_3(c, \alpha, \delta)  :=
		\frac{K(a)}{c}
		\left(\frac{c (2-\delta) }{(c-1) g(2-\delta)}\right)^2 \times 
		\left\{1-\left[1 - g(2-\delta) (c-1)\right]
		\exp\left(g(2-\delta)(c-1)\right)\right\}.
	\end{IEEEeqnarray*}
\end{theorem}

\medskip

The proof is given in Appendix \ref{appendix_proof_th_sqrt}. 
Numerically minimizing the bound $B_2(c, \alpha, \delta)$ 
over~$\delta$ yields  $\bar{B}_2(c, \alpha)$, shown in the central part of Figure~\ref{f:1}.

\subsection{Combined bound and comparison with numerical results}

Finally, we combine the results of Theorems~\ref{thm:linear} 
and~\ref{thm:root}, to obtain useful bounds essentially for all 
values of $\alpha \in (0,2)$, $\alpha \neq 1$, and  {$\c > 50$}:
\begin{IEEEeqnarray*}{c} 
	\Delta(Z_{(c,\infty)}, Z)\leq
	B_4(c, \alpha):= \min \left\lbrace B_1(c,\alpha), \bar{B}_2(c, \alpha)\right\rbrace.
\end{IEEEeqnarray*}
The resulting numerical bound 
is shown on the right-hand side of 
Figure \ref{f:1} (on a log-log scale). 

\begin{figure}[ht!] 
	\centering
	\includegraphics[width=0.4\linewidth]{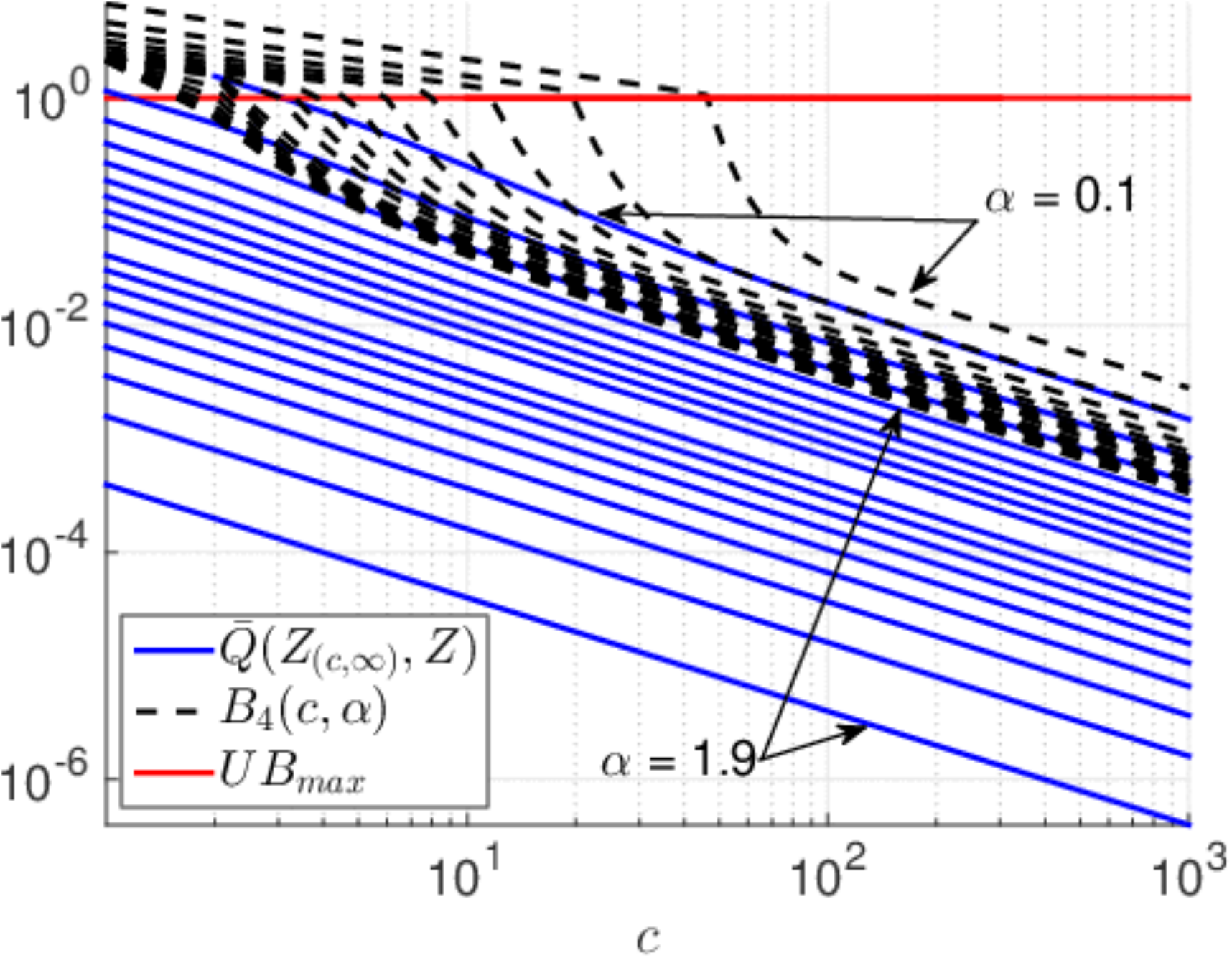}
	\caption{
		Bounds on $\Delta(Z_{(c, \infty)}, Z)$: comparison between  $\bar{Q}(Z_{(c, \infty)}, Z)$ and $B_4(c, \alpha)$, for $\alpha=0.1,\,0.2,\ldots,\,1.9$, $\alpha \neq 1$,
		plotted for $3\leq \c\leq 300$. 
[The red horizontal line at 1 shows the 
maximum possible value of the Kolmogorov distance.]
	}
	\label{f:num_vs_th_integral_residual}
\end{figure}

Figure~\ref{f:num_vs_th_integral_residual}  
shows a comparison between  the theoretical 
bound $B_4(c,\alpha)$ and 
the numerical estimate
$\bar{Q}(Z_{(c,\infty)},Z)$ 
of $\bar{I}(Z_{(c, \infty)}, Z)$,
produced through 
the Matlab routine \texttt{quadgk},
which implements the Gauss-Kronrod method;
see \cite{Riabiz2017b} and \cite{Riabiz2018}
for more details.
This method also produces an approximate upper bound on 
the absolute error $|\bar{I}(Z_{(c, \infty)}, Z)-\bar{Q}(Z_{(c, \infty)}, Z)|$,
which can be used to construct approximate error bands. 
But for $c\geq 3$ these are
negligibly small,
so we do not show them here. 
Observe that $B_4(c,\alpha)$ appears to have the exact same asymptotic 
rate as~$\bar{Q}(Z_{(c, \infty)}, Z)$. 

\clearpage

\section{Bounds on $\alpha$-stable Distribution Approximations}
\label{section:residual_contribution}

Here we develop bounds on the distances 
$\Delta(X, \hat{X})$ 
and $\Delta(X, X_{(0,c)})$ 
defined in~\eqref{eq:kolmog_stable} 
and~\eqref{eq:kolmog_truncated}, respectively.
In terms of inference, ultimately, 
it is these Kolmogorov distances that we wish
to make ``small'' by appropriately choosing the 
value of the parameter~$c$.


\subsection{Nonasymptotic bound on $\Delta(X, \hat{X})$}	
The bound on $\Delta(X, \hat{X})$ stated next,
is established by using the smoothing lemma \eqref{eq:int_berry_esseen} 
and the bound in  Theorem~\ref{thm:linear}. Its proof is given
in Appendix~\ref{app:proof_hat_X_symm}.

\medskip

\begin{theorem}
\label{theorem:X_hat_symm} 
Let 
	$\Delta(X, \hat{X})$  be the Kolmogorov distance between $X$ and $\hat{X}$, as in \eqref{eq:kolmog_stable},
	under the same assumptions and in the same notation as 
	Theorem \ref{thm:linear}. 
	Let $N\geq 1$, and introduce $N$ arbitrary
	points $u_i$ on $[0,1]$,
	 $$0 =: u_0 < u_1 < \cdots < u_N := 1,$$ 	
	together with the corresponding
	values of the logarithm 
	of the CF $\omega_{X_{(0,c)}}(u)$ defined
	in \eqref{eq:cf_X0c_symm},
$$f_0: = 0, \qquad f_i := \log(\omega_{X_{(0,c)}}(u_i)), \;\;\;
	i=1,2,\ldots,N.$$
	Also let,
	\begin{IEEEeqnarray}{rCl} 
		m_i &:=& \frac{f_{i+1} - f_i}{u_{i+1} - u_i},
		\qquad
		q_i := - m_i u_i + f_i. 	
		\label{eq:q_m}
	\end{IEEEeqnarray}	
	\noindent
	Then, for any $\c>1$, $\Delta(X,\hat{X})$  is bounded above by,
	\begin{IEEEeqnarray*}{C} 
		B_5(c, \alpha, N):= cK(a)
		\times
		\left\lbrace 
		\sum_{i=0}^{N-1}\frac{e^{q_i}}{\tilde m_i}\left[ e^{\tilde m_i u_{i+1} }\left( u_{i+1} - \frac{1}{\tilde m_i}\right)  
		- 
		e^{\tilde m_i u_{i} }\left( u_{i} - \frac{1}{\tilde m_i}\right)
		\right]
		+ 	
		\frac{e^{\tilde k_{(1,\infty)} }}
		{a \big(\tilde l_{(1,\infty)}\big)^{2/a}} \Gamma\left(\frac{2}{a}, \, \tilde l_{(1,\infty)} \right)
		\right\rbrace,
		\IEEEeqnarraynumspace
	\end{IEEEeqnarray*}
	\noindent
	where,
	\begin{IEEEeqnarray}{rCl} 
		\tilde m_i &:=& m_i +(c-1)\barg, \nonumber \\
		k_{(1,\infty)} &:=&  -c((1-\exp(-1))+\Gamma(1-a,1)) <0, \label{eq:k_1_infty}\\
		\tilde k_{(1,\infty)} &:=& k_{(1,\infty)} - (c-1)(e^{-1} -1), \nonumber \\
		\tilde l_{(1,\infty)} &:=& (c-1)\bar{\gamma}(a). \nonumber
	\end{IEEEeqnarray}
	and $\barg$ as in \eqref{eq:bar_g} and $\bar{\gamma}(a)$ as in \eqref{eq:bar_gamma_a}. 
\end{theorem}

\medskip

\begin{figure*}[ht!] 
	\centerline{
		\includegraphics[width=0.3\linewidth]{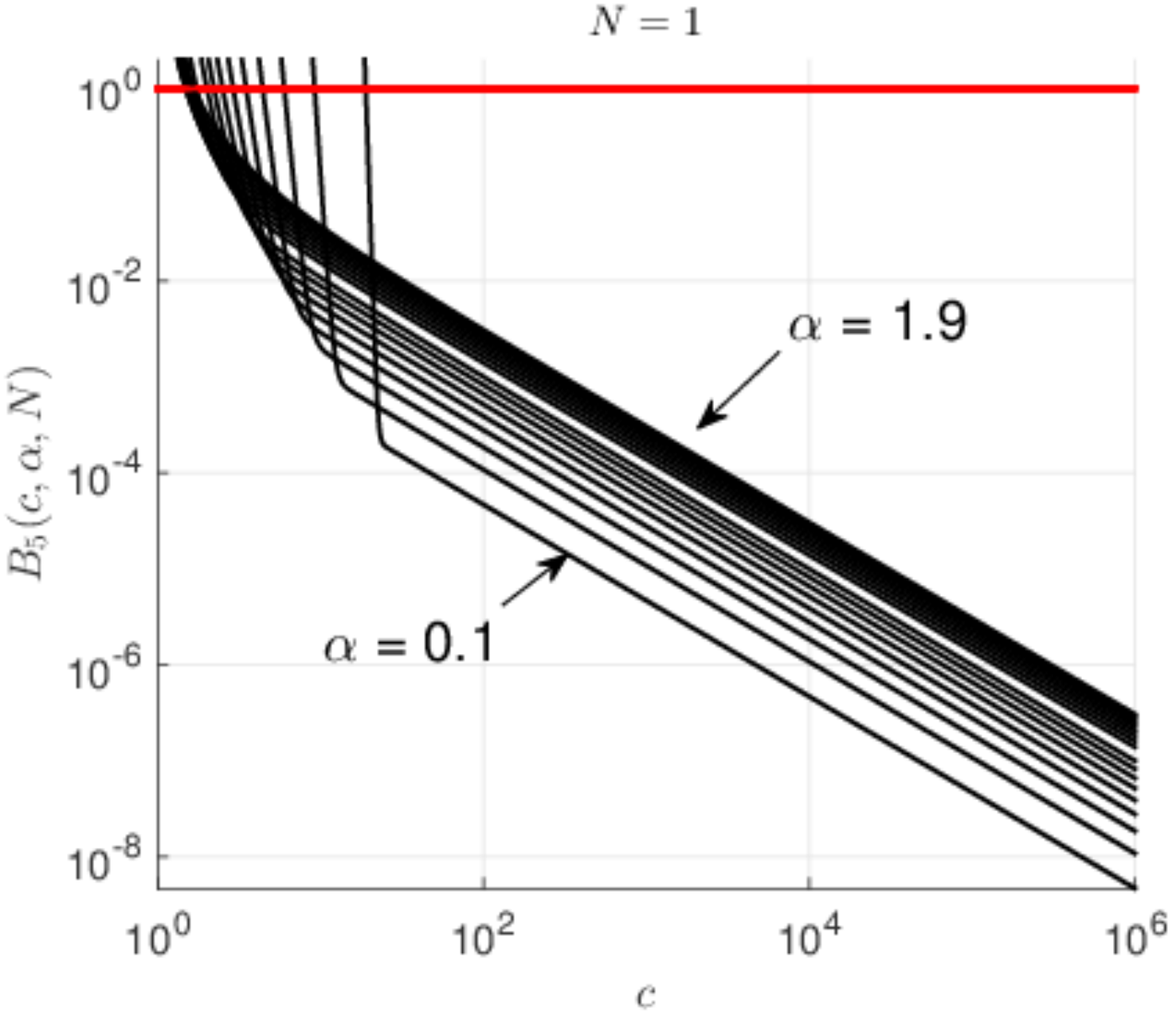}	
		\includegraphics[width=0.3\linewidth]{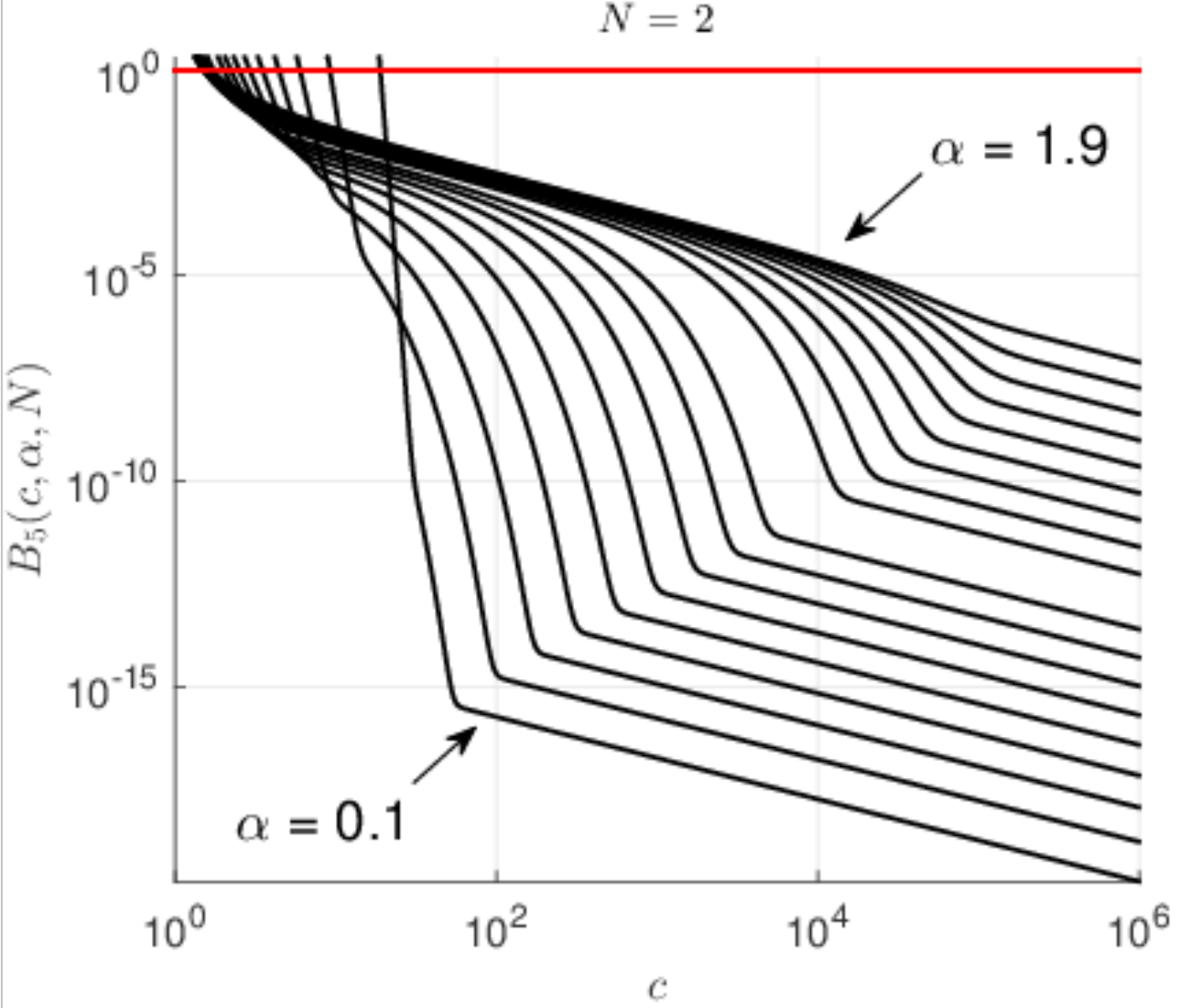}	
		\includegraphics[width=0.3\linewidth]{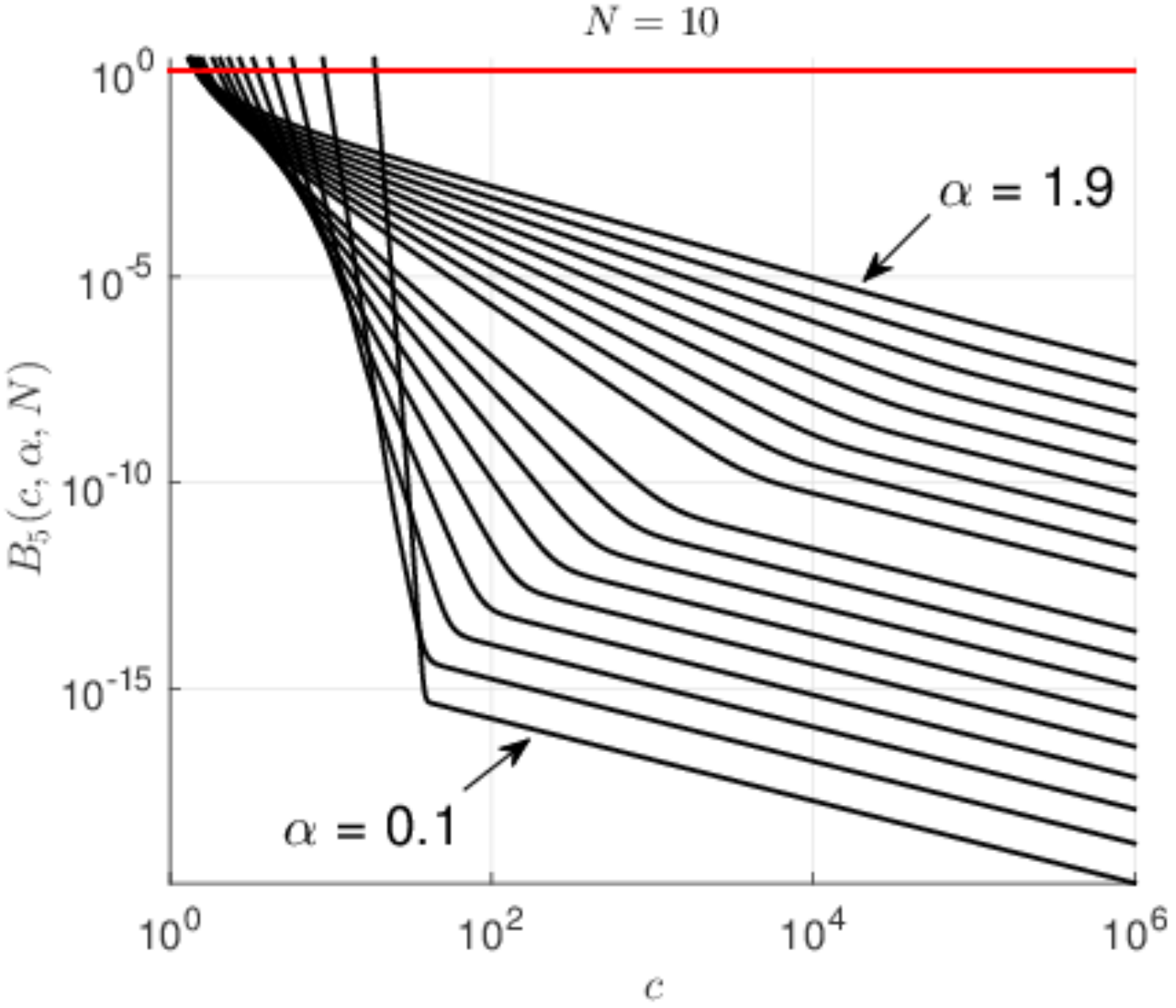}	
	}
	\caption{Bounds on $\Delta(X,\hat{X})$. 
		Each curve shows the bound $B_5(c,\alpha, N)$
		for  $\alpha=0.1,\,0.9,\ldots,\,1.9$, $\alpha \neq 1$,
		plotted against $1\leq \c\leq 10^6$.
		The $N$ points $u_0, \ldots, u_{N-1}$ are 
	logarithmically spaced on [0,1].
[The red horizontal line at 1 shows the 
maximum possible value of the Kolmogorov distance.]
	}	
	\label{fig:piecewise_lin_bounds_X_0c_change_N}
\end{figure*}

The values $\{u_i\}$ and $\{f_i\}$ serve to define a piece-wise 
linear envelope on $\omega_{X_{(0,c)}}(u)$ for $u \in[0,1]$, which is used in 
the proof; see Appendices~\ref{subsection:q_properties} 
and \ref{app:proof_hat_X_symm}.
Increasing $N$ improves (i.e., decreases) the value of
$B_5(c,\alpha,N)$, but the improvement becomes negligible
for $N\geq 10$ and logarithmically spaced points, as shown 
in Figure~\ref{fig:piecewise_lin_bounds_X_0c_change_N}, where bounds
with three different 
values of $N = 1,\,2,\,10$ are compared.

\begin{figure}[ht!] 
	\centering
	\includegraphics[width=0.4\linewidth]{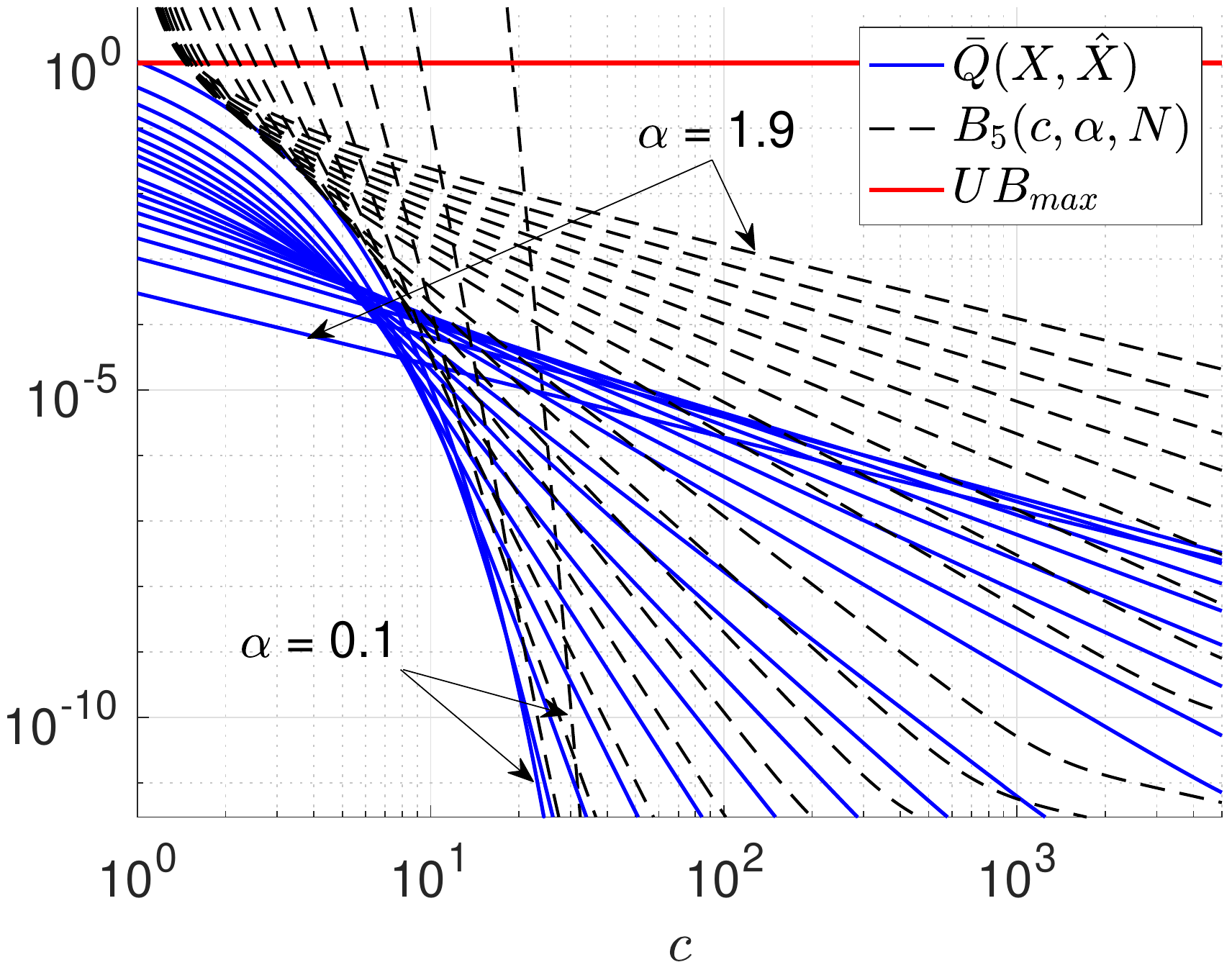} 
	\caption{
		{\small
			Bounds on $\Delta(X, \hat{X})$. The blue solid lines represent $\bar{Q}(X, \hat{X})$; 
			the black dashed lines represent $B_5(c, \alpha, N)$ with $N=10$ and logarithmically spaced points.
			The values are plotted  for $\alpha=0.1,0.2,\ldots,1.9$, $\alpha \neq 1$,
			and	for $1\leq \c\leq 5000$.
[The red horizontal line at 1 shows the 
maximum possible value of the Kolmogorov distance.]
	} }
	\label{fig:numerical_thoretical_bounds_hat_X} 
\end{figure}

In Figure \ref{fig:numerical_thoretical_bounds_hat_X} we compare 
the numerical estimates 
$\bar{Q}(X, \hat{X})$
for $\bar{I}(X, \hat{X})$ obtained in 
\cite{Riabiz2017b},
with the bound $B_5(c, \alpha, N)$ 
of Theorem~\ref{theorem:X_hat_symm}
with $N=10$.
Note that this bound 
correctly captures the dependence on $\alpha$,
and that 
the approximation error is lower for smaller values of~$\alpha$, 
a reversal of the trend shown in Figure~\ref{f:1}. 
One reason for this {is} that,
as $\alpha$ decreases, the
relative significance of the residual term becomes smaller,
when compared with the heavy-tailed initial terms in the PSR.
We also observe that the rate of convergence is dramatically better for 
smaller~$\alpha$, again in contrast with the analysis of the residual 
approximation in Figure~\ref{f:1}.
Finally, it seems that $B_5(c, \alpha, N)$ has the same asymptotic behaviour as
$B_1(c, \alpha)$ for $c \rightarrow \infty$, see also Remark 
\ref{remark:plots_appendix} in Appendix \ref{app:proof_hat_X_symm}.
However, these two bounds have reversed asymptotic ordering with respect~to~$\alpha$. 

\subsection{Nonasymptotic bound on $\Delta(X, X_{(0,c)})$}	

The following bound on $\Delta(X, X_{(0,c)})$ 
is similar
to results obtained in \cite{LedouxPaulauskas1996}. 
Its proof, given in Appendix~\ref{appendix:proof:X0c},
is based on direct computations and does not rely
on the smoothing lemma.

\medskip

\begin{proposition}\label{prop:X_0c_symm_conv}
	Let $\Delta(X, X_{(0,c)})$ be defined as in \eqref{eq:kolmog_truncated}. Under the same assumptions and notation as Theorem \ref{thm:linear},
	\begin{IEEEeqnarray*}{rCl} 
		\Delta(X, X_{(0,c)}) \leq B_6(c, \alpha) := 	\frac{\exp(-1/2)}{\sqrt{(2\pi)}}  
		&\times &
		\sqrt{\Gamma\left(        \frac{\alpha+4}{\alpha}\right)
			\left(\frac{\alpha}{4-\alpha} c^{\frac{\alpha -4}{\alpha}} + \left(\frac{\alpha}{2-\alpha}c^{\frac{\alpha-2}{\alpha}}\right)^2\right)}. 
		\IEEEeqnarraynumspace
	\end{IEEEeqnarray*}
\end{proposition}

\medskip

\begin{figure}[ht!] 
	\centering
	\includegraphics[width=0.4\linewidth]{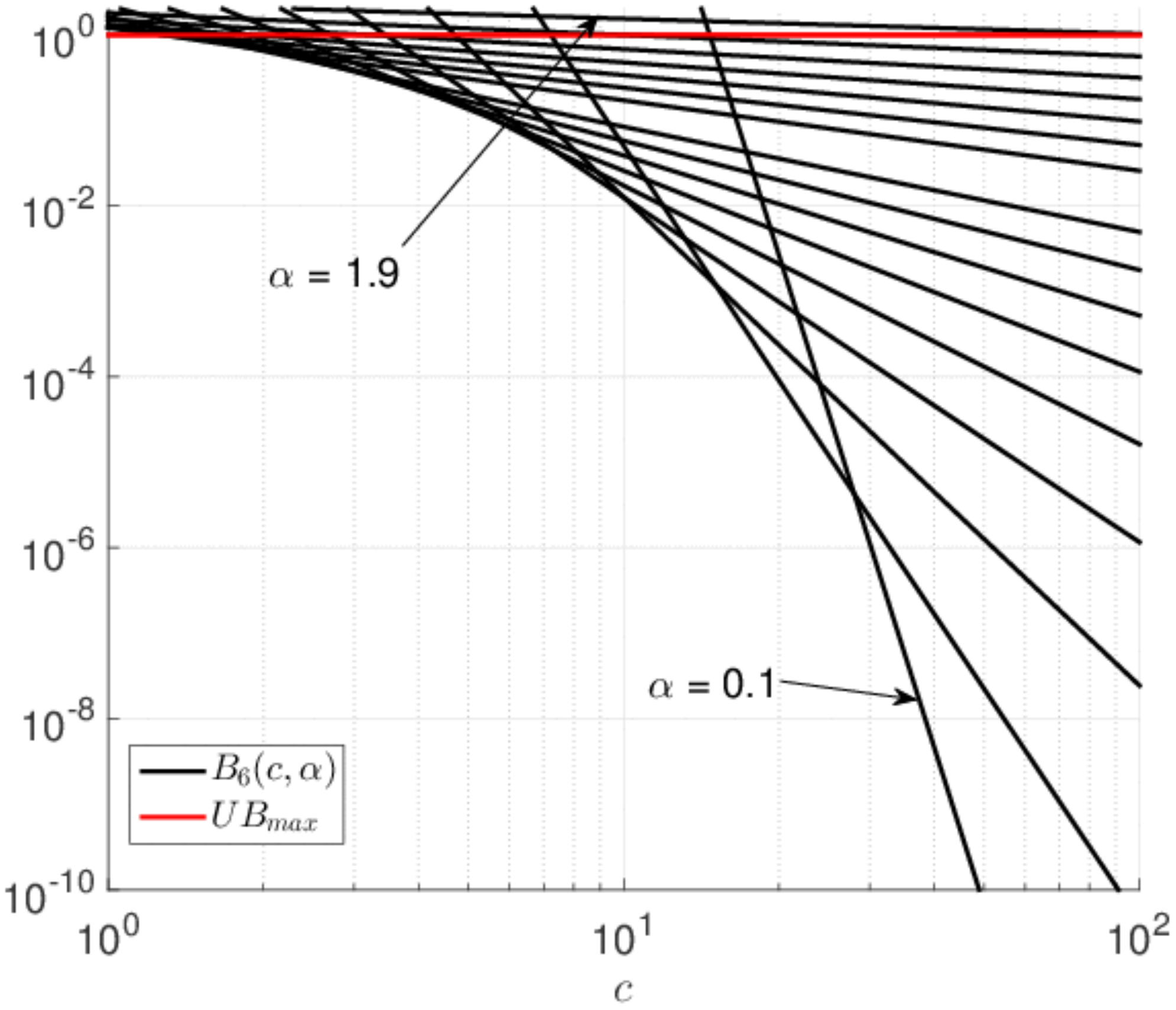} 
	\caption{Bounds on $\Delta(X,X_{(0,c)})$. 
		Each curve represents the values of the bound $B_6(c,\alpha)$
		for  $\alpha=0.1,\,0.9,\ldots,\,1.9$, $\alpha \neq 1$,
		plotted against $1\leq \c\leq 10^2$.
[The red horizontal line at 1 shows the 
maximum possible value of the Kolmogorov distance.]
}
	\label{fig:bound_B7}
\end{figure}

\begin{figure}[ht!] 
	\centering
	\includegraphics[width=0.4\linewidth]{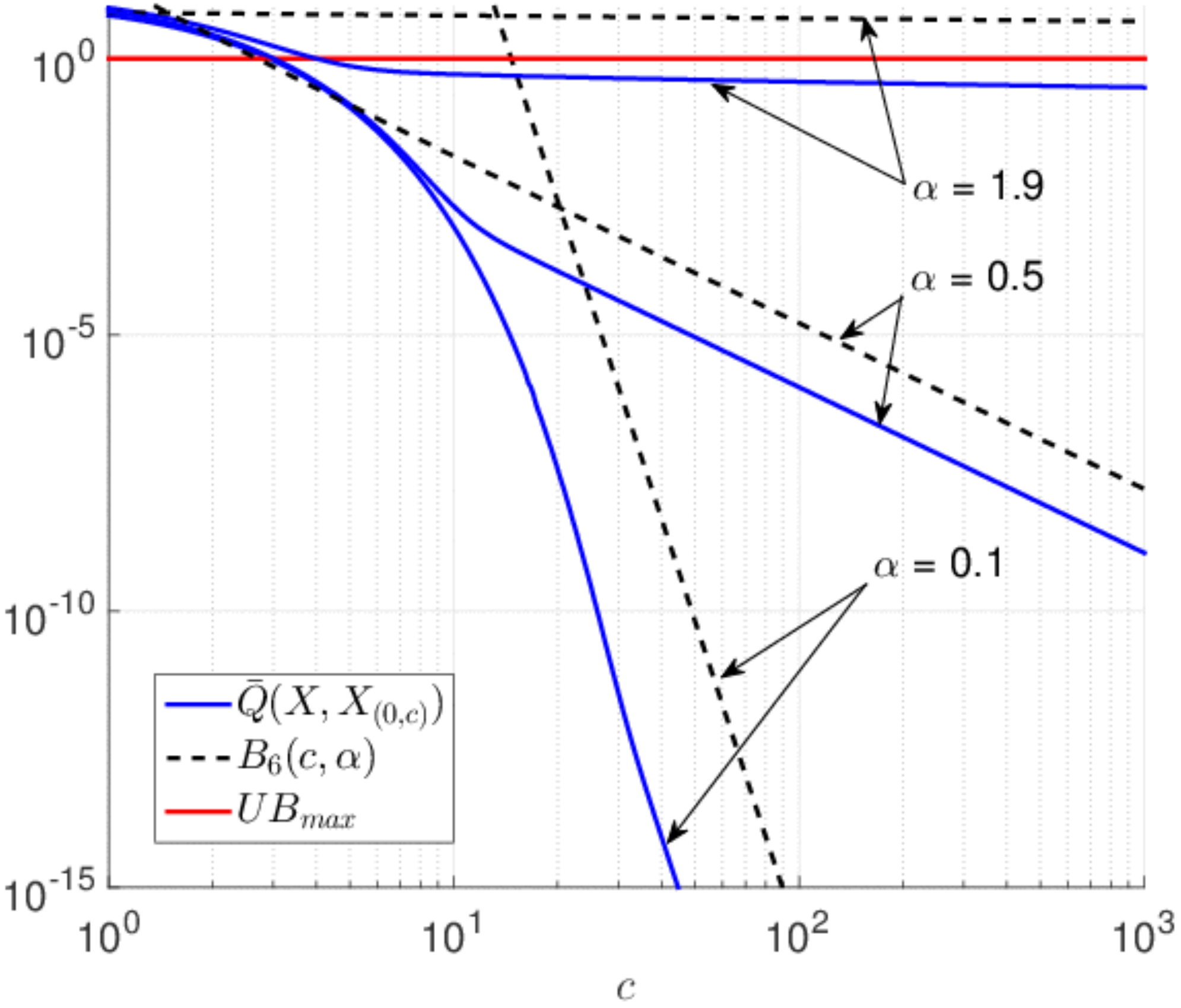} 
	\caption{
		Bounds on $\Delta(X,X_{(0,c)})$. 
		The blue solid lines represent $\bar{Q}(X, X_{(0,c)})$; 
		the black dashed lines represent $B_6(c, \alpha)$.
		The values are plotted  for $\alpha=0.1,\, 0.5,\, 1.9$,
		and	for $1\leq \c\leq 1000$.
[The red horizontal line at 1 shows the 
maximum possible value of the Kolmogorov distance.]
	}
	\label{fig:comparison_Q_B7_X_0c}
\end{figure}

It is easy to see that
the bound $B_6(c, \alpha)$ is of $O(c^{{(\alpha-2)}/{\alpha}})$,
and a corresponding lower bound of the same order is also established in
\cite{LedouxPaulauskas1996}.
Figure~\ref{fig:bound_B7} illustrates this bound,
and Figure \ref{fig:comparison_Q_B7_X_0c}
compares it to the numerical estimates 
$\bar{Q}(X,X_{(0,c)})$ of $\bar{I}(X,X_{(0,c)})$ 
as in \eqref{eq:int_berry_esseen}. 

\subsection{Comparison of $\bar{I}(X, \hat{X})$ and $\bar{I}(X, X_{(0,c)})$}	

Finally, we establish a result that compares 
the approximation of an $\alpha$-stable RV $X$
by $(i)$ the truncated PSR $X_{(0,c)}$, or
$(ii)$ by $\hat{X}=X_{(0,c)}+\hat{R}_{(0,c)}$,
which is the truncated PSR plus a Gaussian 
approximation $\hat{R}_{(0,c)}$
to the residual $R_{(0,c)}$.
Specifically, in Proposition~\ref{prop:res_contribution}, 
proved in Appendix~\ref{appendix:proof_contribution},
we compare the bounds $\bar{I}(X, X_{(0,c)})$ and 
$\bar{I}(X, \hat{X})$.
The result
indicates that,
in the symmetric case 
$W_1 \sim\mathcal{N}(0, \sigma^2_W)$,
adding a Gaussian approximation
will likely provide a better approximation 
to the $\alpha$-stable distribution than the 
truncated PSR alone, 
for most values of the truncation parameter.

\medskip

\begin{proposition} \label{prop:res_contribution}
	Let $\bar{I}(X, \hat{X})$ and $\bar{I}(X, X_{(0,c)})$ be defined as in~\eqref{eq:int_berry_esseen}, under the same assumptions and notation as in Theorem \ref{thm:linear}. 
	Then, for any $\alpha \in(0,2)$, $\alpha \neq 1$,
we have that,
$$\bar{I}(X, \hat{X})<\bar{I}(X, X_{(0,c)}),$$
for all,
	\begin{IEEEeqnarray*}{l} 
		c> c(\alpha) := \frac{\log(2)}{\gamma(1-\alpha/2,1) + e^{-1} - 1}.
		\IEEEeqnarraynumspace
	\end{IEEEeqnarray*}
\end{proposition} 

\medskip

\begin{figure}[t!] 
	\centering
	\includegraphics[width=0.4\linewidth]{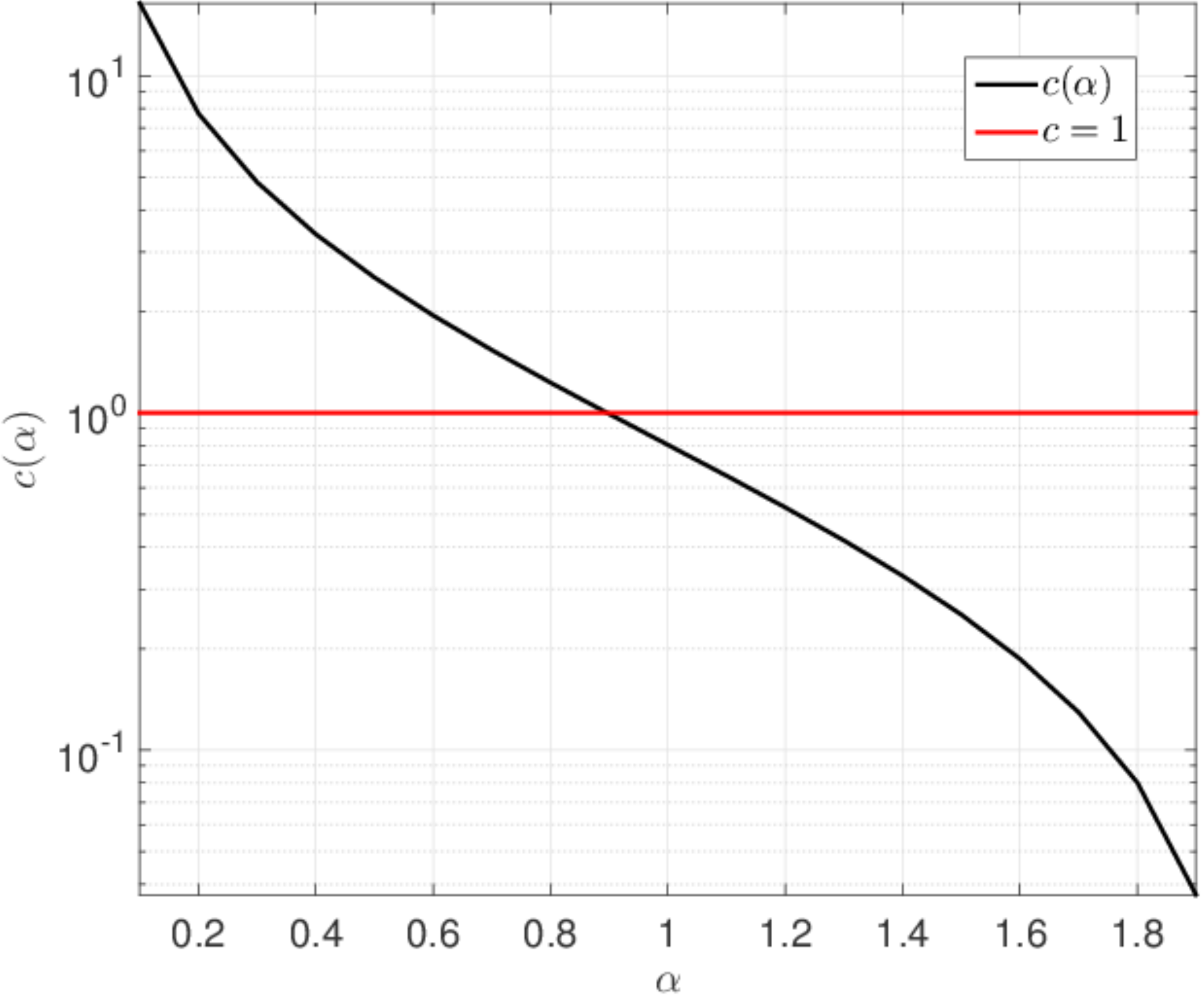} 
	\caption{The function $c = c(\alpha)$ plotted against 
$\alpha = 0.1, \ldots, 1.9.$ 
The red horizontal line corresponds to $c = 1$, which
is the smallest minimum value of $c$ that we consider
relevant for practical purposes.}
	\label{f:c_alpha}
\end{figure}

\begin{figure}[t] 
	\centerline{
		\includegraphics[width=0.4\linewidth]{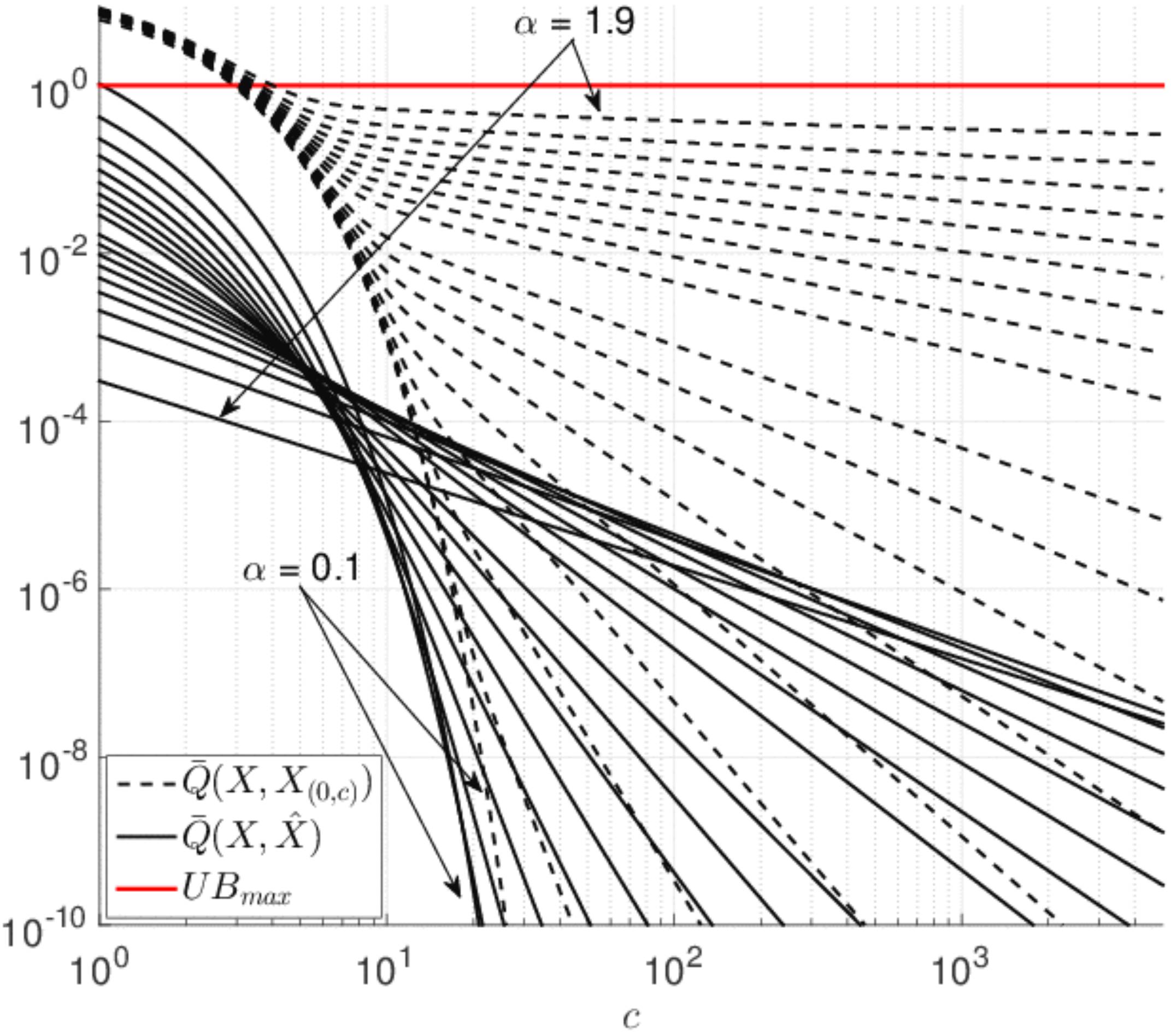}		}
	\caption{Numerical bounds on $\Delta(X, \hat{X})$ and on $\Delta(X, X_{(0,c)})$. The solid lines correspond to 
		$\bar{Q}(X, \hat{X})$ and the dashed lines to $\bar{Q}(X, X_{(0,c)})$ for $\alpha=0.1,\,0.2,\ldots,\,1.9$, $\alpha \neq 1$,
		plotted against $1\leq \c\leq 5000$. 
[The red horizontal line at 1 shows the 
maximum possible value of the Kolmogorov distance.]
}
	\label{fig:residual_improvement}
\end{figure}

Proposition~\ref{prop:res_contribution} suggests that
the Gaussian residual approximation produces a smaller approximation
error than simply truncating the series, a result also borne out 
by previous numerical results reported in
\cite[p.~56-57]{Lemke2014}; see also Figure \ref{fig:truncations_c}.
Although the result of the proposition is only valid for
$c>c(\alpha)$, we note that this is not a severe restriction:
$c(\alpha)<16.5$ for all $\alpha>0.1$ and
$c(\alpha)<1$ for $\alpha>1$; see
Figure \ref{f:c_alpha}. 
Moreover, the condition {$c>c(\alpha)$} is only shown to be
sufficient and, in fact, numerical estimates of the integrals  
$\bar{I}(X, \hat{X})$ and $\bar{I}(X, X_{(0,c)})$ show that,
\begin{IEEEeqnarray*}{rCl} 
	\bar{Q}(X,\hat{X}) < \bar{Q}(X, X_{(0,c)}), \qquad 
\mbox{for all}\; c>1, 
	\IEEEeqnarraynumspace
\end{IEEEeqnarray*}
as shown in Figure \ref{fig:residual_improvement}. 

{We remark however that the results in this section only indicate, but do not prove, that  $\hat{X}$ is closer in distribution to $X$ than $X_{(0,c)}$. This could be proved by providing lower bounds on the relative Kolmogorov distances, which is left to future studies.}

\newpage

$\;$

\newpage

\section{Inference for Regression Models}
\label{section:linear_models}

In order to illustrate the utility of the approximate conditionally Gaussian 
framework introduced in this paper, 
we give an example of a Bayesian inference scheme.
Specifically, we consider the problem of estimating
the parameters~$\bm{\lambda}$ of model~\eqref{eq:AR}, 
described in the Introduction,
when $\mathbf{u}\defeq [u_1, \ldots, u_N]$ has 
symmetric $\alpha$-stable components.

\subsection{Parameter inference in $\alpha$-stable regression models}
\label{sec:AR}
Assume, for simplicity, that 
	${\bf x}={\bf G}{\bm\lambda}+{\bf u}$
is fully observed,
and that the matrix $\bf G$ is known.
We augment this model by introducing a 
	set of latent vectors, 
	$\mathbb{T} \defeq 
	\{\bm \Gamma_1,\ldots,\bm \Gamma_N\}$, one for every element 
	of ${\bf u}=[u_1,\ldots,u_N]'$,
	as follows.
Let $c$ be a truncation parameter, and for each $1\leq n\leq N$, let
$\bm \Gamma_n \defeq [\Gamma_{1,n}, \ldots, \Gamma_{N_{(0,c)}(n),n}]$
be event times of a unit rate Poisson process,
where $N_{(0,c)}(n)$ is the number of those $\Gamma_{j,n}$ that
are smaller than $c$. We then use the approximate
representation,
$$u_n\approx\sum_{j=1}^{N_{(0,c)}(n)}W_{j,n}\Gamma_{j,n}^{-1/\alpha}+R_n,$$
where the $W_{j,n}$ are i.i.d.\ with distribution
$\mathcal{N} (0,\sigma^2_W )$,
and, according to \eqref{eq:tot_gauss}, $R_n$ is an independent Gaussian with 
mean and variance
$m_{(c,\infty)}$ and $S^2_{(c,\infty)}$
given by (\ref{eq:m_c_infty}) and (\ref{eq:S_c_infty}),
respectively, so 
that $m_{(c,\infty)}=0$ 
and 
$S^2_{(c,\infty)}=\sigma_W^2(\frac{\alpha}{2-\alpha})c^{(\alpha-2)/\alpha}$.
This, then, leads to the 
conditionally Gaussian representation of each $u_n$ as
in \eqref{eq:proxy_gauss_model},
\begin{align}\label{eq:approx_noise}
u_n | \{\Gamma_{j,n} \in [0,c]\}\overset{\text{approx}}{\sim} 
\mathcal{N} \big(0 , \sigma^2_n    	\big),
\end{align}
where
the variance $\sigma_n^2$ is,
$$\sigma_n^2
=
\sigma_W^2
\left[
\sum_{j:\Gamma_{j,n}\in[0,c]}\Gamma_{j,n}^{-2/\alpha}
+
\Big(\frac{\alpha}{2-\alpha}\Big)c^{(\alpha-2)/\alpha}
\right].
$$
We assume that the truncation parameter $c$ has been chosen
so that the Kolmogorov distance between the true distribution 
of $u_n$ and its approximation is below a certain threshold, {based on Theorem \ref{theorem:X_hat_symm}}. 

Writing $\bm\Sigma$ for the diagonal matrix with 
elements $\sigma^2_n, \, n = 1, \ldots, N$,
\eqref{eq:AR} and \eqref{eq:approx_noise} imply that the likelihood of $\bf x$ can be approximated as,
$$p(\bf x| G \bm \lambda, \, \mathbb{T}) \approx \mathcal{N}
(\bf G \bm \lambda, \,  \Sigma).$$
Regular inference can then be carried out as for Gaussian models, by 
augmenting the set of parameters to be estimated to 
$\{\bm \lambda, \mathbb{T}\}$:
A Metropolis-within-Gibbs 
sampler can be used, which in the $k$-th iteration draws,
{
\begin{align}
\label{eq:theta_step}
\bm \lambda^{(k)} &\sim p(\bm \lambda | \mathbb{T}^{(k-1)}, \bf x),
\\
\label{eq:gamma_step}
\mathbb{T}^{(k)} &\sim p(\mathbb{T} | \bm{\lambda}^{(k)}, \bf x).
\end{align}
}
The Gibbs step is \eqref{eq:theta_step}:
Adopting a
conjugate Gaussian prior leads to a Gaussian
$p(\bm \lambda | \mathbb{T}^{(k-1)}, \bf x)$.
And for sampling $\mathbb{T}$, the full conditional density
in \eqref{eq:gamma_step} can be 
targeted with a Metropolis step. 
The posterior distribution of the parameters $\bm{\lambda}$ can then 
be estimated by looking at the first component of the chain 
$\{\bm{\lambda}^{(k)}, \mathbb{T}^{(k)}\}_{k=1}^{M}$,
after it has been run for a sufficiently large number $M$ of iterations. 
We refer to \cite{Lemke2014a, Lemke2012} for more details and  
simulation results from this scheme.

\

{
 Under the assumption of Gaussian $W_j$, 
both the exact \eqref{eq:gauss_mixture} and the approximate \eqref{eq:gauss_MS_mixture_proxy} representation of each term $u_n$ are continuous scale mixtures of normals. 
However, the difference between the two representations  lies in the possibility of exactly sampling the mean and scale latent variables.
}

{
	In the presented inference scheme, the distribution of the terms $u_n$ is symmetric.\footnote{{The mean mixing component is here set to 0, in order to allow for a choice of the truncation parameter $c$ based on Theorem \ref{theorem:X_hat_symm}.}}  When $\alpha$ is known, an exact continuous scale mixture of normals representation can be deduced for the symmetric stable distribution from the product property \cite[p.~176]{Feller1966}. In the aforementioned case with $\alpha$ known, the scale latent variable has positive stable distribution, and can thus be sampled exactly via the method of \cite{Chambers1976}. 
	This was used in \cite{Godsill1999, Tsionas1999, GodsillKuruoglu1999, Godsill2000} for developing posterior samplers for the  parameters of the stable distribution and the parameters of linear models with stable noise. 
}

{		
	However, 
	for asymmetric stable distributions, or for symmetric distributions with unknown $\alpha$, it appears not to be possible to sample the mean and scale latent variables in closed form. 
	Consequently, in the cases above it is not possible to do inference based on the exact PSR.  On the other hand, it is possible to do  exact inference for the approximate representation. 
	Exploiting this possibility is perhaps the most relevant continuation of this work.
}
	
{
		Finally, we stress that our findings are indeed related to the stable likelihood, and further analytic studies are required in order to establish the quality of  inference procedures based on  truncating the PSR and approximating its residual. 
		 We refer to \cite{Lemke2014} for numerical insights on superior behaviour of Bayesian estimators based on the truncated PSR with accounted residual, with respect to those based on the simply truncated PSR.  
}

\medskip

As anticipated in the Introduction, the PSR and its approximation 
\eqref{proxy_PSR_split} are also relevant for inference of continuous-time 
\levy~processes and linear models driven by \levy~noise. These are more 
challenging in terms of inference than the discrete time model 
\eqref{eq:AR}--\eqref{eq:AR:noisy}, but in some applications they are 
more realistic, e.g., when data are sampled at irregular time intervals. 
Potential extension of our work in this direction are described next.

\subsection{PSR for continuous-time stochastic processes and linear models}
\label{section:CLT:st:integral}

Consider the continuous-time version of the linear model corresponding 
to \eqref{eq:AR}, where we assume that linear observations are made
at discrete times~$\{t_i\}$, as in~\eqref{eq:AR:noisy},
\begin{IEEEeqnarray}{rCl}
	\dif \mathbf{x}(t) &=& \mathbf{A}\mathbf{x}(t) \dif t + \mathbf{h} \dif \ell(t), \label{eq:CAR1}
	\\
	y(t_i) &=& \mathbf{b}'\mathbf{x}(t_i)  + v(t_i),
	\nonumber
\end{IEEEeqnarray}
where $\mathbf{x}(t) = \left[ x_{1}(t), \ldots ,x_{P}(t)\right]'$ is the state, 
$\mathbf{A}$ is a $P\times P$ matrix describing the interaction of the 
components of $\mathbf{x}(t)$, $\mathbf{h}$ is a $P$-dimensional vector 
describing the effects of the noise process $\{\dif \ell(t)\}$, 
$\mathbf{b}$ is a $P$-dimensional vector, and 
$\{v(t_i)\}$ is the observation noise process.
A wide range of results have been developed 
in the literature for the case
when 
$\{\ell(t)\}$ is a Brownian motion 
\cite{Oksendal2003, Harvey1990}, but, as for the discrete-time case, 
such models are not appropriate for certain applications.

Large jumps and heavy tails in the state process,
as often observed in applications,
can be modelled 
by choosing $\{\ell(t)\}$ to be a (non-Gaussian) \levy~process; see \cite{Brockwell2001, Brockwell2004, Brockwell2009} for a first formulation, \cite{TankovCont2015} 
and \cite{Barndorff-NielsenMikoschResnick2001} 
for a review, and \cite{UnserTaftiSun2014, UnserTaftiAminiEtAl2014, AminiUnser2014} for  more recent work from an engineering perspective.
The sub-class of $\alpha$-stable \levy\ processes \cite{Samoradnitsky1994,JanickiWeron1994} 
is of special importance. 
In fact,
the self-similarity of stable \levy~processes \cite{Samoradnitsky1994} implies that transition densities, although still intractable, all come from the same $\alpha$-stable family. Hence, as argued in the Introduction,
$\alpha$-stable \levy~processes 
may be considered to be the natural first choice towards generalising 
the classical Gaussian process framework to the heavy-tailed case.  

LePage series representations are available for \levy~processes. In the $\alpha$-stable symmetric case, $\{\ell(t)\}$ can be expressed,
\begin{IEEEeqnarray}{rCl}
	\ell(t) 	\overset{\mathcal{D}}{=} 
	\sum_{j=1}^{\infty} (\Gamma_j^{-1/\alpha} W_j) \mathds{1}{(V_j< t)}
	, 
	\label{eq:PSR_levy_alpha}
\end{IEEEeqnarray}
where $\{\Gamma_j\}_{j=1}^\infty$ and $\{W_j\}_{j=1}^\infty$ are as before,
and $\{V_j\}_{j=1}^\infty$ are i.i.d.\ uniform RVs on $[0,T]$, where $T
$ is the time horizon considered.
A similar representation is available for the asymmetric case, 
see \cite[Lemma~4.1.1]{Lemke2014a}.
Like in discrete-time,
$\ell(t)$ is conditionally Gaussian,
\begin{IEEEeqnarray*}{rCl}
	\ell(t) | \{\Gamma_j,V_j\}_{j=1}^\infty\sim  
	\mathcal{N}(m_t, S^2_t),
\end{IEEEeqnarray*}
with $m_t$ and $S^2_t$ can be appropriately defined 
as series involving $\{\Gamma_j,V_j\}_{j=1}^\infty$.
Once again, the series in (\ref{eq:PSR_levy_alpha}) and
in the definitions of $m_t$ and $S^2_t$ cannot be computed exactly. 
However, defining $\tilde{W}_j:= W_j \mathds{1}{(V_j< t)}$, we have that {
$\mathbb{E}[\tilde{W}_1^2]< \infty$}, hence Theorem \ref{th:CLT} 
applies and the residual of \eqref{eq:PSR_levy_alpha} is asymptotically 
Gaussian. {If the variables $W_j$ are selected to be Gaussian}, this implies that an overall approximate conditionally 
Gaussian representation of $\ell(t)$ is again available. 
Furthermore, considering $\mathds{1}{(V_j< t)}$ as a thinning operation on the 
Poisson process associated with the RVs $\{\Gamma_j\}$, it is 
straightforward to extend our present results
to the continuous-time setting of $\alpha$-stable \levy~processes. 

The stable distribution of $\ell(t)$ is inherited by 
$\mathbf{x}(t)$, 
so a PSR representation holds for the $\alpha$-stable vector
$\mathbf{x}(t)$; a Gaussian approximation result for its residual 
is discussed in {\cite{Riabiz2017a}}.  Note, however, that our present
{nonasymptotic} results need to be further adapted to $\mathbf{x}(t)$,
due to the fact that we would need to consider the structure 
of the solution to the stochastic differential equation  \eqref{eq:CAR1}, 
and the multivariate nature of $\mathbf{x}(t)$. 
Preliminary versions of these
ideas  have been implemented in \cite{Lemke2011, Lemke2014a, Lemke2015a}, {by accounting for the PSR residual in inference tasks, while theoretical studies on the choice of the threshold parameter $c$ are left to future developments}.

\newpage

\appendices

\section{Transformations for $\sigma$ and $\beta$}
\label{app:sigma_beta}
Suppose that
$W_1 \sim \mathcal{N}(\mu_W, \sigma^2_W)$
or, more generally, that $W_1$ belongs to a location-scale family,
with location $\mu_W$, scale $\sigma_W^2$, and
PDF $f_W(w)$, $w \in \mathbb{R}$.
For a given value of the tail parameter $\alpha$, 
here we describe how 
any pair of values $\sigma >0$ and $\beta \in [-1,1]$
can be obtained via 
the mappings in \eqref{eq:mapping_sigma} and \eqref{eq:mapping_beta},
by appropriate choices of $\mu_W$ and $\sigma_W^2$.

We  introduce the following auxiliary PDF,
\begin{IEEEeqnarray*}{c} 
	\pi(w) := \dfrac{|w|^{\alpha} f_W(w)}{\int_{\mathbb{R}} 
	|w'|^{\alpha} f_W(w') \dif w'} = \dfrac{\tilde{\pi}(w)}{I}, 
\end{IEEEeqnarray*}
where $\tilde{\pi}(w)$ and $I$ denote the unnormalized density and the normalizing constant, respectively. 
Then \eqref{eq:mapping_sigma}, the transformation related to $\sigma$, can be rewritten as,
\begin{IEEEeqnarray*}{c}
	\sigma = \dfrac{\int_{\mathbb{R}}\tilde{\pi}(w) \dif w}{C_\alpha} = \dfrac{I}{C_\alpha}. 
\end{IEEEeqnarray*}
Since $C_\alpha>0$, we have that $\sigma>0$, and it is easy to see
that it is possible to achieve any $\sigma>0$, 
by appropriately choosing~$\sigma_W^2$.

Similarly, we can express \eqref{eq:mapping_beta}, the 
transformation related to $\beta$, as,
\begin{IEEEeqnarray*}{rCl}
	\beta & = & \dfrac{- \int_{-\infty}^0\tilde{\pi}(w) \dif w + \int_{0}^{\infty}\tilde{\pi}(w) \dif w}{I} \nonumber
	\\
	& = & - \int_{-\infty}^0{\pi}(w) \dif w + \int_{0}^{\infty}{\pi}(w) \dif w \nonumber
	\\
	& = & - \left( 1- I^+ \right) + I^+ \nonumber
	\\
	& = & 2I^+ - 1, 
	\IEEEeqnarraynumspace
\end{IEEEeqnarray*}
where $I^+$ is the probability of $\mathbb{R}^+$ under $\pi$.
Then any $\beta \in [-1,1]$ can be obtained by choosing the parameters of the 
distribution of $W_1$ to give the required value of $I^+ \in [0,1]$.
Since that $|w|^\alpha$ is a symmetric function, it is clear 
from \eqref{eq:mapping_beta} that $\beta=0$ when 
$f_W(w)$ is an even function.
In the Gaussian case, this corresponds to $\mu_W =0$. Similarly, $I^+>0.5$ 
(i.e., $\mu_W>0$) leads to positive skewness $\beta>0$, while $I^+<0.5$ 
(i.e., $\mu_W<0$) leads to $\beta<0$.   
A combined choice of the scale and location parameters is required
to achieve the limiting cases $\beta = -1$ ($\mu_W <0, \sigma_W =0$)
and $\beta = 1$ ($\mu_W >0, \sigma_W =0$).

\section{Proof of Lemma \ref{lem:finite_res_momemts}}
\label{app:finite_res_momemts}

We make use 
of the observation that $R_{(c,d)}$ can be viewed as a compound Poisson 
process. Hence, to compute expectations with respect to its distribution, 
we first condition on $N_{(c,d)}$, the random number of terms in $R_{(c,d)}$, 
and take the expectation over $N_{(c,d)}$.
Using the expression \eqref{eq:Rcd_new} for $R_{(c,d)}$,
\begin{IEEEeqnarray*}{rCl} 
	\phi_{R_{(c,d)}}(s) 
	&=&
	\mathbb{E}\big[ \exp \left( i s R_{(c,d)}\right) \big] 
	\\
	&=&
	\mathbb{E}\big[\mathbb{E}\big[ \exp\big(i s \big(\sum_{j=1}^{N_{(c,d)}} 
	W_j \Gamma_j^{-1/\alpha}
	-\mathbb{E}[W_1]\sum_{j=1}^{d}b_j^{(\alpha)} \big) \big)| {N_{(c,d)}}\big] \big]
	\\
	&=&
	\mathbb{E}\big[\mathbb{E}\big[ \exp\big(i s \big(\sum_{j=1}^{N_{(c,d)}} 
	W_j U_j^{-1/\alpha}
	-\mathbb{E}[W_1]\sum_{j=1}^{d}b_j^{(\alpha)} \big) \big)| {N_{(c,d)}}\big] \big]
	\\
	&=&
	\mathbb{E}\left[\mathbb{E}\left[ e^{i s \sum_{j=1}^{N_{(c,d)}} Y_j
	}\middle| {N_{(c,d)}}\right] \right]
	\times
	\exp\big(-i s \mathbb{E}[W_1]\sum_{j=1}^{d}b_j^{(\alpha)} \big)
	\\
	&=&
	\sum_{{n=0}}^{\infty}\mathbb{E}\left[ e^{i s \sum
		_{{j=1}}^n Y_j}\right] \mathbb{P}\left({N_{(c,d)}} = n\right) 
	\times
	\exp\left( -i s B \right) 
	\\& = &
	\sum_{n=0}^{\infty}(\phi_{Y_1}(s))^n \frac{(d-c)^n}{n!}e^{-(d-c)} 
	\times
	\exp\left( -i s B \right)
	\\& = &
	\exp\left((d-c)(\phi_{Y_1}(s)-1)	
	-i s B \right), 
\end{IEEEeqnarray*}

\noindent
where $\phi_{Y_1}(s)$ is the CF of $Y_1 = W_1 U_1^{-1/\alpha}$. 

Since we assume $\mathbb{E}[W_1^2]$ is finite, it follows that $\mathbb{E}[Y_1^2]$ is finite,
and hence the first and second moments of 
$R_{(c,d)}$ are finite as well and can be computed
\cite[Lemma~XV.4.2]{Feller1966}
by taking derivatives of its CF at zero:
\begin{IEEEeqnarray}{rCl} 
	m_{(c,d)} 
	&=& 
	\mathbb{E}[R_{(c,d)}]
	\nonumber\\
	&=& 
	(-i)\phi'_{R_{(c,d)}}(0)
	\nonumber\\
	&=&
	(d-c)\mathbb{E}[Y_1] -B,
	\label{eq:non_asympt_mean}
	\\
	S^2_{(c,d)}&= & 
	\mathbb{E}[R_{(c,d)}^2]
	-(\mathbb{E}[R_{(c,d)}])^2
	\nonumber\\
	&=&(-i)^2\phi''_{R_{(c,d)}}(0)
	-\big[(-i)\phi'_{R_{(c,d)}}(0)]^2
	\nonumber\\
	&=&
	(d-c) \mathbb{E}[Y_1^2].
	\label{eq:non_asympt_var}
\end{IEEEeqnarray}
\noindent
Since 
$U_1  \overset{\text{i.i.d.}}{\sim} \mathcal{U}(c,d)$, we have that, 
for $k=1,2$, 
\begin{IEEEeqnarray*}{rCl} 
	\mathbb{E}{[}U_1^{-k/\alpha}{]}
	&=&
	\frac{1}{d-c} \int_c^d U^{-k/\alpha} \dif U \nonumber 
	\\
	&=& 
	\frac{1}{d-c}\left[\frac{\alpha}{\alpha-k}U^{-k/\alpha+1}\right]^d_c \nonumber 
	\\
	&= &
	\frac{1}{d-c}\frac{\alpha}{\alpha-{k}}\left(d^{\frac{\alpha-k}{\alpha}}-c^{\frac{\alpha-k}{\alpha}}\right), 
	\label{eq:moments_uniform}
\end{IEEEeqnarray*}
\noindent    	    		
and since $W_1$ is independent of $U_1$,
\begin{IEEEeqnarray}{c} 
	\mathbb{E}[Y_1^k]  = \mathbb{E}[W_1^k] \frac{1}{(d-c)}\frac{\alpha}{\alpha-k}\left(d^{\frac{{\alpha-k}}{\alpha}}-  c^{\frac{{\alpha-k}}{\alpha}}\right). 		
	\IEEEeqnarraynumspace
	\label{eq:Yk_moments}
\end{IEEEeqnarray}
Substituting \eqref{eq:Yk_moments} into the moment expressions \eqref{eq:non_asympt_mean} and \eqref{eq:non_asympt_var} gives~(\ref{eq:m_cd}) 
and~(\ref{eq:S_cd}) as claimed.

\section{Proof of Lemma \ref{lemma:cf_series_R_cd}}\label{app:series_cf}

As in the proof of Theorem~\ref{th:CLT} we note that the CF of $Z_{(c,\infty)}$
can be expressed in terms of the CF of $Y_1=W_1U_1^{-1/\alpha}$.
Therefore, we begin by expanding $\phi_{Y_1}(s)$ as a Taylor series.

Suppose $W_1\sim{\cal N}(\mu_W,\sigma_W^2)$. Using~(\ref{eq:Yk_moments})
and the well-known formula for the moments of the normal distribution,
it is easy to check that,
\begin{IEEEeqnarray}{rCl} 
	\limsup_{k \rightarrow \infty} \frac{1}{k}\mathbb{E}[|Y_1|^k]^{1/k} 
< \infty.
	\label{eq:suff_cond_analiticity}
\end{IEEEeqnarray}	
Therefore, $\phi_{Y_1}(s)$ is analytic around $s=0$ and admits
the Taylor expansion
\cite[p.~514]{Feller1966},
\begin{IEEEeqnarray}{rCl} 
	\phi_{Y_1}\left({s}\right) 
	& = &1 + i\mathbb{E}[Y_1]s -\frac{1}{2}\mathbb{E}[Y_1^2]s^2 + \sum_{k=3}^{\infty} \frac{i^k}{k!} \mathbb{E}[Y_1^k]s^k.
	\IEEEeqnarraynumspace
	\label{eq:phi_Y}
\end{IEEEeqnarray}
\noindent
Then, from \eqref{eq:cf_R_cd},
\eqref{eq:moments_Y} 
and \eqref{eq:Yk_moments}, for $c>0$, we have,
\begin{IEEEeqnarray}{rCl} 
	\phi_{R_{(c,d)}}(s)  
	& = &  
	\exp\left( ism_{(c, d)}
	- \frac{s^2S^2_{(c, d)}}{2} +
	\sum_{k=3}^\infty r_k s^k\right), 
	\IEEEeqnarraynumspace
	\label{eq:cf_rcd_series}
\end{IEEEeqnarray}
where,
\begin{IEEEeqnarray}{rCl} 
	r_k 
	& : = &  \frac{i^k}{k!}\mathbb{E}[W_1^k]\frac{\alpha}{\alpha-k}\left(d^{(\alpha-k)/\alpha}-c^{(\alpha-k)/\alpha}\right).
	\IEEEeqnarraynumspace
	\label{eq:r_k}
\end{IEEEeqnarray}	
Taking the limit as $d \rightarrow \infty$, 
\begin{IEEEeqnarray}{rCl} 
	\phi_{R_{(c,\infty)}}(s)  
	& = &  
	\exp\left( ism_{(c, \infty)}
	- \frac{s^2S^2_{(c, \infty)}}{2} +\sum_{k=3}^\infty \bar{r}_k s^k\right), 
	\label{eq:cf_R_cd_series}
	\IEEEeqnarraynumspace
\end{IEEEeqnarray}
\noindent 
where,
\begin{IEEEeqnarray}{rCl} 
	\bar{r}_k 
	= \frac{i^k}{k!}\mathbb{E}[W_1^k]\frac{\alpha}{k - \alpha }c^{(\alpha-k)/\alpha}, 
\quad k \geq 3.
	\IEEEeqnarraynumspace
	\label{eq:bar_r1}
\end{IEEEeqnarray}
The justification for taking the term-by-term limit of the series 
in~(\ref{eq:phi_Y}) follows a standard argument. Since it is a power
series, it converges absolutely. Now, for any fixed $s$, the coefficients
$r_k$ are continuous functions of $d$, and they are uniformly bounded
in absolute value for $d>\max\{1,c\}$. Therefore, the series converges 
uniformly in $d$, which implies that we can take its term-by-term limit.

From \eqref{eq:cf_res_standardized} and \eqref{eq:S_c_infty}, 
\begin{IEEEeqnarray*}{rCl} 
	\phi_{Z_{(c,\infty)}}(s)  & = &  
	\exp\Big(
	- \frac{s^2}{2} +\lim_{d\rightarrow \infty}
	\sum_{k=3}^\infty z_k s^k\Big)
	\\
	& = &  
	\exp\Big(
	- \frac{s^2}{2} +\sum_{k=3}^\infty \bar{z}_k s^k\Big), \nonumber
	\label{eq:cf_Z_cd_series}
	\IEEEeqnarraynumspace
\end{IEEEeqnarray*}
\noindent
where the term-by-term limit can be justified as before, 
\begin{IEEEeqnarray*}{
		rCl} 
	z_k &:=& \frac{i^k}{k!}\frac{(d-c)\mathbb{E}[Y_1^k]}{\left( (d-c) \mathbb{E}[Y_1^2]\right)^{k/2} }\\
	& = & 
	\frac{i^k}{k!}\frac{\mathbb{E}[W_1^k]\frac{\alpha}{\alpha-k}\left(d^{(\alpha-k)/\alpha}-c^{(\alpha-k)/\alpha}\right)}
	{\left(\mathbb{E}[W_1^2]\frac{\alpha}{\alpha-2}(d^{(\alpha-2)/\alpha}-c^{(\alpha-2)/\alpha})\right)^{k/2}},
	\label{eq:z_k}
	\IEEEeqnarraynumspace
\end{IEEEeqnarray*}
and,
\begin{IEEEeqnarray*}{
		rCl} 
	\bar{z}_k &:=& 
	\frac{i^k}{k!}\frac{\mathbb{E}[W_1^k]\frac{\alpha}{k - \alpha}}
	{\left(\mathbb{E}[W_1^2]\frac{\alpha}{2- \alpha}\right)^{k/2}} c^{1 - k /2},
\quad k \geq 3.
	\IEEEeqnarraynumspace
\end{IEEEeqnarray*}

\section{Proof of Lemma \ref{lemma:cf_integral_R_cd}}\label{app:integral_cf}

\subsection{CF of $R_{(c, \infty)}$}
For convenience, we write 
$Y_1 = W_1 t(U_1),$
where $t(x):= x^{-1/\alpha}$. Then we have,
\begin{IEEEeqnarray*}{rCl} 
	\phi_{Y_1}(s)&=&E[e^{is Y_1}]\\
	&=&E[e^{is W_1t(U_1)}]\\
	&=&\int_{c}^{d}\int_{\mathbb{R}}e^{i s wt(u)}p(w,u)\dif w \dif u\\
	&=&\int_{c}^{d}\int_{\mathbb{R}}e^{is wt(u)}p(w)p(u)\dif w \dif u\\
	&=&\int_{c}^{d}\frac{\phi_{w}(s t(u))}{d-c}\dif u,
\end{IEEEeqnarray*}
\noindent 
where $p(w,u)$ is the joint density of the random vector $(W_1,U_1)$,
and $p(w)$ and $p(u)$ are the respective marginals.
Now substituting the expression for the CF of a Gaussian  
$W_1\sim\mathcal{N}(\mu_W, \sigma^2_W)$, 
\begin{IEEEeqnarray*}{rCl} 
	\phi_{Y_1}(s)&=&\int_{c}^{d}\frac{\exp(is t(u)\mu_W-\sigma_W^2s^2t(u)^2/2)}{d-c}\dif u\\
	&=&\frac{1}{d-c}\int_{t(c)}^{t(d)}\frac{\exp(is t\mu_W-\sigma_W^2s^2t^2/2)}{\dif t/\dif u}\dif t\\
	&=&\frac{\alpha}{d-c}\int_{d^{-1/\alpha}}^{c^{-1/\alpha}}\exp(is t\mu_W-\sigma_W^2s^2t^2/2)t^{-\alpha-1}\dif t,
\end{IEEEeqnarray*}
and recalling \eqref{eq:cf_R_cd}, 
\begin{IEEEeqnarray}{rCl} 
	\log
	(\phi_{R_{(c,d)}} (s))  
	& = &
	c-d
	+    
	{\alpha}\int_{d^{-1/\alpha}}^{c^{-1/\alpha}}e^{(is t\mu_W-\sigma_W^2s^2t^2/2)}t^{-\alpha-1}\dif t - isB.
	\IEEEeqnarraynumspace
	\label{eq:cf_Rcd_integral_1}
\end{IEEEeqnarray}
By simple algebra we can re-write,
\begin{IEEEeqnarray}{rCl} 
	\log
	(\phi_{R_{(c,d)}} (s))  
	& = &
	{\alpha}\int_{d^{-1/\alpha}}^{c^{-1/\alpha}}\left( e^{(is t\mu_W-\sigma_W^2s^2t^2/2)} - 1 -is\mu_Wt\right) t^{-\alpha-1}\dif t	
	- is \mu_W \frac{\alpha}{\alpha -1}c^{\frac{\alpha-1}{\alpha}},
	\IEEEeqnarraynumspace
	\label{eq:cf_Rcd_integral_2}
\end{IEEEeqnarray}
\noindent
and taking the limit $d\rightarrow \infty $, 
\begin{IEEEeqnarray*}{rCl} 
	\log
	(\phi_{R_{(c,\infty)}} (s))  
	& = &
	{\alpha}\int_{0}^{c^{-1/\alpha}}\left( e^{(is t\mu_W-\sigma_W^2s^2t^2/2)} - 1 -is\mu_Wt\right) t^{-\alpha-1}\dif t	
	- is \mu_W \frac{\alpha}{\alpha -1}c^{\frac{\alpha-1}{\alpha}}. 
\end{IEEEeqnarray*}

\subsection{CFs of $X_{(0,c)}$ and $X$}
Notice that no assumptions on the finiteness of the moments 
of $R_{(c,d)}$ or $Y_1$ were made above. Therefore,
we can use the expression of the CF of $R_{(c, d)}$ with 
$c\downarrow 0$ to obtain an analogous expression for 
the CF of $X_{(0,c)}$. In this case, the term
$i B s$ appearing in \eqref{eq:cf_Rcd_integral_1} 
is included in the  PSR residual $R_{(c,d)}$, 
see \eqref{eq:reminder_PSR_finite}. Hence we have,
\begin{IEEEeqnarray}{rCl} 
	\log
	(\phi_{X_{(0,c)}}(s))  
	& = &
	-c
	+    
	\int_{c^{-1/\alpha}}^{\infty}e^{(is t\mu_W-\sigma_W^2s^2t^2/2)}t^{-\alpha-1}\dif t 	\label{eq:phi_X_0c_gauss}
	\\
	& = &
	{\alpha}\int_{c^{-1/\alpha}}^{\infty}\left( e^{(is t\mu_W-\sigma_W^2s^2t^2/2)} - 1 \right) t^{-\alpha-1}\dif t.	
	\IEEEeqnarraynumspace
	\nonumber
\end{IEEEeqnarray}

\noindent 
Finally we obtain an integral expression for the CF of the full PSR $X$,
\begin{IEEEeqnarray}{rCl} 
	\log
	(\phi_{X} (s)) & = & 
	\log
	(\phi_{X_{(0,c)}}(s))  + 
	\log
	(\phi_{R_{(c,\infty)}} (s))  
							\nonumber
	\\
	& = &
	{\alpha}\int_{0}^{\infty}\left( e^{(is t\mu_W-\sigma_W^2s^2t^2/2)} - 1 \right) t^{-\alpha-1}\dif t. 	\label{eq:phi_X_0inf_gauss}
\end{IEEEeqnarray}

\section{Proof of Lemma \ref{lemma:cf_res_symm} }\label{appendix:proofs_cf_symm}

\subsection{CF of $Z_{(c, \infty)}$}

Starting from the expression for $\phi_{R_{(c,d)}}(s)$ 
in \eqref{eq:cf_Rcd_integral_1}, we take $\mu_W=0$ and perform 
the change of variables,
\begin{eqnarray}
	v=\sigma_W^2s^2t^2/2, 
	\label{eq:ch_var_t_s}
\end{eqnarray}
so that $t ={\sqrt{2v}}/{\sigma_W|s|}$,
$dv/dt=\sigma_W^2s^2t=\sqrt{2v}\sigma_W |s|$,
and the indefinite integral corresponding to the 
definite integral in~\eqref{eq:cf_Rcd_integral_1}
becomes,
\begin{IEEEeqnarray}{rCl} 
	\int e^{(-\sigma_W^2s^2t^2/2)}t^{-\alpha-1}\dif t
	&=&
	|s|^\alpha\sigma_W^{\alpha}2^{-(\alpha+2)/2}\int\exp(-v)v^{-\alpha/2-1}\dif v \nonumber
	\\&=&
	|s|^\alpha\sigma_W^{\alpha}2^{-(\alpha+2)/2}\left[-\frac{1}{a} v^{-a} e^{-v} - \int \frac{1}{a}v^{-a} e^{-v} \dif v \right] \nonumber
	\\
	&=&|s|^\alpha\sigma_W^{\alpha}2^{-(\alpha+2)/\alpha}\frac{1}{a}
	\Big(-\exp(-v)v^{-a}-\gamma(1-a,v)\Big
)\nonumber
	\\
	&=&\frac{t^{-\alpha}}{\alpha}\Big(-\exp(-v)-v^a\gamma(1-a,v)\Big),
	\IEEEeqnarraynumspace
	\label{eq:indef_integr_symm} 
\end{IEEEeqnarray} 
where $\gamma(s,x)$ is the lower incomplete gamma function \eqref{eq:lower_inc_gamma} and $a=\alpha/2$, as in \eqref{eq:change_var}. 

We evaluate the indefinite integral  
\eqref{eq:indef_integr_symm} 
in the upper limit
$t = c^{-1/\alpha}$  and $v = v_c 
:=\sigma_W^2s^2c^{-2/\alpha}/2$ and in the lower limit
$t = d^{-1/\alpha}$ and $v = v_d:=\sigma_W^2s^2d^{-2/\alpha}/2$, 
respectively, as in \eqref{eq:cf_Rcd_integral_2}; the
definite integral then becomes,
\begin{IEEEeqnarray*}{rCl}
	& &\frac{c}{\alpha}\left(-\exp(-v_c)-v_c^{a}\gamma(1-a,v_c)\right)
	-\frac{d}{\alpha}\left(-\exp(-v_d)-v_d^{a}\gamma(1-a,v_d)\right).
\end{IEEEeqnarray*} 
\noindent
Then, as $d\rightarrow \infty$, 
$v_d\rightarrow0$ and  $v_d^a \gamma(1-a, v_d)\rightarrow0$, 
hence the definite integral,
\begin{IEEEeqnarray*}{rCl}
	\to\frac{c}{\alpha}\left(-\exp(-v_c)-v_c^{a}\gamma(1-a,v_c)\right)+ \frac{d}{\alpha},
\end{IEEEeqnarray*} 
and substituting in \eqref{eq:cf_Rcd_integral_1}, 
recalling that $B=0$ when $\mu_W =0$, we have,
\begin{IEEEeqnarray}{rCl}
	\phi_{R_{(c,\infty)}}(s) &=&\omega_{R_{(c,\infty)}}(v_c)\nonumber
	\\
	&=& \exp(c(1-\exp(-v_c)-v_c^{a}\gamma(1-a,v_c))).
	\IEEEeqnarraynumspace
	\label{eq:phi_R_c_inf_symm}
\end{IEEEeqnarray} 
Notice that $v_c$ corresponds to $u$ in \eqref{eq:change_var}, in view of \eqref{eq:S_c_infty} and the fact that $\mathbb{E}[W^2_1] = \sigma_W^2$.
Then the first part of the statement follows by recalling that, from \eqref{eq:res:unnormalized:exact}, $\phi_{Z_{(c, \infty)}} (s)= \phi_{R_{(c, \infty)}}(s/S_{(c,\infty)})$.  
Finally, using the change of variables~\eqref{eq:change_var}, $\phi_{R_{(c, \infty)}} (s/S_{(c, \infty)})=\omega_{R_{(c, \infty)}} (u/S^2_{(c,\infty)}) = \psi_{R_{(c, \infty)}} (w)$.

\subsection{CF of $X_{(0,c)}$}
To prove the second part of the lemma, we
evaluate the indefinite integral~\eqref{eq:indef_integr_symm} in the 
upper limit $t = \infty$ and $v = \infty$, and in the lower limit
$t = c^{-1/\alpha}$ and $v = v_c =\sigma_W^2s^2c^{-2/\alpha}/2$, 
as in \eqref{eq:phi_X_0c_gauss}. Recalling from \eqref{eq:ch_var_t_s} that
$t^{-\alpha} v^a = (|s| \sigma_W /\sqrt{2})^\alpha$, where $\alpha = 2a$ as in  \eqref{eq:change_var}, and also recalling the definition of the gamma 
function~\eqref{eq:gamma}, the corresponding definite integral becomes,
{{\small 
		\begin{IEEEeqnarray*}{rCl}
			-\frac{1}{\alpha}\left( \frac{|s|\sigma_W}{\sqrt{2}}\right) ^\alpha \Gamma(1-a)
			-\frac{c}{\alpha}\left(-\exp(-v_c)-v_c^{a}\gamma(1-a,v_c)\right),
		\end{IEEEeqnarray*} 
		\par}}
\noindent
{{\normalsize
		and, substituting in  \eqref{eq:phi_X_0c_gauss} and recalling the definition of $v_c$, we obtain,
		\begin{IEEEeqnarray*}{rCl}
			\phi_{X_{(0,c)}}(s) & = & \omega_{X_{(0,c)}}(v_c)
			\\  
			&=&\exp(-c(1-\exp(-v_c)-u_c^a\gamma(1-a,v_c)+v_c^{a}\Gamma(1-a)))
			\\&=&\exp(-c(1-\exp(-v_c)+v_c^a\Gamma(1-a,v_c))), 
		\end{IEEEeqnarray*} 
		as claimed.
}}

\section{Properties of gamma functions}
\label{appendix:properties_gamma_functions}

\subsection{Inequalities for $\gamma(s,x)$}
Recall the definition
of the lower incomplete gamma function in \eqref{eq:lower_inc_gamma}. 

\medskip

\begin{lemma}
	{\em \cite[Theorem~4.1]{neuman:13}}
	\label{lem:gammaUB}
	For all $x>0$ and $s\in(0,1]$
	\begin{IEEEeqnarray*}{c} 
		\gamma(s,x)\leq\frac{x^s}{s(s+1)}(1+se^{-x}).
	\end{IEEEeqnarray*}
\end{lemma}

\begin{lemma}
	{\em \cite[Ineq.~(8.10.2)]{NIST:DLMF}}
	\label{lem:gammaLB}
	For all $x>0$ and $s>0$
	\begin{IEEEeqnarray}{c} 
		\gamma(s,x)\geq\frac{x^{s-1}}{s}(1-e^{-x}). 
		\label{eq:gammaLB}	
	\end{IEEEeqnarray}
\end{lemma}
\noindent
Combining Lemmas~\ref{lem:gammaUB} and \ref{lem:gammaLB}
we obtain that
\begin{IEEEeqnarray}{c} 
	\lim_{x\downarrow 0}\frac{\gamma(s,x)}{x^s}=\frac{1}{s},
	\label{eq:gammaAS}
\end{IEEEeqnarray}
and also the following:

\medskip

\begin{lemma}
	\label{lem:gammaDOUBLE}
	For all $x>0$ and $s\in(0,1]$:
	\begin{IEEEeqnarray*}{c} 
		-\frac{x}{1+s}
		\leq
		e^{-x}-\frac{s\gamma(s,x)}{x^s}
		\leq
		-\frac{x(1-x)}{2}.
	\end{IEEEeqnarray*}
\end{lemma}

\begin{IEEEproof}
	The bounds in 
	Lemmas~\ref{lem:gammaUB} and~\ref{lem:gammaLB}
	immediately give,
	\begin{IEEEeqnarray*}{c} 
		\frac{1-e^{-x}}{x}
		\leq
		\frac{s\gamma(s,x)}{x^s}
		\leq
		\frac{1+se^{-x}}{1+s},
	\end{IEEEeqnarray*}
	and subtracting all three sides from $e^{-x}$
	and simplifying,
	\begin{IEEEeqnarray}{c} 
		\frac{e^{-x}-1}{1+s}
		\leq
		e^{-x}-\frac{s\gamma(s,x)}{x^s}
		\leq
		\frac{(x+1)e^{-x}-1}{x}.
		\label{eq:prelem}
		\IEEEeqnarraynumspace
	\end{IEEEeqnarray}
	Applying the elementary inequality $e^{-x}\geq 1-x$,
	$x\geq 0$,
	to the left-hand side gives the lower bound 
	in the statement, and similarly applying the
	inequality $e^{-x}\leq 1-x+\frac{x^2}{2}$,
	$x\geq 0$, to the right-hand side gives the 
	corresponding upper bound.
\end{IEEEproof}


\subsection{Inequalities for $\Gamma(s,x)$}
Recall the definition 
of the upper incomplete gamma function in \eqref{eq:upper_inc_gamma}. 

\medskip

\begin{lemma}
	\label{lem:GammaRepr}
	{\em \cite[Eq.~(8.6.7)]{NIST:DLMF}}
	For all $s>0,x>0$, $\Gamma(s,x)$ admits the representation:
	\begin{IEEEeqnarray*}{c} 
		\Gamma(s,x)=x^s
		\int_0^\infty \exp(st-xe^t)\dif t.
	\end{IEEEeqnarray*}
\end{lemma}

\begin{lemma}
	\label{lem:GammaUB}
	For all $x>s>0$:
	\begin{IEEEeqnarray*}{c} 
		\frac{\Gamma(s,x)}{x^se^{-x}}\leq\frac{1}{x-s}.
	\end{IEEEeqnarray*}
\end{lemma}

\begin{IEEEproof}
	From Lemma~\ref{lem:GammaRepr},
	\begin{IEEEeqnarray}{c} 
		\frac{\Gamma(s,x)}{x^se^{-x}}
		=
		\int_0^\infty \exp(st-xe^t+x) \dif t,
		\label{eq:Gammarep}
		\IEEEeqnarraynumspace
	\end{IEEEeqnarray}
	and since $e^t\geq t+1$ for $t\geq 0$,	
	\begin{IEEEeqnarray*}{c} 
		\frac{\Gamma(s,x)}{x^se^{-x}}
		\leq
		\int_0^\infty \exp(-t(x-s))\dif t=\frac{1}{x-s}.
	\end{IEEEeqnarray*}
\end{IEEEproof}

The bound in the last lemma applies 
when $x>s$. When
$x$ may be smaller
than $s$ we have a somewhat weaker bound,
but this time uniformly in $x$:

\medskip

\begin{lemma}
	\label{lem:GammaUB2}
	Let $\delta\in(0,1)$ be arbitrary.
	For all $s>0$ and all $x\geq \delta s$:
	\begin{IEEEeqnarray*}{c} 
		\frac{\Gamma(s,x)}{x^se^{-x}}\leq
		\sqrt{\frac{2\pi}{\delta s}}
		\exp\left[ \frac{\delta s}{2}(\delta^{-1}-1)^2\right] .
	\end{IEEEeqnarray*}
\end{lemma}

\begin{IEEEproof}
	Starting again with~(\ref{eq:Gammarep})
	and noting that the integrand is decreasing in $x$,
	we have that, for $x\geq\delta s$,
	\begin{IEEEeqnarray*}{c} 
		\frac{\Gamma(s,x)}{x^se^{-x}}
		\leq
		\int_0^\infty \exp\left[ s(t-\delta e^t+\delta)\right] \dif t,
	\end{IEEEeqnarray*}
	and since $e^t\geq 1+t+\frac{t^2}{2}$ for $t\geq 0$,
	\begin{IEEEeqnarray*}{c} 
		\frac{\Gamma(s,x)}{x^se^{-x}}
		\leq
		\int_0^\infty \exp\left\{-\frac{\delta s}{2}\Big[t^2-2t(
		\delta^{-1}-1)\Big]\right\}\dif t.
	\end{IEEEeqnarray*}
	Completing the square in the exponent,
	{\small
		\begin{IEEEeqnarray*}{rCl}
			\frac{\Gamma(s,x)}{x^se^{-x}}
			& \leq &	 
			\exp\left[ \frac{\delta s}{2}(\delta^{-1}-1)^2\right]  
			\int_0^\infty \exp\left\{-\frac{\delta s}{2}\Big[t-(
			\delta^{-1}-1)\Big]^2\right\}\dif t\\
			&\leq&
			\sqrt{\frac{2\pi}{\delta s}}
			\exp\left[ \frac{\delta s}{2}\Big(\frac{1}{\delta}-1\Big)^2\right] 
			F_Z\left( (\delta^{-1} -1)(\delta s)^{1/2}\right), 	
		\end{IEEEeqnarray*}
		\par}
	\noindent
	where $F_Z(\cdot)$ denotes the standard normal CDF, which is,
	of course, no greater than~1.
\end{IEEEproof}


\subsection{Asymptotics and derivatives} 
\label{subsect:inc_gamma_asymptotic}

When $x\in \mathbb{R}$, the upper incomplete gamma function \eqref{eq:upper_inc_gamma} has the following asymptotic behaviour for $x\rightarrow \infty$,
\begin{IEEEeqnarray}{rCl}
	\Gamma(s,x) \sim  x^{s-1} e^{-x}.
	\label{eq:Gamma_asymptotic}
\end{IEEEeqnarray}
This can be proved using series expansions of $\Gamma(s,x)$, 
see \cite{AbramowitzStegun1964} or \cite{NIST:DLMF}. 


Finally, by the fundamental theorem of calculus and the definitions \eqref{eq:lower_inc_gamma} and \eqref{eq:upper_inc_gamma}, we have,
\begin{IEEEeqnarray}{rCl}
	\frac{\dif \gamma(s, x)}{\dif x }
	= - 	\frac{\dif \Gamma(s, x)}{\dif x } = x^{s-1} e^{-x}. 
	\label{eq:inc_gamma_derivatives}  
\end{IEEEeqnarray}

\section{Preliminary results for the proofs of Section~\ref{section:residual_convergence} and \ref{section:residual_contribution}}\label{appendix:preliminary:res:convergence}


Here we list a number of auxiliary lemmas that will be needed in the
proofs of Theorems~\ref{thm:linear} and~\ref{thm:root}. The lemma below,
stated without proof, is a simple calculus exercise. 

\medskip

\begin{lemma}
	\label{lem:feller}
	If $0\leq x\leq C$ and $0\leq y\leq C$,
	then for any (not necessarily integer) $n\geq 1$,
	\begin{IEEEeqnarray*}{c} 
		|x^n-y^n|\leq n|x-y|C^{n-1}.		
	\end{IEEEeqnarray*}
\end{lemma}



\subsection{The function $g(w)$} \label{subsection:g_properties}


When $W_j \sim \mathcal{N}(0, \sigma^2_W)$, 
the CF of the PSR standardized residual
is given in \eqref{eq:asympt_cf_res_symm}.
For convenience we write,
$\psi_{Z_{(\c,\infty)}}(\w)= f(\w)^{\c}$, with,
\begin{IEEEeqnarray}{c} 
	f(\w) := \exp 
	\left(1-e^{-\w}-\w^a\gamma\left(1-a,\w\right)
	\right),
	\label{eq:fdefn}
\end{IEEEeqnarray} 
and we also define,
\begin{IEEEeqnarray}{c} 
	g(\w):=\log f(\w)=1-e^{-\w}-\w^a\gamma(1-a,\w). 
	\label{eq:gdefn}
	\IEEEeqnarraynumspace
\end{IEEEeqnarray}

First note that $g(0)=0$
and that, from Lemma~\ref{lem:gammaLB}, we have
$g(\w)\leq 0$ for all $\w\geq 0$ and in fact,
\begin{IEEEeqnarray}{c} 
	g(\w)\leq-\frac{1-e^{-\w}}{\eta}
	\leq
	-\frac{w}{\eta}\left(1-\frac{w}{2}\right),
	\label{eq:gUB1}
\end{IEEEeqnarray}
where the second step follows 
from the fact that
$e^{-x}\leq 1-x+x^2/2$, for $x\geq 0$.
Therefore,
\begin{IEEEeqnarray}{c} 
	g(\w)
	\leq -\frac{\w}{2\eta},
	\;\;\;0\leq w\leq 1.
	\label{eq:gUB} 
\end{IEEEeqnarray}
For $\w\geq 1$ we have the bound:

\medskip

\begin{lemma}
	\label{lem:finer}
	For all $\w\geq 1$,
	\begin{IEEEeqnarray}{c} 
		g(\w)
		\leq 1-e^{-1}-\bar{\gamma}(a)\w^a,
		\label{eq:gUB2}
	\end{IEEEeqnarray}
	where $\bar{\gamma}(a)=\gamma(1-a,1).$
\end{lemma}

\begin{IEEEproof}
	The statement is equivalent to,
	\begin{IEEEeqnarray*}{c} 
		H(\w):=\gamma(1-a,\w)-\gamma(1-a,1)-\frac{e^{-1}-e^{-\w}}{\w^a}\geq 0,
	\end{IEEEeqnarray*}
	for all $\w\geq 1$.
	Since $H(1)=0$ and, using \eqref{eq:inc_gamma_derivatives}, the derivative,
	\begin{IEEEeqnarray*}{c} 
		H'(\w)=a\w^{-a-1}(e^{-1}-e^{-\w}),
	\end{IEEEeqnarray*}
	is always nonnegative,  the result follows.
\end{IEEEproof}

\medskip

Differentiating,
\begin{IEEEeqnarray}{c} 
	g'(\w)=-a\w^{a-1}\gamma(1-a,\w),
	\label{eq:gprime1}
	\IEEEeqnarraynumspace
\end{IEEEeqnarray}
so that
$g'(0+)=-1/\eta$ 
by~(\ref{eq:gammaAS}),
and,
\begin{IEEEeqnarray}{c} 
	0\geq g'(\w)\geq
	-\frac{1+(1-a)e^{-\w}}{\eta(2-a)}\geq -\frac{1}{\eta},
	\label{eq:gprime}
	\IEEEeqnarraynumspace
\end{IEEEeqnarray}
for all $\w\geq 0$, where the first inequality
is obvious by~(\ref{eq:gprime1}),
the second follows 
from Lemma~\ref{lem:gammaUB},
and the third from the fact
that $e^{-\w}\leq 1$ always.
In particular, this implies that,
\begin{IEEEeqnarray}{c} 
	g(\w)\geq -\frac{\w}{\eta},\;\;\;\w\geq 0.
	\label{eq:gLB}
\end{IEEEeqnarray}
Differentiating again,
\begin{IEEEeqnarray*}{c} 
	g''(\w)=-\frac{a}{\w}\left[e^{-\w}-(1-a)\w^{a-1}\gamma(1-a,\w)\right],
\end{IEEEeqnarray*}
and applying Lemma~\ref{lem:gammaDOUBLE}, we have that,
\begin{IEEEeqnarray}{c} 
	-a \leq \frac{a(1-\w)}{2}\leq g''(\w)\leq\frac{a}{2-a},
	\label{eq:g2prime}
	\IEEEeqnarraynumspace
\end{IEEEeqnarray}
where the first inequality only holds for $0\leq \w\leq 3$;
on the other hand, since~(\ref{eq:prelem}) holds for all
$\w\geq 0$ and the function $x\mapsto(x+1)e^{-x}$ is 
decreasing, 
for all $w\geq 1$ we have the simple lower bound,
\begin{IEEEeqnarray}{c} 
	g''(\w)
	\geq
	-\frac{a}{\w}\left(\frac{(\w+1)e^{-\w}-1}{\w}\right)
	\geq
	\frac{a(1-2e^{-1})}{\w^2}.
	\label{eq:g2prime2}
	\IEEEeqnarraynumspace
\end{IEEEeqnarray}
In particular, from~(\ref{eq:g2prime}) and~(\ref{eq:g2prime2})
it follows that,
\begin{IEEEeqnarray}{c} 
	0\leq g''(\w)\leq \frac{a}{2-a},\;\;\;\mbox{for all}\;\w\geq 0.
	\label{eq:g2prime3}
	\IEEEeqnarraynumspace
\end{IEEEeqnarray}
This implies that the function $g(w)$ is {convex}, and a tighter upper 
bound than \eqref{eq:gUB} for $w \in [0,1]$ is the linear interpolation 
at the extremes,
\begin{IEEEeqnarray}{c} 
	g(\w)
	\leq \barg w,
	\;\;\;0\leq w\leq 1,
	\label{eq:gUB_improved} 
\end{IEEEeqnarray}
where $\barg= g(1)$, as in \eqref{eq:bar_g}.
More generally, the linear interpolating upper bound for $g(w)$ on  $ [0, \bar{w}]$ is 
\begin{IEEEeqnarray}{c} 
	g(\w)
	\leq \frac{g(\bar{w})}{\bar{w}} w,
	\;\;\;0\leq w\leq \bar{w}.
	\label{eq:gUB_linear} 
\end{IEEEeqnarray}

%

\subsection{The function $q(u)$} \label{subsection:q_properties}


When $W_j \sim \mathcal{N}(0, \sigma^2_W)$, 
the CF of the truncated PSR is given in~\eqref{eq:cf_X0c_symm}.
For convenience we write 
$\omega_{X_{(0,c)}}(u)= r(u)^{\c}$, with,
\begin{IEEEeqnarray*}{c} 
	r(u) := \exp(-(1-e^{-u}+u^a\Gamma(1-a,u))),
\end{IEEEeqnarray*} 
and we also define:
\begin{IEEEeqnarray*}{c} 
	q(u):=\log r(u)=-(1-e^{-u}+u^a\Gamma(1-a,u)). 
	\IEEEeqnarraynumspace
\end{IEEEeqnarray*}

\noindent
It is clear that $q(0)=0$ and  $q(u)<0$ for $u>0$.  
In fact, using~\eqref{eq:inc_gamma_derivatives}, 
\begin{IEEEeqnarray*}{rCl} 
	q'(u)=-au^{a-1}\Gamma(1-a,u)\leq 0,
\end{IEEEeqnarray*}	
hence $q(u)$ is monotonically decreasing.
Moreover, by \eqref{eq:Gamma_asymptotic}, $q(u)$ is asymptotic 
to $-1$ at $u\rightarrow \infty$. 

Since $a-1<0$ and $\Gamma(1-a,u)$ is decreasing in $u$ we have that $q'(u)$ 
is  monotonically increasing towards zero at $u\rightarrow \infty$,
since $\Gamma(1-a,u)$ tends to $\Gamma(1-a)$ as $u\rightarrow 0$. 
Note also that 
$q'(0)$ diverges to $-\infty$, which prevents a Taylor expansion around zero. 

Differentiating once again,
\begin{IEEEeqnarray*}{rCl} 
	q''(u)=a((1-a)u^{a-2}\Gamma(1-a,u)+u^{-1}e^{-u})\geq 0,
\end{IEEEeqnarray*}	
so $q''$ is also monotonically decreasing,
it decreases to zero at $u \rightarrow \infty$, and it
diverges to $+\infty$ at $u=0$.

Therefore, $\log(\omega_{X_{(0,c)}}(u)) = c q(u) $ is convex,
which means it can easily be bounded above by line segments.
Indeed, the following construction will be useful in
the proof of Theorem~\ref{theorem:X_hat_symm}:
We will employ a piecewise linear interpolating bound 
with $N$ segments for $u \in [0,1]$, and a constant bound for $u>1$. 
When $N = 1$, the bound is simply,
\begin{IEEEeqnarray*}{rCl} 
	\log(\omega_{X_{(0,c)}}(u)) = cq(u)\leq L^1(u)
	&:=&
	\begin{cases}
		-cu((1-\exp(-1))+\Gamma(1-a,1)),
		& u \in[0,1],
		\\
		-c((1-\exp(-1))+\Gamma(1-a,1)),
		& u>1,
	\end{cases}
	\\
	&=&
	\begin{cases}
		m_0 u,& u \in[0,1],
		\\
		k_{(1,\infty)},& u>1,
	\end{cases}
\end{IEEEeqnarray*}	
where,
\begin{IEEEeqnarray*}{rCl} 
	m_0 =  k_{(1,\infty)} & := &   -c((1-\exp(-1))+\Gamma(1-a,1))
	\\
	& =& 
	\log(\omega_{X_{(0,c)}}(1)) 
	\\
	& =& 
	\log(cq(1)) 
	<0. 
\end{IEEEeqnarray*}	
Bounding with more than one line segment in the domain~$[0,1]$ 
turns out to be important for tightening the bound on $\Delta(X, \hat{X})$ 
and for capturing its dependence on $\alpha$,
as observed via numerical integration results; cf. 
Section \ref{section:residual_contribution} and 
Appendix \ref{app:proof_hat_X_symm}. For $N\geq 1$, 
we select $N$ points $u_i$ in $[0,1]$,
$$0 =: u_0 < u_1 < \ldots < u_N := 1,$$ 
and the 
respective values of $\log\omega_{X_{(0,c)}}$,
\begin{IEEEeqnarray*}{rCl} 
	f_0: = 0 >  f_1 := \log(\omega_{X_{(0,c)}}(u_1))  
	> \ldots > f_N := \log(\omega_{X_{(0,c)}}(u_N)). 
\end{IEEEeqnarray*}	
The equation of the $i$-th line segment, for $i = 0, \ldots, N-1$ is $y = (m_i u + q_i)\mathds{1}_{A_i}{(u)}$, with $A_i :=[u_i, u_{i+1}]$ and,
\begin{IEEEeqnarray*}{rCl} 
	m_i &:=& \frac{f_{i+1} - f_i}{u_{i+1} - u_i},
	\qquad
	q_i := - m_i u_i + f_i. 	
\end{IEEEeqnarray*}	
The general upper bound on $\log(\omega_{X_{(0,c)}}(u))$ then becomes,
\begin{IEEEeqnarray}{rCl} 
	\log(\omega_{X_{(0,c)}}(u)) \leq L^N(u)
	&:=&
	\begin{cases}
		\sum_{i = 0}^{N-1}  (m_i u + q_i)\mathds{1}_{A_i}{(u)} ,& u \in[0,1],
		\\
		k_{(1,\infty)},& u>1,
	\end{cases}
	\label{eq:bound_log_phi_X0c_piece_lin}
	\IEEEeqnarraynumspace
\end{IEEEeqnarray}	
with
$k_{(1,\infty)}$ 
as above. 


\section{Proof of Theorem \ref{thm:linear}}\label{appendix:proofs_th_lin}


\textit{Step I.} 
Let $W_j \sim \mathcal{N}(0, \sigma^2_W)$
so that the CF of the PSR standardized residual $Z_{(c,\infty)}$ 
is given in \eqref{eq:asympt_cf_res_symm},
and write
$\psi_{Z_{(\c,\infty)}}(\w)= f(\w)^{\c}$, 
with $f(w)$ defined as in \eqref{eq:fdefn}. 
We apply the smoothing lemma \eqref{eq:int_berry_finite} to $	\Delta(Z_{(c, \infty)},Z)$ as in~\eqref{eq:kolmog_residual}.
Given that the standard Gaussian PDF is uniformly bounded by 
$m:=1/\sqrt{2\pi}<2/5$, 
for any 
$\Theta>0$ equation~\eqref{eq:int_berry_finite} gives,
\begin{IEEEeqnarray}{c} 
	\pi\Delta(Z_{(c, \infty)},Z)
	\leq
	\int_{-\Theta}^\Theta|\phi_{Z_{(c, \infty)}}(s)-\phi_{Z}(s)|\frac{1}{|s|}\dif s
	+\frac{9.6}{\Theta}, 
	\label{eq:essen_res_s}
	\IEEEeqnarraynumspace
\end{IEEEeqnarray}
with $\phi_Z(s)$ as in \eqref{eq:Z:lim_dc}.  
Letting $\Theta\to\infty$ and changing variables as in
\eqref{eq:change_var} we have,
\begin{IEEEeqnarray}{rCl} 
	\pi\Delta(Z_{(c, \infty)},Z)
	&\leq& 
	\int_{0}^{\infty}
	\big|f(\w)^\c
	-\phi_Z\big(
	\sqrt{{2\w}/{\eta}}
	\big)^\c\big|
	\frac{1}{\w}
	\dif \w \nonumber\\
	&=& \pi\bar{I}(Z_{(c, \infty)},Z). 
	\label{eq:esseen}
	\IEEEeqnarraynumspace
\end{IEEEeqnarray}

\noindent
\textit{Step II.} 
We 
apply Lemma~\ref{lem:feller} to the 
integrand in~(\ref{eq:esseen}), 
with
$x=f(\w)$, $y=\phi_Z(\sqrt{2\w/\eta}) = \exp(-w/\eta)$ by \eqref{eq:cf_std_gauss_w}, and $n = \c$. 
Let $g(w) = \log(f(w))$ as in \eqref{eq:gdefn}, and  $\bar{\gamma}(a)$ and $\barg$ as in \eqref{eq:bar_gamma_a} and \eqref{eq:bar_g}, respectively,  and define,
\begin{IEEEeqnarray}{c} 
	h(\w):=\begin{cases}	
		-\barg w,&\;\w\in[0,1],\\
		e^{-1}-1+\bar{\gamma}(a)\w^a,&
		\;w>1.
	\end{cases}
	\label{eq:hdef}
\end{IEEEeqnarray}
From~(\ref{eq:gUB_improved}) and~(\ref{eq:gUB2}) it follows
that $x\leq \exp(-h(\w)) $ for all $\w\geq 0$.
Since, by \eqref{eq:gLB}, we have $y\leq x$, for all $w \geq 0$, 
it follows that $y\leq \exp(-h(\w))$, for all $w \geq 0$.
Therefore, we can take,
$C=\exp(-h(\w))$ in Lemma~\ref{lem:feller},
and substituting the resulting bound in~(\ref{eq:esseen}) 
we get,
\begin{IEEEeqnarray}{rCl} 
	\pi\Delta(Z_{(c, \infty)},Z)
	& \leq &
	\int_0^{\infty} \frac{\c}{\w}
	\left|f(\w)-\phi_Z\left(\sqrt{\frac{2\w}{\eta}}\right)\right|
	e^{-(\c-1)h(\w)}\dif \w.
	\IEEEeqnarraynumspace
	\label{eq:esseen2}
\end{IEEEeqnarray}

\noindent
\textit{Step III.} 
In order to bound the absolute difference
in the integrand in \eqref{eq:esseen2}, 
we write $f$ 
as a quadratic Taylor expansion. Noting that 
$f'=g'\exp(g)$ and $f''=(g''+g'^2)\exp(g)$, where $g$ is defined in \eqref{eq:gdefn}, and
recalling that $g(0)=0$ and $g'(0)=-1/\eta$,
we have,
\begin{IEEEeqnarray*}{rCl} 
	\left|f(\w)-1+\frac{\w}{\eta}\right|
	&\leq&
	\frac{\w^2}{2}
	\sup_{v\geq 0}
	\left[\left(|g''(v)|+g'(v)^2\right)e^{g(v)}\right]\\
	&\leq&
	\frac{\w^2}{2}
	\left[\frac{a}{2-a}+\frac{1}{\eta^2}\right],
\end{IEEEeqnarray*}
where in the second inequality we used~(\ref{eq:gprime}),~(\ref{eq:g2prime3}),
and the fact that $g(w)\leq 0$ for all $w\geq 0$.
From the standard quadratic expansion for the exponential
function we similarly have,
\begin{IEEEeqnarray*}{c} 
	\left|\exp\left(-\frac{\w}{\eta}\right) -1+\frac{\w}{\eta}\right|\leq
	\frac{\w^2}{2\eta^2}.
\end{IEEEeqnarray*}
Combining the last two bounds,
\begin{IEEEeqnarray}{c} 
	\left|f(\w)-\phi_Z\left(\sqrt{\frac{2\w}{\eta}}\right)\right|
	\leq
	\frac{\w^2}{2}
	\left[\frac{a}{2-a}+\frac{2}{\eta^2}\right].
	\label{eq:expansion}
	\IEEEeqnarraynumspace
\end{IEEEeqnarray}

\noindent
\textit{Step IV.}  
Finally, substituting \eqref{eq:expansion} in~(\ref{eq:esseen2}),
\begin{IEEEeqnarray}{rCl} 
	\Delta(Z_{(c, \infty)},Z)
	&\leq&
	\frac{c}{\pi}
	\left[\frac{a}{2(2-a)}+\frac{1}{\eta^2}\right]
	\int_0^{\infty} 
	\w
	\exp\big(-(\c-1)h(\w)\big)\dif \w
	\nonumber\\
	&=&
	cK(a)
	[I^Z(\c)+J^Z(\c)],
	\label{eq:IJ}
	\IEEEeqnarraynumspace
\end{IEEEeqnarray}
where $K(c)$ is defined in \eqref{eq:K_c_a}, and where we first integrate over $[0,1]$,

\begin{IEEEeqnarray}{rCl} 
	I^Z(\c)&:=&\int_0^1 \w \exp\big(\w(\c-1)\barg\big)\dif \w \nonumber 	\\
	& =& \frac{1}{(\c-1)^2 \barg^2} 
	+\frac{1}{(\c-1)\barg}\left(1-\frac{1}{(\c-1)\barg}\right) \exp\left((\c-1)\barg\right), 
	\IEEEeqnarraynumspace
	\label{eq:I}
\end{IEEEeqnarray}
\noindent
and then over $[1,\infty)$,
\begin{IEEEeqnarray}{rCl}
	J^Z(\c)
	&:=&
	e^{(1-e^{-1})(\c-1)} 
	\int_{1}^{\infty} 
	\w
	\exp\Big\{-(\c-1)\bar{\gamma}(a)\w^a\Big\}\dif \w\nonumber\\
	&=&
	\frac{e^{(1-e^{-1})(\c-1)}}
	{a [(\c-1)\bar{\gamma}(a)]^{2/a}}
	\int_{(\c-1)\bar{\gamma}(a)}^{\infty}
	u^{\frac{2}{a}-1}
	e^{-u}
	\dif u\nonumber\\
	&=&
	e^{(1-e^{-1})(\c-1)}
	\frac{ \Gamma(2/a, (\c-1)\bar{\gamma}(a)) }
	{a[(\c-1)\bar{\gamma}(a)]^{2/a}},
	\label{eq:Jfirst}
	\IEEEeqnarraynumspace
\end{IEEEeqnarray}
where
$\Gamma(s,x)$
is the upper incomplete gamma function defined in \eqref{eq:upper_inc_gamma}. 
Combining~(\ref{eq:IJ}),~(\ref{eq:I}) and~(\ref{eq:Jfirst}) 
gives the claimed result, upon re-writing,
\begin{IEEEeqnarray}{rCl}
	B_1(c,\alpha) 
	&=&
	\left(cK(a) (I^Z(c)  + J^Z(c)) \right) \nonumber
	\\
	&=&   \left(\frac{cK(a)}{c-1} \big((c-1)I^Z(c)  + (c-1)J^Z(c)\big) \right).
	\IEEEeqnarraynumspace
	\label{eq:B1_split}
\end{IEEEeqnarray}

\section{Proof of Corollary~\ref{cor_alternative}}\label{appendix:proof:corollary:lin}
Clearly, it suffices to prove the first asymptotic assertion of the
corollary. To that end, we examine each of the two terms in
the expression \eqref{eq:B1_split} for $B_1(c, \alpha)$.
First we note that,
		\begin{IEEEeqnarray*}{rCl}
			(c-1)I^Z(c) 
			& = &		\frac{1}{(\c-1) \barg^2}
			+ \left(\frac{1}{\barg}-\frac{1}{(\c-1) \barg^2}\right)
			\exp\left((\c-1)\barg\right) \\
			& \sim &\frac{1}{(c-1)\barg^2},
		\end{IEEEeqnarray*}
		\par
since $\barg <0$.
And also,
\begin{IEEEeqnarray*}{rCl}
	(c-1)J^Z(c) 
	& = &  
	\frac{
		(\c-1)\exp\left\{(\c-1)
		(1-e^{-1})
		\right\}
	}
	{a[(\c-1)\bar{\gamma}(a)]^{2/a}} 
	\Gamma\left({2}/{a}, (\c-1)\bar{\gamma}(a)\right) 
	\\
	& \sim &
	\frac{
		(\c-1)\exp\left\{(\c-1)
		(1-e^{-1})
		\right\}
	}
	{a[(\c-1)\bar{\gamma}(a)]^{2/a}} 
	\times
	[(\c-1)\bar{\gamma}(a)]^{2/a - 1} 
	\exp(-(c-1) \bar{\gamma}(a))
	\\
	& = &
	\frac{1}{a \bar{\gamma}(a) }
	\exp\left\{(\c-1)
	(1-e^{-1} - \bar{\gamma}(a))
	\right\},
\end{IEEEeqnarray*}
\noindent
where we used the asymptotic property~\eqref{eq:Gamma_asymptotic}.
Observe that the exponent in the last expression is negative: 
Using~\eqref{eq:gammaLB} and the fact that
$a\in(0,1)$, we have,
\begin{IEEEeqnarray*}{rCl}
	1 - e^{-1} - \bar{\gamma}(a) < (1 - e^{-1}) \left( 1 - \frac{1}{1-a}\right) <0.
\end{IEEEeqnarray*}
Therefore, combining these two estimates with~\eqref{eq:B1_split},
we indeed have,
		\begin{IEEEeqnarray*}{rCl}
			B_1(c, \alpha)& \sim &
			K(a)
			\frac{1}{(c-1)g^2(1)}.
		\end{IEEEeqnarray*}


\section{Proof of Theorem~\ref{thm:root}}\label{appendix_proof_th_sqrt}

As in the proof of Theorem \ref{thm:linear},
we start from \eqref{eq:essen_res_s}
and  perform the change of variables \eqref{eq:change_var}
to obtain that, for any $\Theta>0$,
\begin{IEEEeqnarray}{rCl} 
	\pi\Delta(Z_{(c, \infty)},Z)
	& \leq & 
	\int_{0}^{\eta \Theta^2/2\c}
	\left|f(w)^\c
	-\phi_Z\left(
	\sqrt{\frac{2w}{\eta}}
	\right)^\c\right|
	\frac{1}{w}
	\dif w +  \frac{9.6}{\Theta}.
	\label{eq:esseenR}
	\IEEEeqnarraynumspace
\end{IEEEeqnarray}
\noindent
But here, instead of letting $\Theta \rightarrow \infty$, we choose
$\Theta=\sqrt{2(2-\delta)\c/\eta}$, and 
we apply Lemma~\ref{lem:feller} to the 
integrand in~(\ref{eq:esseenR}), with
$x=f(\w)$ and $y=\phi_Z(\sqrt{2\w/\eta}) =\exp(-\w/\eta)$.
From \eqref{eq:gLB} $y\leq x$, for all $w \geq 0$,
and using~(\ref{eq:gUB_linear}) with $\bar{w} = 2 - \delta$
we have that $x\leq\exp\big(\frac{g(2-\delta)}{2 - \delta} w\big)$.  We can then take $C=\exp\big(\frac{g(2-\delta)}{2 - \delta} w\big)$ in 
Lemma~\ref{lem:feller} applied to~(\ref{eq:esseenR}),
to obtain,
\begin{IEEEeqnarray}{l} 
	\pi\Delta(Z_{(c, \infty)},Z) -\frac{9.6\sqrt{\eta}}{\sqrt{2(2-\delta)\c}} 
	\leq 
	\int_0^{2-\delta} \frac{\c}{\w}
	\left|f(\w)-\phi_Z\left(\sqrt{\frac{2\w}{\eta}}\right)\right|
	\exp\left(\frac{ (\c-1)g(2-\delta)\w}{2-\delta}\right)\dif w. 
	\label{eq:esseen2R}
	\IEEEeqnarraynumspace
\end{IEEEeqnarray}
\noindent
Recalling the earlier expansion~(\ref{eq:expansion}),
substituting in~(\ref{eq:esseen2R}), and integrating,
yields,
\begin{IEEEeqnarray*}{rCl} 
	\pi\Delta(Z_{(c, \infty)},Z) -\frac{9.6\sqrt{\eta}}{\sqrt{2(2-\delta)\c}}
	&\leq& \c
	\left[\frac{a}{2(2-a)}+\frac{1}{\eta^2}\right]\int_0^{2-\delta}\w\exp\left(\frac{(c-1) g(2-\delta) w}{2-\delta}\right)\dif\w
	\\
	&= &
	\frac{1}{\c}
	\left[\frac{a}{2(2-a)}+\frac{1}{\eta^2}\right]
	\left(\frac{c (2-\delta) }{(c-1) g(2-\delta)}\right)^2 \times 
	\\
	& &
	\times 
	\left\{1-\left[1 - g(2-\delta) (c-1)\right]
	\exp\left(g(2-\delta)(c-1)\right)\right\},
\end{IEEEeqnarray*}
as claimed.


\section{Proof of Theorem \ref{theorem:X_hat_symm}}\label{app:proof_hat_X_symm}

Starting from the smoothing lemma \eqref{eq:int_berry_esseen}
stating that $\Delta(X,\hat{X})\leq \bar{I}(X, \hat{X})$,
we will proceed to bound $\bar{I}(X, \hat{X})$.
Using the expressions \eqref{eq:true_stable_CF} and 
\eqref{eq:proxy_stable_cf} and performing
the change of variables \eqref{eq:change_var}, we can 
express,
\begin{IEEEeqnarray}{rCl}
	\bar{I}(X, \hat{X})
	&
	=& 
	\frac{1}{\pi}\int_{-\infty}^\infty \frac{\left|\phi_{{X_{(0,c)}}}(s)\right| \left|\phi_{R_{(c,\infty)}}(s) - \phi_{\hat{R}}(s) \right|}{|s|}             
	\dif s \nonumber
	\\
	& = &
	\frac{1}{\pi}\int_{0}^\infty \frac{\left|\omega_{{X_{(0,c)}}}(u)\right| \left|\omega_{R_{(c,\infty)}}(u) - \omega_{\hat{R}}(u) \right|}{u}             
	\dif u 
	\nonumber 
	\\
	& = &
	\frac{1}{\pi}\int_{0}^\infty \frac{\left|\omega_{{X_{(0,c)}}}(u)\right| \left|\psi_{Z_{(c,\infty)}}(u) - \psi_{Z}(u) \right|}{u}             
	\dif u,
	\label{eq:truncated_plus_residual}
\end{IEEEeqnarray}
where 
$\omega_{X_{(0,c)}}(u)$,
$\omega_{R_{(c, \infty)}}(u)$ 
and $\omega_{\hat{R}}(u)$ are
defined in \eqref{eq:cf_X0c_symm},  
\eqref{eq:phi_R_c_inf_symm} and \eqref{eq:residual_proxy_u}
respectively; the last equality follows again 
by  \eqref{eq:change_var}, with $\psi_{Z_{(c, \infty)}}(w)$ 
as in \eqref{eq:asympt_cf_res_symm} and $\psi_Z(w)$ 
as in \eqref{eq:cf_std_gauss_w}. 
The proof is in the following steps. 

The term $|\psi_{Z_{(c,\infty)}}(u) - \psi_{Z}(u)|$ in
\eqref{eq:truncated_plus_residual} has already been bounded over the intervals 
$[0,1]$ and $(1,\infty)$ in the proof of Theorem~\ref{thm:linear};
see Appendix~\ref{appendix:proofs_th_lin}. 
Combining equations~\eqref{eq:esseen2} and~\eqref{eq:expansion}, and recalling the definition of $K(a)$ in \eqref{eq:K_c_a}, we obtain,
\begin{IEEEeqnarray}{rCl} 
	\left|\psi_{Z_{(c,\infty)}}(u) - \psi_{Z}(u) \right|
	&\leq &{cu^2}\left(\frac{a}{2(2-a)}+\frac{1}{\eta^2}\right) \exp(-(c-1)h(u))\nonumber
	\\
	&= &{cu^2\pi}K(a) \exp(-(c-1)h(u)), \label{eq:bound_yiannis_wrt_u}
	\IEEEeqnarraynumspace
\end{IEEEeqnarray}	
with,
\begin{IEEEeqnarray*}{rCl} 
	h(u)&=&\begin{cases}-\barg u,&u\in[0,1]\\
		e^{-1}-1+u^a\bar{\gamma}(a),&u>1.\end{cases}
\end{IEEEeqnarray*}	

For the term $|\omega_{X_{(0,c)}}(u)|$ in \eqref{eq:truncated_plus_residual},
we recall the bounds obtained in Appendix \ref{subsection:q_properties} based 
on the fact that it is log-convex.
From \eqref{eq:bound_log_phi_X0c_piece_lin} and \eqref{eq:bound_yiannis_wrt_u},
the numerator of the integrand in \eqref{eq:truncated_plus_residual} 
is bounded as, 
\begin{IEEEeqnarray*}{rCl} 
	\left|\omega_{X_{(0,c)}}(u)\right|
	\left|\psi_{Z_{(c,\infty)}}(u) - \psi_{Z}(u) \right|
	&\leq &  cu^2 \pi K(a)\exp(\tilde h(u)),
\end{IEEEeqnarray*}	
with,
\begin{IEEEeqnarray*}{rCl} 
	\tilde  h(u)
	&=&
	\begin{cases}
		\sum_{i = 0}^{N-1}  \left( (m_i +(c-1)\barg) u + q_i\right) \mathds{1}_{A_i}{(u)}, & u\in[0,1]\\
		k_{(1,\infty)}-(c-1)(e^{-1}-1+u^a\bar{\gamma}(a)), & u>1
	\end{cases}
	\\
	&=&
	\begin{cases}
		\sum_{i = 0}^{N-1}  \left( \tilde{m}_i u + q_i\right) \mathds{1}_{A_i}{(u)}, & u\in[0,1]\\
		\tilde k_{(1,\infty)}-\tilde l_{(1,\infty)}u^a, & u>1,
	\end{cases}
\end{IEEEeqnarray*}	
where $m_i$ and $q_i$ are as in \eqref{eq:q_m};
$k_{(1,\infty)}$ is as in \eqref{eq:k_1_infty}, and,
\begin{IEEEeqnarray*}{rCl} 
	\tilde m_i &:=& m_i  +(c-1)\barg, \\ 
	\tilde k_{(1,\infty)} & = & k_{(1,\infty)} - (c-1)(e^{-1} -1), \\
	\tilde l_{(1,\infty)} & =&  (c-1)\bar{\gamma}(a).
\end{IEEEeqnarray*}	
Finally, substituting into the integral \eqref{eq:truncated_plus_residual}
\begin{IEEEeqnarray}{rCl} 
	\bar{I}(X, \hat{X}) &=&cK(a)\int_0^\infty u\exp(-\tilde h(u))\dif u
	\nonumber\\
	&=& c K(a)(I_N^{{X}}(c)+J^{{X}}(c)),
	\label{eq:twoIJ}
\end{IEEEeqnarray}	
where $I_N^{{X}}(c)$ and $J^{{X}}(c)$ denote the integrals over $(0,1)$ 
and $(1,\infty)$, respectively. Computing the integral $I^X_N(c)$,
\begin{IEEEeqnarray*}{rCl} 
	I_N^{{X}}(c)
	&=&
	\sum_{i=0}^{N-1}\int_{A_i} u\exp(\tilde m_i u + q_i)\dif u
	\\
	&=&
	\sum_{i=0}^{N-1}\exp(q_i)\int_{A_i} u\exp(\tilde m_i u)\dif u
	\\
	&=&
	\sum_{i=0}^{N-1}\frac{\exp(q_i)}{\tilde m_i}\left[ \exp(\tilde m_i u )\left( u - \frac{1}{\tilde m_i}\right) \right]_{u_i}^{u_{i+1}} 
	\\
	&=&
	\sum_{i=0}^{N-1}\frac{e^{q_i}}{\tilde m_i}\left[ e^{\tilde m_i u_{i+1} }\left( u_{i+1} - \frac{1}{\tilde m_i}\right)  
	-
	e^{\tilde m_i u_{i} }\left( u_{i} - \frac{1}{\tilde m_i}\right)
	\right].
\end{IEEEeqnarray*}
Observe that, when $N=1$, $i =0$, $q_0 = 0$, $u_0 = 0$, $u_1 = 1$, $\tilde m_0 = m_0  +(c-1)\barg = \log(\omega_{X_{(0,c)}}(1)) + (c-1)\barg$, and the last equation becomes,
\begin{IEEEeqnarray*}{rCl} 
	I_1^{{X}}(c) 
	&=& 
	\frac{1}{\tilde m_0}\left[ e^{\tilde m_0  }\left( 1 - \frac{1}{\tilde m_0}\right)  
	+ 
	\frac{1}{\tilde m_0}
	\right]
	\\
	&=& 
	\frac{1}{\tilde m_0}
	\left[ e^{\tilde m_0  } + \frac{1}{\tilde m_0} \left( 1 - e^{\tilde m_0  } \right)  
	\right].
\end{IEEEeqnarray*}	
Similarly for the integral $J^X(c)$, 
\begin{IEEEeqnarray*}{rCl} 
	J^{{X}}(c)
	&=&
	\int_1^\infty u\exp\big(\tilde k_{(1,\infty)}-\tilde l_{(1,\infty)}u^a\big) \dif u
	\\
	&=&
	\exp{\big(\tilde k_{(1,\infty)}\big)}\int_1^\infty u\exp\big(-\tilde l_{(1,\infty)}u^a\big) \dif u
	\\
	&=&
	\frac{\exp{\big(\tilde k_{(1,\infty)}\big)}}
	{a \big(\tilde l_{(1,\infty)}\big)^{2/a}}\int_{\tilde l_{(1,\infty)}}^\infty t^{\frac{2}{a} -1} \exp(-t) \dif t
	\\
	&=&
	\frac{\exp{\big(\tilde k_{(1,\infty)}\big)}}
	{a \big(\tilde l_{(1,\infty)}\big)^{2/a}} \Gamma\big({2}/{a},\tilde l_{(1,\infty)} \big). 
\end{IEEEeqnarray*}	
Substituting these in (\ref{eq:twoIJ}) yields exactly the claimed bound.
\hfill{\footnotesize $\blacksquare$}

\medskip

\begin{figure*}[ht!] 
	\centerline{
		\includegraphics[width=0.3\linewidth]{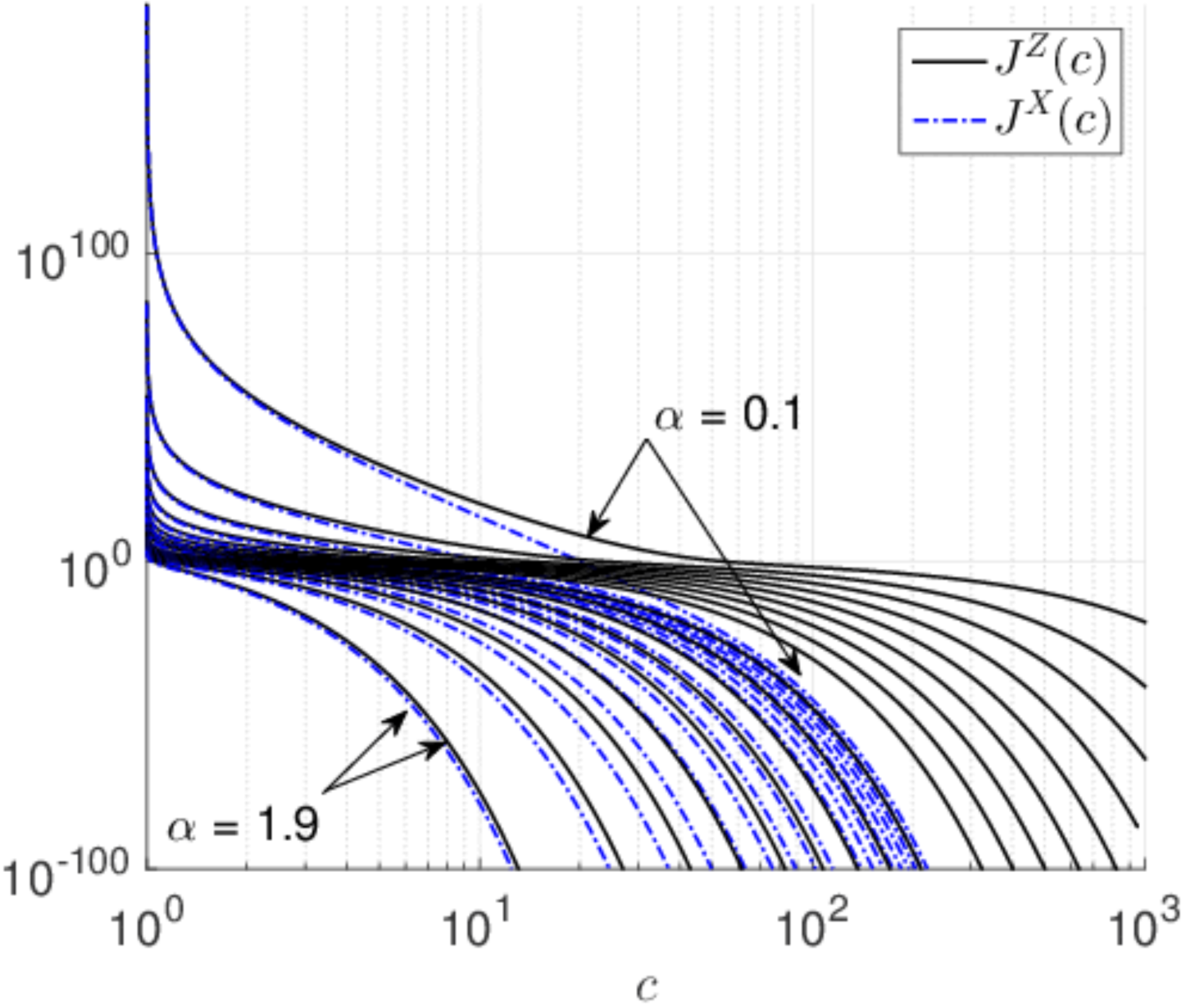} 
		\includegraphics[width=0.3\linewidth]{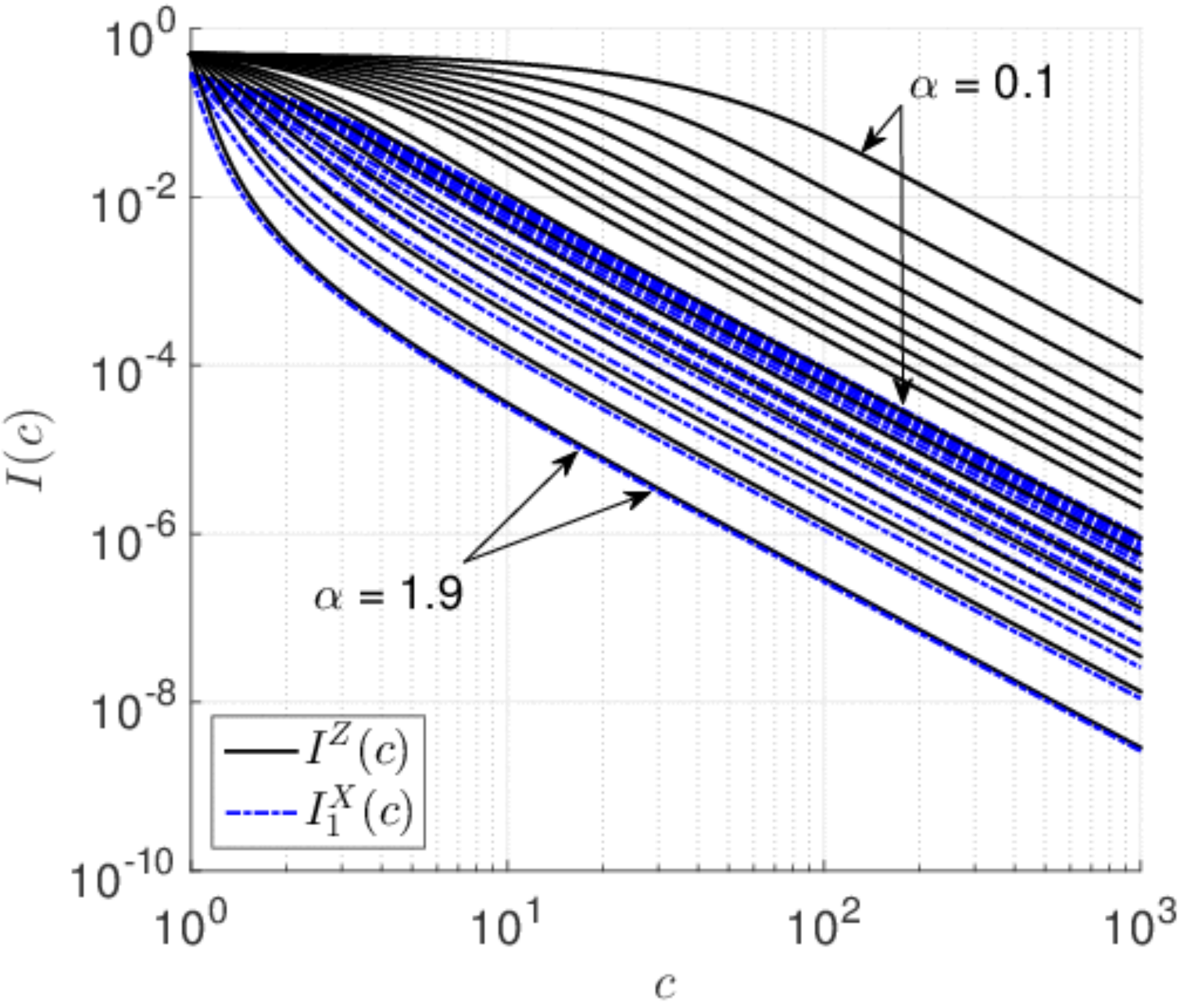}	
		\includegraphics[width=0.3\linewidth]{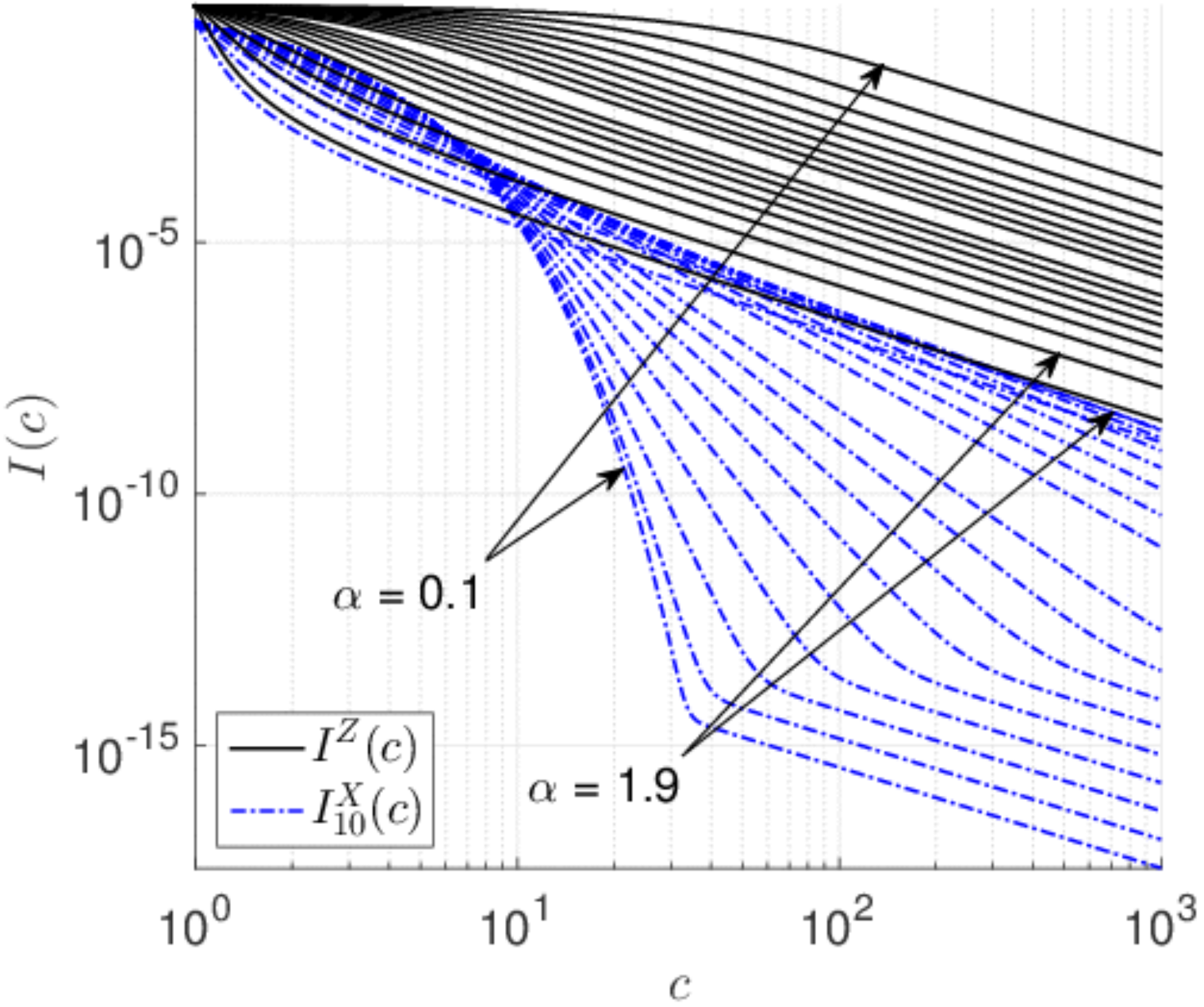}		
	}		
	\caption{Comparison of the terms $I^Z(c)$ and $I_N^{{X}}(c)$, and $J^Z(c)$ and $J^{{X}}(c)$, for $N = 1,10$, and $\alpha = 0.1, 0.2, \ldots, 1.9$, plotted against $c \in [1, 10^3]$. }	
	\label{fig:comparison_I_J_X_0c}
\end{figure*}

\begin{figure*}[ht!] 
	\centerline{
		\includegraphics[width=0.3\linewidth]{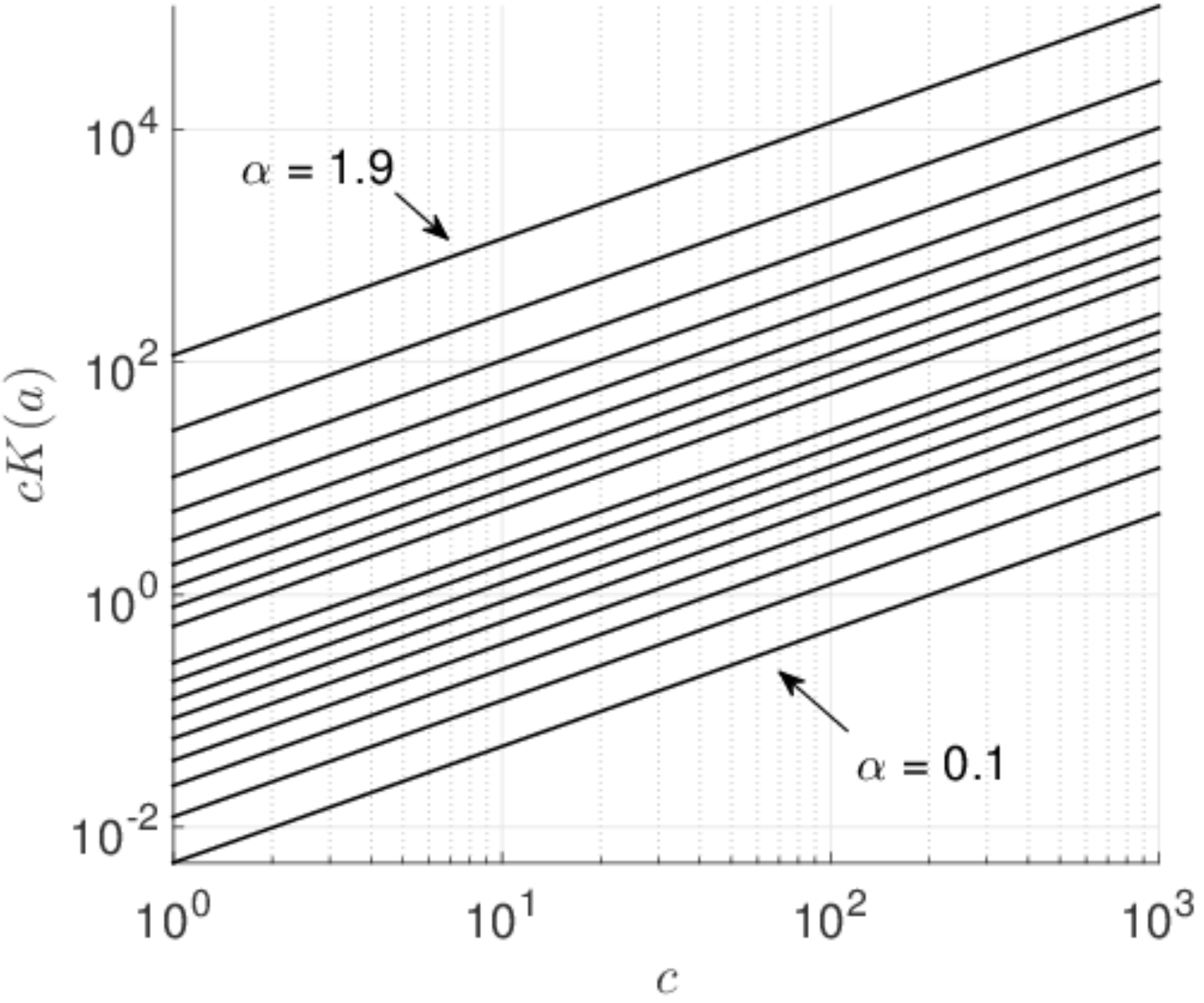}	
		\includegraphics[width=0.3\linewidth]{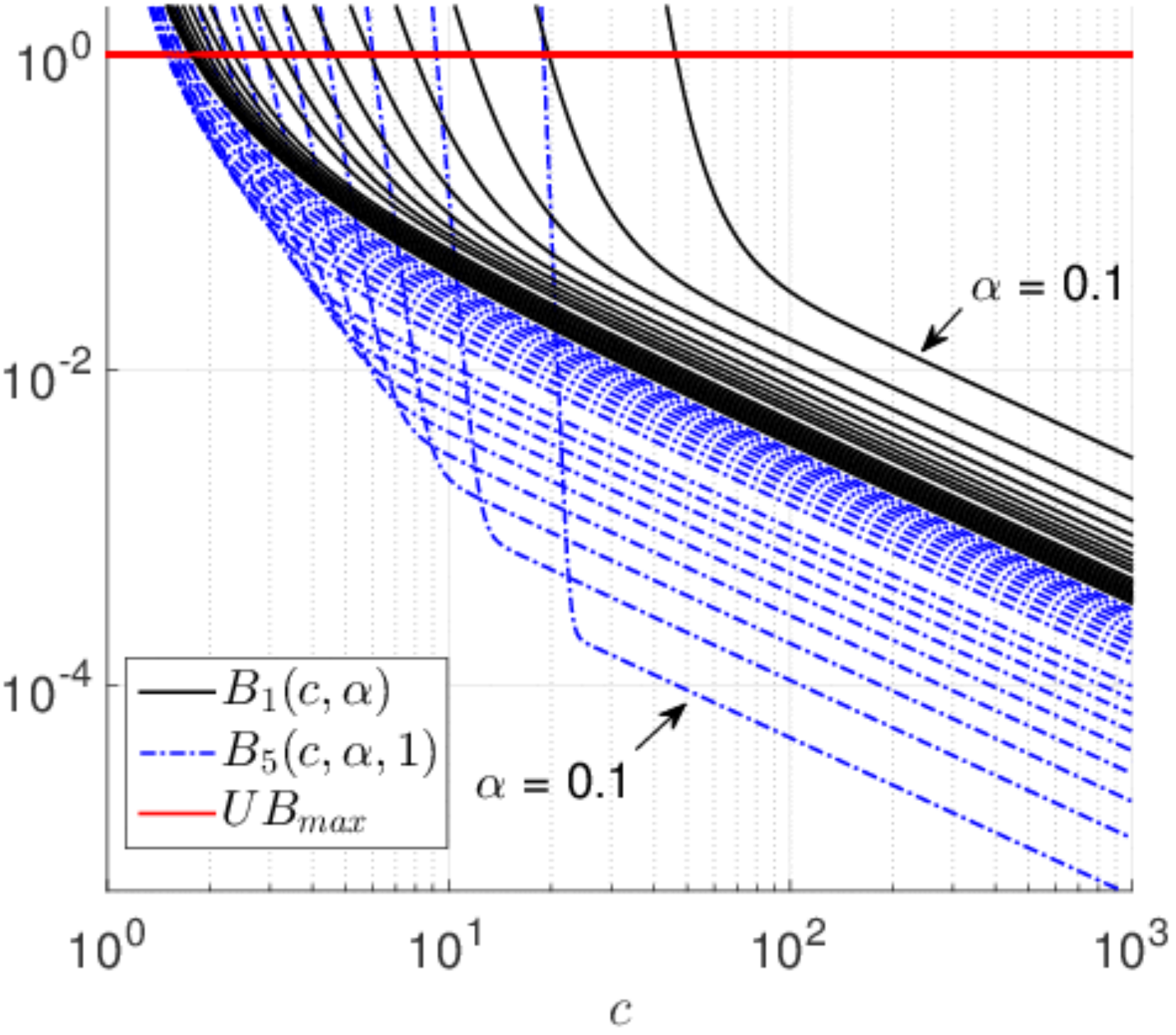}	
		\includegraphics[width=0.3\linewidth]{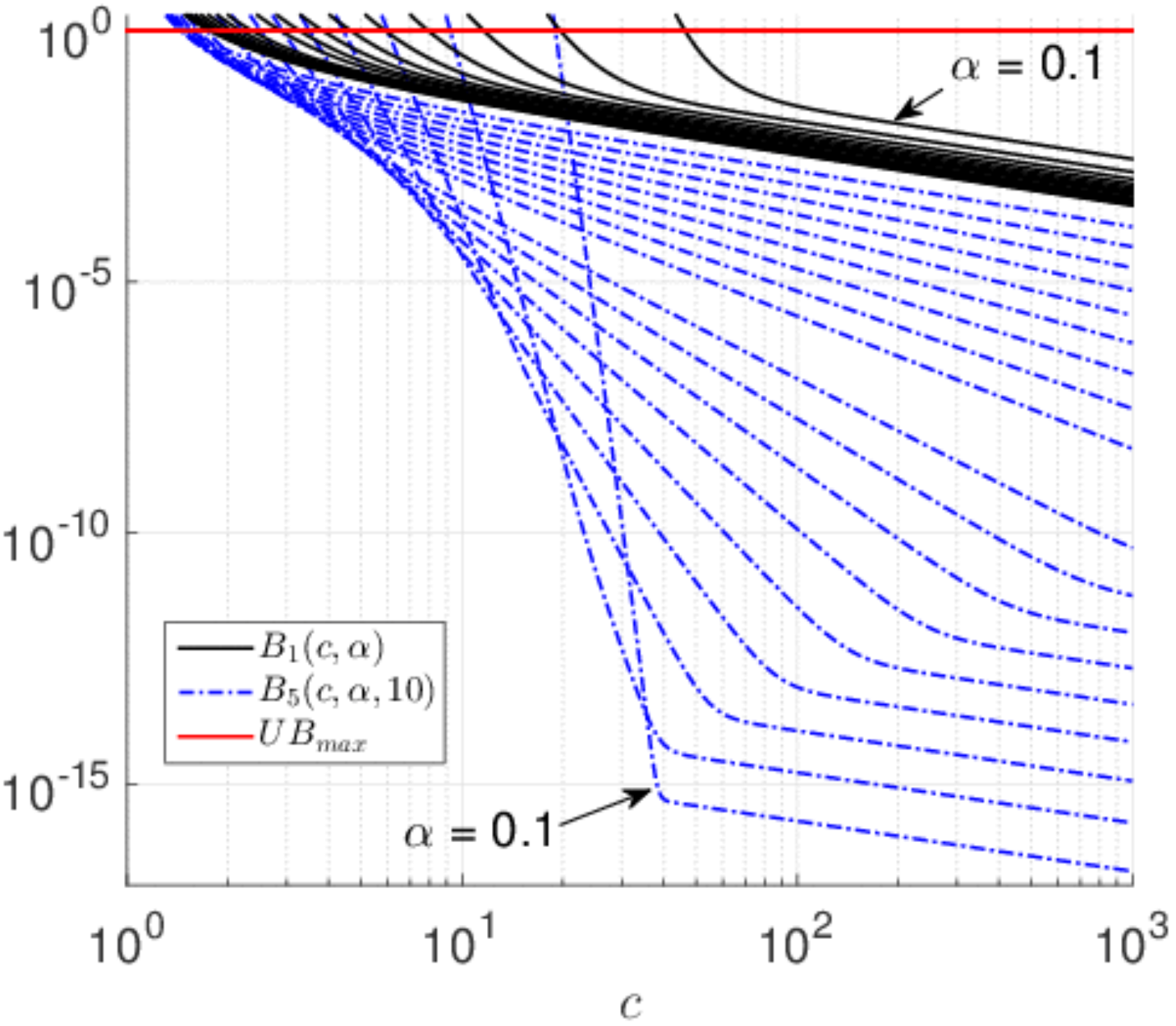}	
	}
	\caption{Left: The function $cK(a)$, with $a = \alpha/2$, 
plotted against $1 \leq c \leq 10^3$. Centre and right: comparison of $B_4(c, \alpha) = cK(a)(I^Z(c) + J^Z(c))$ and $B_5(c, \alpha, N) =cK(a)(I_N^{{X}}(c) + J^{{X}}(c))$, for $N = 1,10$, and $\alpha = 0.1, 0.2, \ldots, 1.9$, plotted against $1 \leq c \leq 10^3$. 
[The red horizontal line at 1 shows the 
maximum possible value of the Kolmogorov distance.]}	
	\label{fig:comparison_B1_B5}
\end{figure*}

\begin{remark}\label{remark:plots_appendix}
Observe that the bounds $I_N^{{X}}(c)$ and $J^{{X}}(c)$ are smaller 
than $I^Z(c)$ and $J^Z(c)$, because the latter correspond to the former when
$\omega_{X_{(0,c)}}(u) \equiv 1$;
Figure \ref{fig:comparison_I_J_X_0c} illustrates their difference.
The asymptotic rates of both $J^Z(c)$ and $J^{{X}}(c)$ 
depend on $\alpha$, and that larger values of $\alpha$ give smaller bounds. 
On the other hand, the asymptotic rates of both $I^Z(c)$ and 
$I_N^{X}(c)$ are independent on $\alpha$, but their nonasymptotic 
behaviour does depend of $\alpha$ for a wide range of values of~$c$.
Moreover, for $N>1$, lower values of $\alpha$ now give lower values of 
the bound $I_N^{{X}}(c)$,
whereas for $N=1$, the dependence of $I_N^{{X}}(c)$ 
on $\alpha$ is the same as that of $I^Z(c)$. But even in this case, 
the dependence of the combined bound $B_5(c, \alpha, N)$ on $\alpha$ is the opposite 
to that of $B_1(c, \alpha)$, as illustrated in Figure~\ref{fig:comparison_B1_B5}.
{ This can be justified by the fact that the growth of $I^X$ in $\alpha$ is much slower that that of $I^Z$, so while 
 the coefficient $cK(a)$ rectifies the former, it fails to do so with the latter. 
}

\end{remark}

\section{Proof of Proposition \ref{prop:X_0c_symm_conv}}\label{appendix:proof:X0c}

We adapt the proof strategy used for a similar result
in \cite{LedouxPaulauskas1996} on the convergence of the 
truncated PSR with 
$W_1\sim\mathcal{N}(0,1)$. The difference here is that
the number of terms in the truncated PSR is random, 
and that we allow the variance $\sigma_W^2$ to not
necessarily be equal to~1.

For $x>0$,
\begin{IEEEeqnarray*}{rCl} 
	\Delta(X, X_{(0,c)}) & =  &\sup_{x \in \mathbb{R}}\left| F_X(x) - F_{X_{(0,c)}}(x)\right|
	\\
	& = &
	\sup_{x \in \mathbb{R}}	
	\left|
	\mathbb{E}[\mathds{1}_{(X\leq x)}]  - \mathbb{E}[\mathds{1}_{(X_{(0,c)}\leq x)}]
	\right|. 
\end{IEEEeqnarray*}
\noindent 
Conditioning first on the number $N_{(0,c)}$ of terms in
the truncated PSR and the Poisson event times $\{\Gamma_j\}$,
we can expand,
\begin{IEEEeqnarray*}{rCl} 
	\mathbb{E}[\mathds{1}_{(X\leq x)}] 	
	&=& \mathbb{E}\left[\mathds{1}_{ \big(\sum_{j =1}^\infty \Gamma_j^{-1/\alpha} W_j \leq x \big) } \right]  
	\\
	&=& \mathbb{E}\left[\mathbb{E}\left[\mathds{1}_{ \big(\sum_{j =1}^\infty \Gamma_j^{-1/\alpha} W_j \leq x \big) } \middle|  N_{(0,c)}, \{\Gamma_j\}_{j=1}^\infty\right]  \right]
	\\
	&=& \mathbb{E}\left[\mathbb{P}\left(\sum_{j =1}^\infty \Gamma_j^{-1/\alpha} W_j \leq x   \middle| N_{(0,c)}, \{\Gamma_j\}_{j=1}^\infty \right)  \right]
	\\
	&=& \mathbb{E}\left[ F_Z\left( \frac{x}{\sigma_W S} \right) \right], 
\end{IEEEeqnarray*}
\noindent
where $F_Z(\cdot)$ is again the CDF of the standard normal distribution 
and $S^2$ is defined in \eqref{eq:inf_variance_X}. 
Similarly we can compute,
\begin{IEEEeqnarray*}{rCl} 
	\mathbb{E}[\mathds{1}_{(X_{(0,c)}\leq x)}] 	&=&  \mathbb{E}\left[ F_Z\left( \frac{x}{\sigma_W S_{N_{(0,c)}}} \right) \right], 
\end{IEEEeqnarray*}
where,
\begin{IEEEeqnarray*}{rCl} 
	S^2_{N_{(0,c)}} := \sum_{j=1}^{N_{(0,c)}} \Gamma_j^{-2/\alpha},
	\IEEEeqnarraynumspace
\end{IEEEeqnarray*}
\noindent 
with the convention that $S^2_{N_{(0,c)}} = 0$  if $N_{(0,c)} =0$, and $F_Z(-\infty) = 0,\, F_Z(\infty) = 1$. Then, given that $\Gamma_j>0$ for any $j\geq1$ with probability 1, it follows that $ S^2_{N_{(0,c)}} < S^2$, and, for $x>0$,
\begin{IEEEeqnarray*}{rCl} 
	\Delta(X, X_{(0,c)}) & =  &\sup_{x \in \mathbb{R}}\mathbb{E}\left[F_Z\left( \frac{x}{\sigma_W S_{N_{(0,c)}}} \right) - F_Z\left( \frac{x}{\sigma_W S} \right) \right]. 
\end{IEEEeqnarray*}
The argument of the expectation above can be bounded as,
\begin{IEEEeqnarray*}{rCl} 
	F_Z\left( \frac{x}{\sigma_W S_{N_{(0,c)}}} \right) - F_Z\left( \frac{x}{\sigma_W S} \right) 	
	&=& \frac{1}{\sqrt{(2\pi)}}
	\int^{ x/(\sigma_W S_{N_{(0,c)}})}_{ x/(\sigma_W S)} \exp\left(- {u^2}/{2} \right) \dif u
	\\
	&\leq& 
	\frac{x}{\sigma_W\sqrt{(2\pi)}} \left(\frac{1}{S_{N_{(0,c)}}} - \frac{1}{S} \right)\exp\left(- \frac{x^2}{2\sigma^2_W S^2} \right)
	\\
	&\leq& 
	\frac{1}{\sigma_W\sqrt{(2\pi)}} \left(\frac{S - S_{N_{(0,c)}}}{ S_{N_{(0,c)}} S} \right)S\sigma_W \exp(-1/2)
	\\
	&=& 
	\frac{\exp(-1/2)}{\sqrt{(2\pi)}} \left(\frac{S^2 - S_{N_{(0,c)}}S }{ S_{N_{(0,c)}} S} \right)
	\\
	&\leq& 
	\frac{\exp(-1/2)}{\sqrt{(2\pi)}} \left(\frac{S^2 - S^2_{N_{(0,c)}}}{ S^2_{N_{(0,c)}} } \right),
\end{IEEEeqnarray*}
where in the first inequality we used the constant bound of the integrand 
on the integration domain; in the second inequality we used the fact that the 
mode of the Rayleigh distribution is achieved for $x = \sigma_WS$; and
in the third inequality we used (twice) the fact that $ S^2_{N_{(0,c)}} < S^2$. Notice that this upper bound does not depend on $\sigma_W$. Finally, taking the expectation, 
\begin{IEEEeqnarray*}{rCl} 
	\Delta(X, X_{(0,c)}) 	
	&\leq &  	
	\frac{\exp(-1/2)}{\sqrt{(2\pi)}} \mathbb{E}\left[{ S^{-2}_{N_{(0,c)}} } (S^2 - S^2_{N_{(0,c)}}) \right]
	\\
	&\leq &  	
	\frac{\exp(-1/2)}{\sqrt{(2\pi)}} \left(\mathbb{E}\left[{ S^{-4}_{N_{(0,c)}} }\right]\right)^{1/2}
	\left(\mathbb{E}\left[ (S^2 - S^2_{N_{(0,c)}})^2 \right]\right)^{1/2}
	\\
	&\leq &  	
	\frac{\exp(-1/2)}{\sqrt{(2\pi)}} \left(\mathbb{E}\left[\Gamma_1^{4/\alpha}\right]\right)^{1/2}
	\left(\mathbb{E}\left[ (S^2 - S^2_{N_{(0,c)}})^2 \right]\right)^{1/2}
	\\
	& = &  	
	\frac{\exp(-1/2)}{\sqrt{(2\pi)}} \left(\Gamma \left(\frac{\alpha+4}{\alpha}\right)\right)^{1/2}
	\times
	\left(\frac{\alpha}{4-\alpha} c^{\frac{\alpha -4}{\alpha}} +
	\left(\frac{\alpha}{2-\alpha}c^{\frac{\alpha-2}{\alpha}}\right)^2\right)^{1/2},
\end{IEEEeqnarray*}
\noindent 
where in the second inequality we used the Cauchy-Schwartz inequality; 
and in the third inequality we used the fact that 
$S^2_{N_{(0,c)}} > \Gamma_1^{-2/\alpha}$, where $\Gamma_1$ is the smallest 
of the $\Gamma_j$ variables, hence exponentially distributed. Thus,
\begin{IEEEeqnarray*}{rCl} 
	\mathbb{E}\left[\Gamma_1^{4/\alpha}\right] = \int_0^\infty u^{(4/\alpha + 1) -1} e^{-u} \dif{u} = \Gamma\left(\frac{4}{\alpha} +1\right),	
	\IEEEeqnarraynumspace
\end{IEEEeqnarray*}
with $\Gamma(\cdot)$ the gamma function \eqref{eq:gamma}.
On the other hand, 
$ S^2 - S^2_{N_{(0,c)}} =  \sum_{j = N_{(0,c)} +1}^\infty\Gamma_j^{-2/\alpha} = \lim_{d\rightarrow \infty }\sum_{j : \Gamma_j \in (c,d)}\Gamma_j^{-2/\alpha}$, and we know that such limit exists through the PSR, being the second moment of the PSR residual with deterministic $W_j =1$. Hence,
\begin{IEEEeqnarray*}{rCl} 
	\mathbb{E}\left[ (S^2 - S^2_{N_{(0,c)}})^2 \right] 
	&=& \lim_{d\rightarrow \infty}\mathbb{E}\left[ \sum_{j : \Gamma_j \in (c,d)}\Gamma_j^{-2/\alpha} \sum_{i : \Gamma_i \in (c,d)}\Gamma_i^{-2/\alpha}\right] \nonumber
	\\
	&=& \lim_{d\rightarrow \infty}\mathbb{E}\left[\mathbb{E}\left[ \sum_{j=1}^{N_{(c,d)}}\Gamma_j^{-2/\alpha} \sum_{i=1}^{N_{(c,d)}}\Gamma_i^{-2/\alpha}\middle| N_{(c,d)}\right]\right]  \nonumber
	\\
	&=& \lim_{d\rightarrow \infty}\mathbb{E}\left[\mathbb{E}\left[ \sum_{j=1}^{N_{(c,d)}}U_j^{-2/\alpha} \sum_{i=1}^{N_{(c,d)}}U_i^{-2/\alpha}\middle| N_{(c,d)}\right]\right]  \nonumber
	\\
	&=& \lim_{d\rightarrow \infty}
	\mathbb{E}[N_{(c,d)}]\mathbb{E}[U_1^{-4/\alpha}]
	+ 
	\mathbb{E}[N^2_{(c,d)} - N_{(c,d)}](\mathbb{E}[U_1^{-2/\alpha}])^2 \nonumber
	\\
	&=& \frac{\alpha}{4-\alpha} c^{\frac{\alpha -4}{\alpha}} + \left(\frac{\alpha}{2-\alpha}c^{\frac{\alpha-2}{\alpha}}\right)^2,
\end{IEEEeqnarray*}

\noindent
where the $\{U_j\}$ are i.i.d.\ uniformly distributed RVs on $(c,d)$, 
as in \eqref{eq:U_j_c_d}. 
\hfill{\footnotesize $\blacksquare$}


\section{Proof of Proposition \ref{prop:res_contribution}} \label{appendix:proof_contribution}

Recall that 
$\bar{I}(X, \hat{X})$ is given by 
\eqref{eq:truncated_plus_residual}.
For 
$\bar{I}(X, X_{(0,c)})$,
using \eqref{eq:true_stable_CF} and
\eqref{eq:proxy_stable_cf}, and performing
the change of variables \eqref{eq:change_var}, we 
similarly have,
\begin{eqnarray}
	\bar{I}(X, X_{(0,c)})
	= 
	\frac{1}{\pi}\int_{0}^\infty \frac{\left|\omega_{{X_{(0,c)}}}(u)\right| \left|\psi_{Z_{(c,\infty)}}(u) - 1 \right|}{u}             
	\dif u.
	\label{eq:not_truncated}
\end{eqnarray}

We proceed by comparing the integrands in 
\eqref{eq:truncated_plus_residual}.
	and (\ref{eq:not_truncated}).
Write $\tilde{c}(a)= c(\alpha)$, with $a=\alpha/2$.
It suffices to show that, for all $c>\tilde{c}(a)$
and all $u>0$,
\begin{IEEEeqnarray}{rCl} 
	\left|\psi_{Z_{(c,\infty)}}(u) - \psi_{Z}(u) \right| 	< \left|\psi_{Z_{(c,\infty)}}(u) - 1 \right|,
	\label{eq:c_alpha_inequality} 
	\IEEEeqnarraynumspace
\end{IEEEeqnarray}
where $\log(\psi_{Z}(u)) =-cu/\eta$,
and
$\psi_{Z_{(c, \infty)}}$, defined in \eqref{eq:asympt_cf_res_symm},
satisfies,
\begin{IEEEeqnarray*}{rCl} 
	\log(\psi_{Z_{(c, \infty)}}(u))\nonumber
	&=& c(1-\exp(-u)-u^{a}\gamma(1-a,u))
	\\
	& = & cg(u),
\end{IEEEeqnarray*}
with $g$ defined in \eqref{eq:gdefn}.
Using the fact that $g(u) <0$, see 
\eqref{eq:gUB1}, we have $|\omega_{R_{(c,\infty)}}(u) - 1| = 1 - \omega_{R_{(c,\infty)}}(u)$. Furthermore, using \eqref{eq:gLB}, we have 
$|\omega_{R_{(c,\infty)}}(u) - \omega_{\hat{R}}(u) | = \omega_{R_{(c,\infty)}}(u) - \omega_{\hat{R}}(u)$, and \eqref{eq:c_alpha_inequality} becomes,    
\begin{IEEEeqnarray*}{rCl} 
	v(u) := 1- 2\exp(c g(u)) +  \exp(-cu/\eta) > 0.
\end{IEEEeqnarray*}
Using \eqref{eq:gUB}, we have that, for $u \in (0,1]$
\begin{IEEEeqnarray*}{rCl} 
	v(u)
	&\geq&
	(1 - 2\exp(-cu/2\eta) + \exp(-cu/\eta))
	\\
	&=&(1 - \exp(-cu/2\eta))^2 > 0, 
\end{IEEEeqnarray*}
hence \eqref{eq:c_alpha_inequality} holds for any $c>0$ when $u\leq 1$. 
We then consider the case $u>1$. By the monotonicity of $g$, 
see \eqref{eq:gprime}, we have that $g(u) < g(1)$ for $u>1$, leading to,
\begin{IEEEeqnarray*}{rCl} 
	v(u) > 1 -2\exp(cg(1)) + \exp(-uc/\eta) := z(u), \quad u>1,
\end{IEEEeqnarray*}
where $z(u)$ is a lower bound on $v(u)$ for $u>1$. 
We know that $z(1) = v(1) >0$ as shown above,
and also,
\begin{IEEEeqnarray*}{rCl} 
	z'(u)
	&=&
	-c/\eta \exp(-uc/\eta) < 0, 
\end{IEEEeqnarray*}
implying that the lower bound $z(u)$ on $v(u)$ is decaying. 
Furthermore,
	$\lim_{u \rightarrow \infty} z(u) = 1 - 2\exp(cg(1))$,
so that, for all $u\geq 1$,
$$v(u)\geq
	1 - 2\exp(cg(1)),$$
and the right-hand side above is itself positive
as long as
$c>-{\log(2)}/{g(1)} = \tilde c(a)$, as required.




\ifCLASSOPTIONcaptionsoff
  \newpage
\fi




\newpage

\bibliographystyle{plain}
\bibliography{alpha_Stable_CLT_arxiv_v4}
\end{document}